\documentclass[10pt]{amsart}
\usepackage{amsmath,amssymb,amsfonts,amscd}
\input xy
\xyoption{all}

\newtheorem{theorem}{Theorem}[subsection]
\newtheorem{lemma}[theorem]{Lemma}
\newtheorem{proposition}[theorem]{Proposition}
\newtheorem{definition}[theorem]{Definition}
\newtheorem{corollary}[theorem]{Corollary}

\newtheorem{example}[theorem]{Example}

\DeclareMathOperator{\M}{\mathcal{M}}
\DeclareMathOperator{\LM}{\mathcal{LM}}
\DeclareMathOperator{\RM}{\mathcal{RM}}
\DeclareMathOperator{\tr}{Tr}
\DeclareMathOperator{\rM}{rM}

\DeclareMathOperator{\lL}{\ell L}
\DeclareMathOperator{\ran}{ran}
\DeclareMathOperator{\MAX}{MAX}
\DeclareMathOperator{\MIN}{MIN}
\DeclareMathOperator{\itimes}{\check{\otimes}}
\DeclareMathOperator{\ptimes}{\hat{\otimes}}
\DeclareMathOperator{\htimes}{\otimes_h}
\DeclareMathOperator{\stimes}{\overline{\otimes}}

\DeclareMathOperator{\Id}{Id}
\DeclareMathOperator{\spn}{span}

\DeclareMathOperator{\wks}{wk*}
\DeclareMathOperator{\wk}{wk}
\DeclareMathOperator{\Ml}{{\mathcal M}_\ell}
\DeclareMathOperator{\Al}{{\mathcal A}_\ell}
\DeclareMathOperator{\Mr}{{\mathcal M}_r}
\DeclareMathOperator{\Ar}{{\mathcal A}_r}
\DeclareMathOperator{\Cl}{\mathcal{C}_\ell}
\DeclareMathOperator{\Cr}{\mathcal{C}_r}
\DeclareMathOperator{\Her}{Her}
\DeclareMathOperator{\wlim}{wk*-\lim}

\DeclareMathOperator{\I}{I}
\DeclareMathOperator{\II}{II}
\DeclareMathOperator{\III}{III}
\DeclareMathOperator{\op}{op}
\DeclareMathOperator{\eps}{\varepsilon}
\DeclareMathOperator{\ball}{ball}
\DeclareMathOperator{\Cent}{Cent}

\begin{document}
\title[The Calculus of One-Sided $M$-Ideals and Multipliers]{The 
Calculus of One-Sided $M$-Ideals and Multipliers \\
in Operator Spaces}
\author{David P. Blecher}
\address{Dept. of Mathematics, Univ. of Houston, Houston, TX 77204}
\email{dblecher@math.uh.edu}
\author{Vrej Zarikian}
\address{Dept. of Mathematics, Univ. of Texas, Austin, TX 78712}
\email{zarikian@math.utexas.edu}
\begin{abstract}
The theory of
one-sided $M$-ideals and multipliers of operator spaces
is simultaneously a generalization of classical $M$-ideals,
ideals in operator algebras, and aspects of the theory
of Hilbert $C^*$-modules and their maps.
Here we give a systematic exposition of this theory;
a reference tool for `noncommutative
functional analysts' who may encounter
a one-sided $M$-ideal or multiplier in their
work.
\end{abstract}
\thanks{Blecher was supported by a grant from the NSF.
Zarikian was supported by a NSF VIGRE Postdoctoral Fellowship.}
\date{\today}
\maketitle

\newpage

\tableofcontents

\newpage

\section{Introduction} \label{i}

The classical theory of $M$-ideals emerged in the early 
seventies with the paper
\cite{AlfsenEffros} of Alfsen and Effros.
Recently in \cite{BEZ} (see also \cite{Zarikian})
the authors together with E. G. Effros developed 
a one-sided variant of this classical theory.
The intention was to create a tool that 
would be useful for, and appropriate to,
`noncommutative functional analysis'.
Our one-sided theory contains the classical $M$-ideals
and summands as particular examples.
Other interesting examples of one-sided $M$-ideals
include right ideals in $C^*$-algebras, submodules of Hilbert
$C^*$-modules, and 
an important class
of right ideals  in nonselfadjoint operator algebras.
(see also \cite{BSZ}). A key point is that
our definitions are stated only in terms
of the underlying operator space structure, and yet very
often encode important algebraic
information. This is already seen in the classical
case: for example, the classical $M$-ideals in a $C^*$-algebra
are exactly the closed two-sided ideals.   
In the present article we give a systematic
and up-to-date
account of the `one-sided $M$-theory'.  
Following a program proposed by 
Effros in 2000, we exhibit here the `calculus of one-sided $M$-ideals'.

One major device used in the analysis of one-sided $M$-ideals
is the recent theory
of one-sided operator space multipliers introduced in \cite{Shilov}, 
which we also take a little further in this article.  
The `left $M$-projections' are exactly the orthogonal 
projections in a certain algebra $\Al(X)$ of 
such multipliers.
In fact, $\Al(X)$ is a $C^*$-algebra, and moreover is
a von Neumann algebra if $X$ is a dual operator space.  
This is a key point in our theory, since amongst other
things, it allows us to manipulate left $M$-projections
in the way that operator theorists and 
$C^*$-algebraists are used to.   Applying
basic von Neumann algebra theory leads to our most interesting
results.

Many of our results consist of `one-sided counterparts' to
the classical calculus of $M$-ideals 
of Banach spaces, and to
the related theory of multipliers 
and centralizers of such spaces (which the reader will find comprehensively 
treated in \cite{Behrends,HWW}).  The proofs, however, are 
far from being simple
imitations of the classical ones, and require quite different,
`noncommutative', arguments.   
Other results which we obtain 
appear to be new even for classical $M$-ideals, or at any rate do not
seem to be in the literature.  We exhibit many results which are
truly `noncommutative' in nature, having no classical
counterpart (e.g.\ the existence of one-sided $M$-ideals in
certain noncommutative $L^p$ spaces).             
We give a large assortment of examples
coming from rather diverse sources (e.g.\
Hilbertian operator spaces, Hilbert $C^*$-modules, 
nonselfadjoint operator algebras, and
low finite-dimensional operator spaces).
These examples show that the noncommutative versions
of approximately half of the 
classical `calculus of $M$-ideals' break down without further 
hypotheses. For example, the set of 
one-sided $M$-ideals in a given  operator space is 
closed with respect to
the `closed span' operation, $\vee$, but is 
not closed with respect to intersection, $\wedge$.
In fact the reason for this
is truly noncommutative, and is related to 
Akemann's result that  the `meet' of two open projections need not be open
 \cite{Akemann69}.  The intersection
of two one-sided $M$-ideals {\em is} a one-sided  $M$-ideal
if one imposes certain conditions considered by Akemann.
 
We now summarize very briefly 
the structure of our article.  
We have intentionally written
a reference tool for a reader who might encounter
a one-sided multiplier or $M$-ideal in their work.
In Section 2, we give definitions and basic results that 
will be used throughout. In Section 3, we study how
projections and partial isometries in the abstract
multiplier algebras manipulate the underlying space $X$.
The basic von Neumann algebra projection calculus
is key here.  
Section 4 is devoted to examples. The lengthy
Section 5 is the apex of our work, and is intended as 
a systematic presentation of the 
basic calculus of one-sided $M$-ideals and multipliers. 
We discuss subspaces, quotients, duality, interpolation, 
tensor products, etc..
Section 6 treats type decomposition
and Morita
equivalence for dual operator spaces. The short 
Section 7 is devoted to the complete $M$-ideals of 
Effros and Ruan (see \cite{ERcmp}), and the associated `centralizers'.
This is very similar to the classical theory, however
we will try to deduce the main results quickly from
results in earlier sections.  In Section 8 we discuss
a few directions for future progress.
In the Appendices we list some background facts
which are used very often,
about Banach spaces and infinite matrices over
operator spaces.

We end this introduction with a little notation. Other notation
will be encountered as we proceed.  
The underlying field is always $\mathbb{C}$. We usually use the 
letters $H, K, L, \cdots$ for Hilbert spaces.  
We write
$\overline{S}$ for the {\em norm closure} of a set $S$, and
$\overline{S}^{\wks}$ for the {\em  weak-$*$ closure}.   
Notationally, a problem arises in 
the use of the symbol $*$, which we use
for three different things in this paper.  Namely,
we have the dual space $X^*$ of a space $X$, the
adjoint or involution $S^*$ of an operator on a Hilbert space, and
the adjoint operator $R^* : Y^* \rightarrow Z^*$ of
an operator $R : Z \rightarrow Y$. We are forced by reasons
of personal
taste to leave it to the reader to determine which is
meant in any given formula. To alleviate some of the 
pressure, we use the symbol $T^\star$ for the involution
in the aforementioned $C^*$-algebra $\Al(X)$.    
We will use the word {\em projection} both for 
an idempotent linear map, and for an orthogonal 
projection in an operator algebra.

An operator
space is a linear subspace of $B(H)$ for a Hilbert space $H$.
In this paper all operator spaces are norm complete.
Equivalently, there is an abstract characterization due to
Ruan. Basic notation and facts about
operator spaces may be found in \cite{ERbook,Paulsen,Pisier}, for example
(there is also a forthcoming text on related matters
by the first author and Le Merdy).
We recall that if $X$ is
an operator space, then so is $X^*$; its matrix norms
come from the identification $M_n(X^*) \cong CB(X,M_n)$.
A {\em dual operator space} is one which is completely
isometrically isomorphic to the operator space dual of
another operator space. Any $\sigma$-weakly closed
subspace of $B(H)$ is a dual operator space with
predual determined by the predual of $B(H)$.
Conversely, any dual operator space
is linearly completely
isometrically isomorphic and weak-* homeomorphic to a
$\sigma$-weakly closed subspace of some $B(H)$.   
We also recall from basic operator space theory  
that for a dual operator space $X = Y^*$, the
space $CB(X)$ is canonically a dual operator space too.
A bounded net $T_i \in CB(X)$
converges in the weak-* topology to $T \in CB(X)$ if and
only if $T_i(x) \rightarrow T(x)$ weak-* in $X$, for every
$x \in X$.  We use  $H_c$ and $H_r$ for the
Hilbert column and row spaces associated to a Hilbert space $H$ (see the texts cited above
or Appendix \ref{App.B} below).

For the classical theory of $M$-ideals, the standard source 
is \cite{HWW}.
It would be helpful if the reader was at least vaguely familiar
with \cite{BEZ}, and had access to \cite{Shilov,BSZ}.      

We will use
the term $W^*$-algebra for a $C^*$-algebra with predual.
In view of a well-known theorem of Sakai \cite{Sakai} this
is `the same as' a von Neumann algebra. For a
unital Banach algebra ${\mathcal A}$ we write $\Her({\mathcal A})$ 
for the set of Hermitian elements; these are the 
elements $h \in {\mathcal A}$ such that
$$f(h) \in \mathbb{R} \; \text{for every state} \; f \; \text{on}
\;  \mathcal{A}.$$
Equivalently: $\Vert \exp (i t h) \Vert \leq 1$ for all $t \in \mathbb{R}$.  
The Hermitians in a $C^*$-algebra are exactly the self-adjoint elements.
We will be interested in Hermitians in the Banach algebra $B(X)$.

\smallskip
 
{\bf Acknowledgments.} This paper would not exist without the 
vision of E. G. Effros, who originally proposed the  program 
of generalizing the basic calculus of Banach space $M$-ideals
to `one-sided $M$-ideals' in operator spaces.
We  thank L. G. Brown and H. P. Rosenthal
for steering us toward some helpful papers, and 
Marius Junge and David Sherman for some helpful input during
conversations on
Morita equivalence and noncommutative $L^p$ spaces. 

\section{Preliminaries} \label{I}

\subsection{One-Sided Multipliers} \label{I.A}

Let $X$ be an operator space. Following \cite{Shilov} \S 4 (see also \cite{Wend}), we say that a 
map $T:X \to X$ is a \emph{left multiplier} of $X$ 
if there exists a linear complete isometry 
$\sigma:X \to B(H)$ and an operator $A \in B(H)$ such that
\begin{equation} \label{I.A.1}
	\sigma(Tx) = A\sigma(x)
\end{equation}
for all $x \in X$. In that case, we refer to $(\sigma,A)$ as an \emph{implementing pair} for $T$. 
It is easy to see that every left multiplier of $X$ is linear. We denote by $\Ml(X)$ the set of 
all left multipliers of $X$, and for $T \in \Ml(X)$, we define the \emph{multiplier norm} by
\begin{equation} \label{I.A.2}
	\|T\|_{\Ml(X)} = \inf\{\|A\|: (\sigma,A) \text{ is an implementing pair for } T\}.
\end{equation}
It turns out that the infimum in Equation (\ref{I.A.2}) is always achieved. 
This is a consequence of the existence of the
Arveson-Hamana Shilov boundary or
injective envelope \cite{Arv1,Arv2,Hamana,Shilov,BP01}. 
Clearly,
\begin{equation} \label{I.A.3}
	\|T\|_{cb} \leq \|T\|_{\Ml(X)}
\end{equation}
for all $T \in \Ml(X)$, so that $\Ml(X) \subset CB(X)$. In fact, strict inequality is possible in 
(\ref{I.A.3}) (see \cite{Shilov}), but this will not be of concern to us. 

It can be shown that $\Ml(X)$ is a unital Banach
algebra with respect to $\|\cdot\|_{\Ml(X)}$, and 
the usual composition product. In fact, it is a unital operator algebra with respect to the 
operator space structure induced by the canonical linear isomorphisms
\begin{equation} \label{I.A.4}
	M_n(\Ml(X)) \cong \Ml(C_n(X)).
\end{equation}

Whereas the definition of a left multiplier in terms of the existence of an implementing pair 
is extrinsic, there is an extremely useful intrinsic definition:

\begin{theorem}[\cite{BEZ}, Theorem 4.6] \label{BEZ}
Let $X$ be an operator space and $T:X \to X$ be linear. Then $T$ is an element of the closed unit 
ball of $\Ml(X)$ if and only if the map
\[
	\tau_T^c:C_2(X) \to C_2(X):\begin{bmatrix} x\\ y \end{bmatrix} \mapsto
		\begin{bmatrix} Tx\\ y \end{bmatrix}
\]
is completely contractive. 
\end{theorem}

Similar definitions and results hold for the 
\emph{right multiplier algebra} $\Mr(X)$. Left
multiplication must be replaced by right multiplication, and columns 
must be replaced by rows.  
There is a slight twist though---if
 we regard  $\M_r(X)$ as a subset of $CB(X)$,
we need to put
the reverse of the usual composition multiplication on
this subset of $CB(X)$.  The 
reason for this is that $X$ should be a right 
$\M_r(X)$ module.  Any left multiplier on $X$ commutes with
every right multiplier on $X$.
  
\subsection{One-Sided Adjointable Multipliers} \label{I.B}

Let $X$ be an operator space and $T:X \to X$ be a map. Again following \cite{Shilov} \S 4, we say 
that $T$ is a \emph{left adjointable map} of $X$ if there exists a linear complete isometry 
$\sigma:X \to B(H)$ and a map $S:X \to X$ such that
\begin{equation} \label{I.B.1}
	\sigma(Tx)^*\sigma(y) = \sigma(x)^*\sigma(Sy)
\end{equation}
for all $x, y \in X$. Left adjointable maps are linear and bounded. The collection of all left 
adjointable maps of $X$ will be denoted $\Al(X)$. If $T \in \Al(X)$, then there exists a unique map 
$S:X \to X$ satisfying Equation (\ref{I.B.1}). 
Henceforth, we will denote this map $T^\star$. 
We have
$T^\star \in \Al(X)$, and $(T^\star)^\star = T$. Every left adjointable map of $X$ is a left 
multiplier of $X$. That is, $\Al(X) \subset \Ml(X)$. In fact, $T \in \Al(X)$ if and only if there
exists a linear complete isometry $\sigma:X \to B(H)$ and an operator $A \in B(H)$ such that 
$(\sigma,A)$ is an implementing pair for $T$ and $A^*\sigma(X) \subset 
\sigma(X)$.  Indeed, $(\sigma,A^*)$ 
is then an implementing pair for $T^\star$. Because of this last fact, we will often refer to 
left adjointable maps as \emph{left adjointable multipliers}. For $T \in \Al(X)$,
\begin{equation} \label{I.B.2}
	\|T\|_{\Ml(X)} = \|T\|_{cb} = \|T\|.
\end{equation}
Also, $\Al(X)$ is a $C^*$-algebra with respect to  the involution $\star$ and the usual 
composition product. In fact, $\Al(X)$ is the `diagonal' of the operator algebra $\Ml(X)$ (we 
recall that any nonselfadjoint operator algebra $\mathcal{B} \subset B(H)$ contains a canonical 
$C^*$-algebra called the diagonal, namely $\mathcal{B} \cap \mathcal{B}^*$, and the latter is 
well-defined independently of the representation of $\mathcal{B}$ on a Hilbert space).   

There are intrinsic 
characterizations of certain classes of left adjointable multipliers, analogous 
to Theorem \ref{BEZ}. Namely,

\begin{theorem}[\cite{BEZ}, Corollary 4.8 and Theorem 4.9] \label{Hermitian}
Let $X$ be an operator space and $T:X \to X$ be linear. Then 
\begin{enumerate}
\item[(i)] $T$ is a unitary element of $\Al(X)$ if and only if $\tau_T^c:C_2(X) \to C_2(X)$ is a 
completely isometric surjection.
\item[(ii)] $T$ is a self-adjoint element of $\Al(X)$ if and only if $\tau_T^c$ is a Hermitian 
element of $CB(C_2(X))$.
\end{enumerate}
\end{theorem}

Using Theorem \ref{BEZ}, one may show: 

\begin{theorem}[\cite{BEZ}, Theorem 5.4 and \cite{Dual}, Corollary 3.2] \label{wst}
For an operator space $X$, $\Ml(X^*)$ is a dual operator algebra, and $\Al(X^*)$ is a 
$W^*$-algebra. A bounded net $(T_i)$ in $\Ml(X^*)$ (resp. in $\Al(X^*)$) converges weak-* to $T$ 
in $\Ml(X^*)$ (resp. in $\Al(X^*)$) if and only if 
$T_i  \to T$ in the point-weak-* topology
(that is,  if and only if $T_if \to Tf$ weak-* for all $f \in X^*$).
\end{theorem}
 
Furthermore, every element of $\Al(X^*)$ is weak-* continuous 
(\cite{BEZ}, Theorem 5.5). On the 
other hand, it is unknown if elements of $\Ml(X^*)$ are automatically weak-* continuous.\\

Similar definitions and results hold for the \emph{right adjointable multiplier algebra} $\Ar(X)$. 
Note, however, that in the 
development of $\Ar(X)$ as opposed to $\Al(X)$, there is a slight twist.
Considered as a subset of $CB(X)$, the space $\Ar(X)$  should be
given the opposite of the usual multiplication on $CB(X)$.

The algebras $\Ml(X)$, $\Al(X)$, $\Mr(X)$, and $\Ar(X)$ are invariants
for the operator space $X$. If $\phi:X \to Y$ is a completely isometric surjection, then the map
$\Phi:CB(X) \to CB(Y):T \mapsto \phi \circ T \circ \phi^{-1}$ restricts to a completely isometric
isomorphism $\Ml(X) \to \Ml(Y)$ (resp. $\Mr(X) \to \Mr(Y)$) and to a *-isomorphism $\Al(X) \to \Al(Y)$
(resp $\Ar(X) \to \Ar(Y)$). This can be seen in a number of ways.

Finally, we remark that \cite{BP01} presents a very effective framework 
for the above multiplier algebras, in terms of the injective
envelope.  Since we do not essentially need this perspective,
we shall not emphasize it here, although it
is often useful.    

\subsection{One-Sided $M$- and $L$-Structure} \label{I.C}

Let $X$ be an operator space. Following \cite{BEZ}, we define a \emph{complete left $M$-projection} on
$X$ to be an orthogonal projection $P \in \Al(X)$. Actually, though it is abusive to do so, we will
abbreviate `complete left $M$-projection' by `left $M$-projection'. 
Hopefully the reader will agree 
that the savings in verbiage 
outweigh the possible confusion with \cite{BEZ}. Being a left 
$M$-projection is equivalent to a number of other conditions, many of which are 
often easier to verify 
in practice:

\begin{theorem}[\cite{BEZ}, Proposition 3.2 and Theorem 5.1] \label{left_M}
Let $X$ be an operator space and $P:X \to X$ be a linear idempotent map. Then the following
are equivalent:
\begin{enumerate}
\item[(i)] $P$ is a left $M$-projection.
\item[(ii)] The map $\tau_P^c:C_2(X) \to C_2(X)$ is completely contractive. Equivalently,
by Theorem \ref{BEZ}, $P$ is an element of the closed unit ball of $\Ml(X)$.
\item[(iii)] The map $\tau_P^c$ is a Hermitian element of $CB(C_2(X))$. Equivalently, by
Theorem \ref{Hermitian}, $P \in \Al(X)_{sa}$.
\item[(iv)] The map $\nu_P^c:X \to C_2(X):x \mapsto \begin{bmatrix} Px\\ (\Id - P)x \end{bmatrix}$
is completely isometric.
\item[(v)] The maps $\nu_P^c$ and $\mu_P^c:C_2(X) \to X:\begin{bmatrix} x\\ y \end{bmatrix}
\mapsto Px + (\Id - P)y$ are completely contractive.
\item[(vi)]  There exists  a completely isometric embedding
$\sigma  \colon X \hookrightarrow B(H)$, and a  projection
$e \in B(H)$ such that $\sigma(Px) = e \sigma(x)$
for all $x \in X$.
\end{enumerate}
\end{theorem}

A linear subspace $J$ of $X$ is a \emph{right $M$-summand} of $X$ if it is the range of a
left $M$-projection $P$ on $X$.  The kernel of $P$, which equals 
the range of $I-P$, is also a right $M$-summand; it is called 
the {\em complementary right $M$-summand} to $J$.   We often 
write this complementary summand as $\tilde{J}$, and write
$X = J \oplus_{\rM} \tilde{J}$.    Another way of stating this identity,
by equivalence (iv)
in the last theorem, is that there is a completely isometric embedding 
$\sigma$  of $X$ into $B(L,H \oplus K)$, such that 
$P_K \sigma(J) = P_H \sigma(\tilde{J}) = 0$.    

Since left $M$-projections are (completely) contractive, right 
$M$-summands are automatically closed. Furthermore, since left $M$-projections on a dual operator
space are weak-* continuous, right $M$-summands on such spaces are automatically weak-* closed.
A closed linear subspace $J$ of an operator space $X$ is a \emph{right $M$-ideal}
of $X$ if its double annihilator, $J^{\perp\perp}$, is a right $M$-summand of $X^{**}$. Every
right $M$-summand of $X$ is a right $M$-ideal of $X$, but the converse is false.  \\

Dual to the one-sided $M$-structure of an operator space $X$
 is its one-sided $L$-structure.
A linear idempotent map $P:X \to X$ is a \emph{right $L$-projection} on $X$ if $P^*:X^* \to X^*$ 
is a left $M$-projection on $X^*$. It is useful, as in the case of left $M$-projections, to have
alternative characterizations of right $L$-projections:

\begin{proposition}[\cite{BEZ}, Proposition 3.4 and Corollary 3.5] \label{right_L}
Let $X$ be an operator space and $P:X \to X$ be a linear idempotent map. Then the following are
equivalent:
\begin{enumerate}
\item[(i)] $P$ is a right $L$-projection.
\item[(ii)] The map $\nu_P^r:X \to R_2 \ptimes X:x \mapsto \begin{bmatrix} Px & (\Id - P)x 
\end{bmatrix}$ is a complete isometry.
\item[(iii)] The maps $\nu_P^r$ and $\mu_P^r:R_2 \ptimes X \to X:\begin{bmatrix} x & y \end{bmatrix}
\mapsto Px + (\Id - P)y$ are completely contractive.
\end{enumerate}
Here $\ptimes$ is the operator space projective tensor product.
\end{proposition} 

A linear subspace $J$ of $X$ is a \emph{left $L$-summand} of $X$ if it is the range of a right 
$L$-projection. As with right $M$-summands, the fact that left $L$-summands are closed is
automatic. On the other hand, there is no need to define the concept of a \emph{left $L$-ideal}---a 
closed linear subspace $J$ of $X$ is a left $L$-summand of $X$ if and only if its double
annihilator, $J^{\perp\perp}$, is a left $L$-summand of $X^{**}$ (\cite{BEZ}, Proposition 3.9).

We end this section with some assorted results about the 
one-sided $M$- and $L$-structure of a
general operator space $X$. These can be found in 
\cite{BEZ} or \cite{Zarikian}.  The reader may want to 
take a few moments to absorb these, since they will be used freely 
and usually silently in the sequel.

\begin{itemize}
\item Let $P:X \to X$ be a bounded linear idempotent map. Then $P$ is a left $M$-projection on $X$ if and only
if $P^*$ is a right $L$-projection on $X^*$.
\item A closed linear subspace $J$ of $X$ is a right $M$-ideal of $X$ if and only if $J^\perp$ is
a left $L$-summand of $X^*$. A closed linear subspace $J$ of $X$ is a left $L$-summand of $X$ if and only 
if $J^\perp$ is a right $M$-summand of $X^*$.
\item Every right $M$-summand (resp. left $L$-summand) is the range of a 
unique left $M$-projection (resp. right $L$-projection).
\item If a right $M$-ideal $J$ is the range of a contractive 
projection $P$, then it is in fact a right $M$-summand
(and $P$ is the unique left $M$-projection onto $J$). In
particular, a right
 $M$-ideal which is also a dual Banach space must be a right $M$-summand (cf. Lemma
\ref{App.A.1}). Hence, a weak-* closed right $M$-ideal of $X^*$ is a right $M$-summand of $X^*$.
\item Every right $M$-ideal $J$ of $X$ is `Hahn-Banach smooth', 
which is to say that each $f \in J^*$ has a unique
norm-preserving extension $\tilde{f} \in X^*$. 
Indeed, it is easy to see that any right $M$-ideal $J$ 
is an {\em HB-subspace} (see \cite[p. 44]{HWW}), and therefore
enjoys the properties of such spaces.
In particular $J^*$ may be identified with a closed linear subspace
of $X^*$, none other than 
the complementary left $L$-summand of $J^\perp$. That is,
$X^* = J^\perp \oplus_{\lL} J^*$.
\end{itemize}

Another useful sufficient condition for a  right $M$-ideal 
to be a right $M$-summand  comes from a fact from \cite{GKS}.  
Namely, 
those authors define a class of subspaces of a Banach space 
called {\em $u$-ideals}.   It is easy to see that 
every one-sided $M$-ideal is a $u$-ideal.   On the other hand,
from  \cite[p.14]{GKS} we know that any 
$u$-ideal not containing a linear topologically isomorphic
copy of $c_0$ is a {\em $u$-summand}.  It is easy to see that
any $u$-summand is contractively complemented.  From this
and the second last `bullet' above, we deduce that:

\begin{itemize} 
\item Any one-sided $M$-ideal 
not containing a linear topologically isomorphic
copy of $c_0$ is a one-sided $M$-summand.  
\end{itemize}

\subsection{The One-Sided Cunningham Algebra} \label{I.D}

Let $X$ be an operator space. Following \cite{Zarikian}, we define the \emph{right Cunningham algebra} 
of $X$ to 
be the closed linear span in $CB(X)$ of the right $L$-projections on $X$:
\[
	\Cr(X) = \overline{\spn}\{P: P \text{ is a right $L$-projection on $X$}\}.
\]
The closure here is in the associated norm, namely the 
`cb' norm on $CB(X)$.
For $T \in \Cr(X)$, one has that
\begin{equation} \label{I.D.1}
	\|T\|_{cb} = \|T\|.
\end{equation}
This follows from the fact that for any $T$ in the span of the right $L$-projections on $X$, $T^*$ is 
in the span of the left $M$-projections on $X^*$, and therefore is an element of $\Al(X^*)$, so that
\[
	\|T\|_{cb} = \|T^*\|_{cb} = \|T^*\| = \|T\|
\]
by Equation (\ref{I.B.2}). Clearly, the map $CB(X) \to CB(X^*):T \mapsto T^*$ restricts to a linear 
isometry $\rho:\Cr(X) \to \Al(X^*)$. Because $\Al(X^*)$ is a $W^*$-algebra (Theorem 
\ref{wst}), it is 
the norm-closed linear span of its projections (the left $M$-projections on $X^*$). Since every
left $M$-projection on $X^*$ is weak-* continuous (cf. Section \ref{I.C}), and 
thus the adjoint of
a right $L$-projection on $X$, $\rho$ is surjective.
It follows that $\Cr(X)$ is a (closed) subalgebra of $CB(X)$, and
that the contractive projections in $\Cr(X)$ are exactly the 
right $L$-projections on $X$.
Endowed with the involution defined by
\[
	T^\star = \rho^{-1}(\rho(T)^\star),
\]
$\Cr(X)$ becomes a $C^*$-algebra. Obviously, $\rho$ is then an isometric *-anti-isomorphism (i.e. 
linear, adjoint-preserving, and product-reversing). By Sakai's theorem 
\cite{Sakai}, $\Cr(X)$ is 
actually a $W^*$-algebra with unique predual, 
and $\rho$ is a weak-* homeomorphism. A bounded net $(T_i)$ in $\Cr(X)$ 
converges weak-* to 
$T \in \Cr(X)$ if and only if $T_i(x) \to T(x)$ weakly for all $x \in X$.

Similar definitions and results hold for the \emph{left Cunningham algebra}, $\Cl(X)$.

\section{Spatial Action} \label{II}
 
We have seen that
to each operator space $X$, there are
associated natural $C^*$-algebras $\Al(X)$ and $\Ar(X)$, and $W^*$-algebras $\Cr(X)$ and $\Cl(X)$. For
dual operator spaces, the $C^*$-algebras are in fact $W^*$-algebras. These operator algebras are not
concretely given, which is to say that they do not appear acting on a Hilbert space. On the other hand,
they are not totally `space-free': 
they do act on $X$. So one can hope to employ spatial intuition and arguments.
Of course, one must exercise caution in doing so---$X$ is typically far from being a Hilbert space. In this
section we carefully attempt to understand how 
projections and partial isometries in our abstract
operator algebras manipulate the underlying space $X$.
 
\subsection{Projections} \label{II.A}

First we investigate the partial ordering of projections. Recall that for projections $E$ and $F$ in an
abstract $C^*$-algebra $\mathcal{A}$, $E \leq F$ if and only if $EF = FE = E$. 
If $\mathcal{A}$ is a $*$-subalgebra of $B(H)$,
then $E \leq F$ if and only if $E(H) \subset F(H)$. It is reassuring that the analogous result holds in
our setting.

\begin{proposition} \label{II.A.1}
Let $X$ be an operator space, and $P$ and $Q$ be left $M$-projections (resp. right $L$-projections) on $X$. Then 
$P \leq Q$ if and only if $P(X) \subset Q(X)$.
\end{proposition}

\begin{proof}
Suppose $PQ = QP = P$. Then $P(X) \subset Q(X)$. Conversely, suppose that
 $P(X) \subset Q(X)$. Then $QP = P$.
But then $PQ = (QP)^\star = P^\star = P$.
\end{proof}

Next we investigate the lattice of projections. Recall that for an abstract $W^*$-algebra $\mathcal{A}$,
the projections in $\mathcal{A}$ form a complete lattice with respect to the partial ordering of
projections. If $\mathcal{A} \subset B(H)$ and $\{E_i: i \in I\}$ is a family of projections in
$\mathcal{A}$, then $\bigwedge_{i \in I} E_i$ is the projection onto $\bigcap_{i \in I} E_i(H)$ and
$\bigvee_{i \in I} E_i$ is the projection onto $\overline{\spn}\{\bigcup_{i \in I} E_i(H)\}$. To determine
what happens in our setting, we need abstract (space-free) formulas for the meet and join of projections
in a $W^*$-algebra. Such formulas are well-known to the experts, but are hard to track down in the
literature. We include them here for the reader's convenience.

\begin{lemma} \label{II.A.2}
Let $\mathcal{A} \subset B(H)$ be a $W^*$-algebra and $E, F \in \mathcal{A}$ be projections. Then
\[
	E \vee F = \wlim_{n \to \infty} (E + F)^{1/n}
\]
and
\[
	E \wedge F = \wlim_{n \to \infty} I - (2I - E - F)^{1/n},
\]
where $I \in B(H)$ is the identity.
\end{lemma}

\begin{proof}
Let $T \in \mathcal{A}$ be such that $0 \leq T \leq I$. 
Since $x^{1/n} \nearrow \chi_{(0,1]}(x)$ for all
$x \in [0,1]$, we have
$T^{1/n} \to R(T)$ weak-*, 
where $R(T) \in \mathcal{A}$ is the projection onto the closure 
of the range of $T$ (Borel functional calculus). It follows that if $T \in \mathcal{A}_+$, then
$T^{1/n} \to R(T)$ weak-* as well. In particular, $(E + F)^{1/n} \to R(E + F)$ weak-*. But one easily
verifies that $\ker(E + F) = \ker(E) \cap \ker(F)$, which implies that
\[
	\overline{\ran(E + F)} = \ker(E + F)^\perp = \overline{\ker(E)^\perp + \ker(F)^\perp} =
		\overline{\ran(E) + \ran(F)} = \ran(E \vee F),
\]
so that $R(E + F) = E \vee F$. Hence,
\[
	E \vee F = \wlim_{n \to \infty} (E + F)^{1/n}.
\]
The remaining formula follows from the first formula and the fact that 
$E \wedge F = I - (I - E) \vee (I - F)$.
\end{proof}

We now apply Lemma \ref{II.A.2} to our situation.

\begin{proposition} \label{II.A.3}
Let $X$ be a dual operator space, and $P$ and $Q$ be left $M$-projections on $X$. Then
\begin{enumerate}
\item[(i)] $(P \wedge Q)(X) = P(X) \cap Q(X)$.
\item[(ii)] $(P \vee Q)(X) = \overline{P(X) + Q(X)}^{\wks}$.
\end{enumerate}
\end{proposition}

\begin{proof}
(i) By Proposition \ref{II.A.1}, $(P \wedge Q)(X) \subset P(X) \cap Q(X)$. Conversely, suppose that
$x \in P(X) \cap Q(X)$. 
Then $(2\Id - P - Q)(x) = 0$. Thus, $f(2\Id - P - Q)(x) = 0$ for any polynomial $f$ without constant term. 
Hence, $(2\Id - P - Q)^{1/n}(x) = 0$ for all $n \in 
\mathbb{N}$. But then by Lemma \ref{II.A.2} 
and Theorem \ref{wst}, we have that  
\[
	(P \wedge Q)(x) = \wlim_{n \to \infty} (\Id - (2\Id - P - Q)^{1/n})(x) = x,
\]
which says that $x \in (P \wedge Q)(X)$. 

(ii) Again by Proposition \ref{II.A.1}, $P(X) + Q(X) \subset (P \vee Q)(X)$.
Since $(P \vee Q)(X)$ is weak-* closed,
\[
	\overline{P(X) + Q(X)}^{\wks} \subset (P \vee Q)(X).
\]
Conversely, suppose that
 $x \in (P \vee Q)(X)$. Then by Lemma \ref{II.A.2} and Theorem \ref{wst} again,
we have that
\[
	x = (P \vee Q)(x) = \wlim_{n \to \infty} (P + Q)^{1/n}(x).
\]
Now for each fixed $n \in \mathbb{N}$, there exist polynomials $f_k$ without constant terms such 
that $f_k(P + Q) \to (P + Q)^{1/n}$. Therefore,
\[
	(P + Q)^{1/n}(x) = \lim_{k \to \infty} f_k(P + Q)(x) \in \overline{P(X) + Q(X)}.
\]
Since this is true for each $n \in \mathbb{N}$, $x \in \overline{P(X) + Q(X)}^{\wks}$.
\end{proof}

\begin{corollary} \label{II.A.4}
Let $X$ be an operator space, and $P$ and $Q$ be right $L$-projections on $X$. Then
\begin{enumerate}
\item[(i)] $(P \wedge Q)(X) = P(X) \cap Q(X)$.
\item[(ii)] $(P \vee Q)(X) = \overline{P(X) + Q(X)}$.
\end{enumerate}
\end{corollary}

The proof is the same as that of Proposition \ref{II.A.3}, 
modulo the obvious modifications. Alternatively,
one can deduce Corollary \ref{II.A.4} from Proposition \ref{II.A.3} by exploiting the duality between
one-sided $M$-structure and one-sided $L$-structure.\\

The passage from the finite to the infinite is straightforward.

\begin{proposition} \label{II.A.5}
\begin{enumerate}
\item[(i)] If $X$ is a dual operator space and $\{P_i: i \in I\}$ is a family of left $M$-projections on
$X$, then $\left(\bigwedge_{i \in I} P_i\right)(X) = \cap_{i \in I} P_i(X)$ and 
$\left(\bigvee_{i \in I} P_i\right)(X) = \overline{\spn}^{\wks}\left\{\bigcup_{i \in I} P_i(X)\right\}$.
\item[(ii)] If $Y$ is an operator space and $\{Q_j: j \in J\}$ is a family of right $L$-projections on
$Y$, then $\left(\bigwedge_{j \in J} Q_j\right)(Y) = \cap_{j \in J} Q_j(Y)$ and 
$\left(\bigvee_{j \in J} Q_j\right)(Y) = \overline{\spn}\left\{\bigcup_{j \in J} Q_j(Y)\right\}$.
\end{enumerate}
\end{proposition}

\begin{proof}
We will prove half of each of the two assertions. The reader will have no trouble supplying the missing
arguments. 

(i) For each finite subset $F$ of $I$, define $P_F = \bigwedge_{i \in F} P_i$. By Proposition
\ref{II.A.3}, $P_F(X) = \cap_{i \in F} P_i(X)$. Also, $P_F \to P \equiv \bigwedge_{i \in I} P_i$ weak-*.
Since $P \leq P_i$ for all $i \in I$, $P(X) \subset P_i(X)$ for all $i \in I$ (Proposition 
\ref{II.A.1}).
Therefore, $P(X) \subset \cap_{i \in I} P_i(X)$. On the other hand, if $x \in \cap_{i \in I} P_i(X)$,
then $P(x) = \wlim_F P_F(x) = x$, so that $x \in P(X)$. Hence, $\cap_{i \in I} P_i(X) \subset P(X)$.

(ii) For each finite subset $F$ of $J$, define $Q_F = \bigvee_{j \in F} Q_j$. By Corollary 
\ref{II.A.4},
$Q_F(X) = \overline{\spn}\left\{\bigcup_{j \in F} Q_j(X)\right\}$. Also, 
$Q_F \to Q \equiv \bigvee_{j \in J} Q_j$ weak-*. Since $Q_j \leq Q$ for all $j \in J$, $Q_j(X) \subset Q(X)$ 
for all $j \in J$, which implies that $\overline{\spn}\left\{\bigcup_{j \in J} Q_j(X)\right\} \subset Q(X)$. 
On the other hand, if $x \in Q(X)$, then $Q_F(x) \to x$ weakly. Thus,
$x \in \overline{\spn}^{\wk}\left\{\bigcup_{j \in J} Q_j(X)\right\} 
= \overline{\spn}\left\{\bigcup_{j \in J} Q_j(X)\right\}$. Since the choice of $x$ was arbitrary,
$Q(X) \subset \overline{\spn}\left\{\bigcup_{j \in J} Q_j(X)\right\}$.
\end{proof}

\begin{corollary} \label{latis}  
\begin{itemize} \item [(i)]  If $X$ is an operator space, then there 
is a lattice isomorphism between the complete lattice of right $L$-projections
on $X$, and the lattice of left $L$-summands of $X$.
\item [(ii)]  Let $X$ be a dual operator space, and define the 
`meet' of a family of subspaces of $X$ to be their intersection, and the
`join' to be the weak-* closure of their span.  Then there
is a lattice isomorphism between the complete lattice of left $M$-projections
on $X$, and the lattice of right $M$-summands of $X$.
\end{itemize} 
\end{corollary}

\subsection{Partial Isometries} \label{II.B}

Next we investigate partial isometries $W$ in $\Al(X)$. Recall that an element 
$W$ of a $C^*$-algebra $\mathcal{A}$ is
a partial isometry if $W^*W$ is a projection.  We assume familiarity 
with the basic facts about partial isometries on Hilbert spaces (see \cite{KR_II}).

\begin{proposition} \label{II.B.1}
Let $X$ be an operator space and $W \in \Al(X)$ be a partial isometry. Then $P = W^\star W$ and $Q = WW^\star$
are left $M$-projections on $X$. One has that
\[
	\ran(P) = \ran(W^\star), \; \ker(P) = \ker(W), \;
\ran(Q) = \ran(W), \text{ and } \ker(Q) = \ker(W^\star).
\]
In particular, $\ran(W)$ and $\ran(W^\star)$ are closed and
\[
	X = \ran(W^\star) \oplus_{\rM} \ker(W) = \ran(W) \oplus_{\rM} \ker(W^\star).
\]
Finally, $W$ maps $\ran(W^\star)$ completely isometrically onto $\ran(W)$.
\end{proposition}

\begin{proof}
Clearly $P$ and $Q$ are left $M$-projections on $X$.
We have that
\[
	\ran(P) = \ran(W^\star W) \subset \ran(W^\star) = \ran(W^\star WW^\star) = \ran(PW^\star)
		\subset \ran(P).
\]
Likewise, $\ran(Q) = \ran(W^\star)$. Clearly, $Wx = 0$ implies $Px = W^\star Wx = 0$. Conversely, $Px = 0$
implies $Wx = WW^\star Wx = WPx = 0$. Hence, $\ker(P) = \ker(W)$. Likewise, $\ker(Q) = \ker(W^\star)$. Clearly,
$W$ is a linear bijection between $\ran(W^\star)$ and $\ran(W)$ (with inverse $W^\star$). Since both
$W$ and $W^\star$ are completely contractive, the last assertion follows.
\end{proof}

The same result is true if we replace $\Al(X)$ by $\Cr(X)$, left $M$-projections by right $L$-projections, and
$\oplus_{\rM}$ by $\oplus_{\lL}$.
 
\subsection{Murray-Von Neumann Equivalence} \label{II.BB}

Let $X$ be an operator space and let 
$P$ and $Q$ be left $M$-projections on $X$. 
We say that $P$ and $Q$ are {\em (left) Murray-von Neumann
equivalent} (or simply {\em equivalent}),
 and write $P \sim Q$, if there exists a 
partial isometry $V \in \Al(X)$ such that $V^\star V = P$
and $VV^\star = Q$ (i.e. if they are equivalent in $\Al(X)$). 
Suppose $P \sim Q$ (via $V$), and let $J$ and $K$ be the
right $M$-summands corresponding to $P$ and $Q$, respectively. Denote by $\tilde{J}$ the complementary right $M$-summand
of $J$, and by $\tilde{K}$ the complementary right $M$-summand of $K$. 
By Proposition \ref{II.B.1}, we have that $\ker(V) = \tilde{J}$
and $V$ maps $J$ completely isometrically onto $K$.  If, further,
$X$ is a dual operator space, then as we stated after Theorem
\ref{wst} and Theorem \ref{left_M}, $V$ and $V^\star$ are
weak-* continuous, and $J$ and $K$ are  weak-* closed.
Thus $J$ is also weak-* homeomorphic to $K$.

It is clear that if $\Al(X)$ is commutative, then 
Murray-von Neumann equivalent left $M$-projections on $X$ must be equal.
Conversely, it follows easily from the well-known von Neumann
algebra projection calculus that if $P$ and $Q$ are 
noncommuting left $M$-projections on a dual operator space $X$,
then $P$ and $P^\perp \equiv \Id - P$ have
nonzero subprojections which are
Murray-von Neumann equivalent (likewise for $Q$ and $Q^\perp$). Indeed, a characterization
of noncommutativity in von Neumann algebras is the 
existence of two equivalent mutually
orthogonal nonzero projections. 

\subsection{Inner Products on Operator Spaces} \label{II.BBB}

A perspective emphasized in \cite{Shilov} is that a 
general operator space $X$ may be studied using 
the framework of Hilbert $C^*$-modules (see e.g. 
\cite{La}).  Indeed, every operator space $X$ has 
canonical $C^*$-algebra-valued inner products,
which we will call the {\em Shilov inner products}.  
These may be viewed as the natural inner products
on the Hamana triple envelope of $X$, or may be described
in terms of the injective envelope.  Since the reader may not 
be familiar with these notions, we will not go into much detail.
Suffice it to say that a `left Shilov inner product' 
$\langle \cdot , \cdot \rangle$ is a map from 
$X \times X$ into a $C^*$-algebra, which is linear in the second
variable and conjugate linear in the first variable. It is explicitly described
at the top of p. 307 in \cite{Shilov}, and it has the following
property: $T \in \Al(X)$ if and only if
there is a map $S$ on $X$ with
\begin{equation} \label{adeq}
\langle Tx , y \rangle = \langle x , S y \rangle \end{equation}
for all $x, y \in X$.     Also $\Vert \langle x , x \rangle \Vert 
= \Vert x \Vert^2$ for $x \in X$.     In fact we shall 
be more concerned here with the associated notion 
of orthogonality with respect to the inner product: we  
write $x \rightthreetimes y$ if and only if
$\langle x , y \rangle = 0$. This has 
a simple reformulation: 
from the universal property 
of the Hamana triple envelope it is easy to see that 
$x \rightthreetimes y$ if and only if there 
exists a completely isometric embedding $\sigma$ of $X$ into a
$C^*$-algebra with 
\begin{equation} \label{pere}
\sigma(x)^* \sigma(y) = 0.
\end{equation} 
Similarly, for subsets $B, C \subset X$, we write  
$B \rightthreetimes C$ if $x \rightthreetimes y$ for 
all $x \in B, y \in C$.  It can be easily shown that this
is equivalent to the existence of a map $\sigma$ as above
such that (\ref{pere}) holds for every choice of $x \in B, y \in C$.
It follows that if $x \rightthreetimes y$ in an operator
space $X$, and if $x, y \in Y \subset X$, then
$x \rightthreetimes y$ in $Y$. (The converse of this 
is false.)

In \cite[Theorem 5.1]{BEZ} we showed that a linear
idempotent map $P$ on an operator space $X$ is a left $M$-projection
if and only if $P(X) \rightthreetimes (I-P)(X)$.   It follows
that if $J, K$ are linear subspaces of $X$ with
$J + K = X$ and $J \rightthreetimes K$, then $J$ is a
right $M$-summand of $X$ and $K$ is the complementary
right $M$-summand.

\subsection{Polar Decomposition} \label{II.C} 

In this section we determine when a left adjointable multiplier has a `nice' polar decomposition.
We shall see that on a general operator space $X$, not every $T \in \Al(X)$ is
`polar decomposable'. We give necessary and sufficient conditions for this to occur (Theorem
\ref{polar}). One important instance when these conditions are met is when $T$ has closed
range (Theorem \ref{pd}), and we draw a number of conclusions from this. Of course, if
$X$ is a dual operator space, then $T$ automatically has a nice polar decomposition
(Theorem \ref{II.B.2}). Again, this has favorable consequences. As the reader shall see, much of this section 
is a good illustration of the remark in the previous section that a general operator space $X$ 
may be studied using the framework of Hilbert $C^*$-modules.

For  $T \in \Al(X)$, we write $|T|$ for the square
root of $T^\star T$ in $\Al(X)$. Note that if 
$\langle \cdot , \cdot \rangle$ is
the `left Shilov inner product' on $X$ (see Section \ref{II.BBB}),
 then $$\langle T x , T x \rangle = \langle T^\star T x , x \rangle = \langle |T| x , |T| x \rangle,$$
for $x \in X$, using (\ref{adeq}). This together with the fact that 
$\Vert \langle x , x \rangle \Vert = \Vert x \Vert^2$
shows that $\ker(T) = \ker(|T|)$.

As is remarked in passing in \cite{Shilov}, Appendix B, the 
entire discussion in Section 15.3 of \cite{Wegge} on 
the polar decomposition of adjointable operators on a Hilbert 
$C^*$-module carries over 
verbatim to the setting of left adjointable multipliers of an operator space. 
Using arguments identical to those
found in \cite[Section 15.3]{Wegge}, one obtains the following two theorems:

\begin{theorem} \label{polar}
Let $X$ be an operator space and $T \in \Al(X)$. Then the following are equivalent:
\begin{enumerate}
\item[(i)] $T = W|T|$ for some partial isometry $W \in \Al(X)$ such that $\ran(W^\star) = \overline{\ran(|T|)}$ and
$\ran(W) = \overline{\ran(T)}$.
\item[(ii)] $X = \overline{\ran(|T|)} \oplus_{\rM} \ker(|T|) = \overline{\ran(T)} \oplus_{\rM} \ker(T^\star)$.
\item[(iii)] $\overline{\ran(|T|)}$ and $\overline{\ran(T)}$ are right $M$-summands of $X$.
\end{enumerate}
\end{theorem}

Of course, it follows too that if (i)-(iii) of the last theorem hold, then
$\ker(|T|)$ and $\ker(T^\star)$ are right $M$-summands of $X$.

\begin{theorem} \label{pd}  Suppose that $T$ is a left
adjointable multiplier of an operator space $X$.
If $\ran(T)$ is closed in $X$,
then $T^\star$ and $|T|$ have closed ranges as well.  Also,
$\ker(T)$, $\ker(T^\star)$, $\ran(T)$, and $\ran(T^\star)$ are
right M-summands of $X$.  We have:
$$X = \; \ker(T) \oplus_{\rM} \ran(T^\star) \; = \; \ker(T^\star)
\oplus_{\rM}
\ran(T) \; \; .
$$  Finally, $T$ has a polar decomposition $T = V |T|$ for a left
adjointable partial isometry $V$ satisfying
$\ker(V) \; = \; \ker(T) \; , \; \ker(V^\star) \; = \;
\ker(T^\star) \; , \; \ran(V) = \ran(T)$, and $\ran(V^\star) = \ran(T^\star)$.
\end{theorem}

Using Theorem \ref{pd}, one may draw a number of important conclusions.

\begin{corollary} \label{II.B.1b} 
A closed linear subspace $J$ of an operator space $X$ is a right 
$M$-summand if and only if it is the range of some $T \in \Al(X)$.
\end{corollary}

Interestingly, for a left adjointable multiplier $T$ on a dual operator space $X$, this says that
$\ran(T)$ is closed if and only if $\ran(T)$ is weak-* closed.

\begin{corollary} \label{inv}  If $T \in \Al(X)$ for
an operator space $X$, then  the following are
equivalent:
\begin{itemize}
\item [(i)]  $T$ is a bijection,
\item [(ii)]   $T$ is
invertible in $\Al(X)$,
\item [(iii)]  $T$ and $T^\star$ are bounded away from zero. That is,
there exists $\epsilon > 0$ such that
 $\Vert Tx \Vert \geq \epsilon \Vert x \Vert$ for all
$x \in X$, and similarly for $T^\star$.
\end{itemize}
In this case we have $(T^{-1})^\star =
(T^\star)^{-1}$.
\end{corollary}
 
\begin{proof}  If $T$ is bijective, then
by the displayed equation in  \ref{pd}, we have
$\ran(T^\star) = X$ and
$\ker(T^\star)  \; = \;   \{ 0 \}$ .
So $T^\star$ is bijective.   If $x, y \in X$, then
$$\langle T^{-1} x , y \rangle \; = \;
\langle T^{-1} x , T^\star ((T^\star)^{-1}(y)) \rangle
\; = \;
\langle  x , (T^\star)^{-1}(y) \rangle \; \; , $$
where $\langle  \cdot , \cdot  \rangle$ is the 
`inner product' discussed in \ref{II.BBB}.
This gives (ii) as well as the final formula.  That (ii) $\Rightarrow$ (iii)
is clear.  Finally, given (iii), we see
by the displayed equation in  \ref{pd} that
$T$ is onto, and we obtain (i).
\end{proof}
  
The following is related to
a result of Lance for Hilbert $C^*$-modules \cite[Theorem 3.5]{La}:
 
\begin{corollary} \label{uni}  
Let $X$ be an operator space. A surjective
left adjointable isometry $T : X \rightarrow X$
is a unitary in $\Al(X)$.
\end{corollary}

\begin{proof}  By the last result, $T^{-1}$ is left adjointable.
Since $T$ and $T^{-1}$ are elements of norm one in 
the $C^*$-algebra $\Al(X)$, and since their product is the identity in either order,
they must be unitary.   \end{proof}   

Now suppose that $X$ is a dual operator 
space. Then $\Al(X)$ is a $W^*$-algebra, which implies that
every $T \in \Al(X)$ satisfies $T = W|T|$, where $W \in \Al(X)$ is a partial isometry such that
$W^\star W$ is the left support of $|T|$ (the smallest projection $P \in \Al(X)$ such that $P|T| = |T|$)
and $WW^\star$ is the left support of $T$. This simple observation yields the following result.

\begin{theorem} \label{II.B.2}
Let $X$ be a dual operator space and $T \in \Al(X)$. Then $T = W|T|$, where $W \in \Al(X)$ is a partial
isometry such that $\ran(W^\star) = \overline{\ran(|T|)}^{\wks}$ and $\ran(W) = \overline{\ran(T)}^{\wks}$.
Furthermore, $\ker(W) = \ker(|T|)$ and $\ker(W^\star) = \ker(T^\star)$. In particular,
$X = \overline{\ran(T)}^{\wks} \oplus_{\rM} \ker(T^\star)$.
\end{theorem}

\begin{proof}
Let $W \in \Al(X)$ be the partial isometry described in the paragraph preceding the proposition. We have that
\[
	\ran(|T|) = \ran(W^\star T) = \ran(W^\star WW^\star T) \subset \ran(W^\star W).
\]
Since $\ran(W^\star W)$ is weak-* closed, $\overline{\ran(|T|)}^{\wks} \subset \ran(W^\star W)$. On the other 
hand, we have by the proof of Lemma \ref{II.A.2} that
\[
	W^\star W = \wlim_{n \to \infty} |T|^{1/n}.
\]
By spectral theory (i.e.\ polynomial approximation via the functional calculus)
it is easy to see that
$\ran(|T|^{1/n}) \subset \overline{\ran(|T|)}$ 
for all $n \in \mathbb{N}$.  Hence  
$\ran(W^\star W) \subset \overline{\ran(|T|)}^{\wks}$. Now suppose 
that $|T|x = 0$. Then $|T|^{1/n}x = 0$ for all
$n \in \mathbb{N}$, which implies that $W^\star Wx = 0$. Conversely, suppose
that $W^\star Wx = 0$. Then $Wx = 0$,
which implies that $|T|x = T^\star Wx = 0$. Hence, $\ker(W^\star W) = \ker(|T|)$. By Proposition \ref{II.B.1},
$\ran(W^\star) = \overline{\ran(|T|)}^{\wks}$ and $\ker(W) = \ker(|T|)$. This is half of what we need to prove.
Proceeding to the other half, we have that
\[
	\ran(T) = \ran(W|T|) = \ran(WW^\star W|T|) \subset \ran(WW^\star),
\]
from which it follows that $\overline{\ran(T)}^{\wks} \subset \ran(WW^\star)$. On the other hand,
\[
	WW^\star = W(W^\star W)W^\star = \wlim_{n \to \infty} W|T|^{1/n}W^\star.
\]
Since $\ran(|T|^{1/n}) \subset \overline{\ran(|T|)}$,
we have $\ran(W|T|^{1/n}) \subset \overline{\ran(T)}$,
for all $n \in \mathbb{N}$.
Thus  $\ran(WW^\star) \subset \overline{\ran(T)}^{\wks}$. Now suppose
that $T^\star x = 0$. Then $|T|W^\star x = 0$, which
implies that $|T|^{1/n}W^\star x = 0$ for all $n \in \mathbb{N}$. Therefore, $WW^\star x = 0$. Conversely,
suppose that
$WW^\star x = 0$. Then $W^\star x = 0$. Thus, $T^\star x = |T|W^\star x = 0$. Hence, 
$\ker(WW^\star) = \ker(T^\star)$. By Proposition \ref{II.B.1} again, $\ran(W) = \overline{\ran(T)}^{\wks}$ and
$\ker(W^\star) = \ker(T^\star)$.
\end{proof}

\begin{corollary} \label{II.B.2a}
Let $J$ be a linear subspace of a dual operator space $X$. Then the following are equivalent:
\begin{enumerate}
\item[(i)] $J$ is a right $M$-summand of $X$.
\item[(ii)] $J = \overline{\ran(T)}^{\wks}$ for some $T \in \Al(X)$.
\item[(iii)] $J = \ker(T)$ for some $T \in \Al(X)$.
\end{enumerate} 
\end{corollary}

\begin{corollary} \label{II.B.3b}
Let $X$  be an operator space and $T \in \Al(X)$. Then $\overline{\ran(T)}$ is a
right $M$-ideal of $X$.
\end{corollary}

\begin{proof}  
By basic functional analysis we have $$\overline{\ran(T)}^{\perp \perp} 
= \ker(T^*)^{\perp} = ( ^\perp\ran(T^{**}))^{\perp}
=  \overline{\ran(T^{**})}^{\wks}.$$ 
We shall see in Section \ref{IV.C} that $T^{**} \in \Al(X^{**})$.
Thus the result follows from Corollary \ref{II.B.2a}.
\end{proof}

Dual to Theorem \ref{II.B.2} we have the following result:

\begin{corollary} \label{II.B.3}
Let $X$ be an operator space and $T \in \Cr(X)$. Then $T = W|T|$, where $W \in \Cr(X)$ is a partial isometry such
that
\[
	\ran(W^\star) = \overline{\ran(|T|)}, \ran(W) = \overline{\ran(T)}, \ker(W) = \ker(|T|), \text{ and }
		\ker(W^\star) = \ker(T^\star).
\]
In particular, $X = \overline{\ran(T)} \oplus_{\lL} \ker(T^\star)$.
\end{corollary}

The proof is essentially the same as that of Theorem \ref{II.B.2}.

\begin{corollary} \label{II.B.3a}
Let $J$ be a linear subspace of an operator space $X$. Then the following are equivalent:
\begin{enumerate}
\item[(i)] $J$ is a left $L$-summand of $X$.
\item[(ii)] $J = \overline{\ran(T)}$ for some $T \in \Cr(X)$.
\item[(iii)] $J = \ker(T)$ for some $T \in \Cr(X)$.
\end{enumerate}
\end{corollary}

{\bf Remarks.} 1) Corollaries \ref{II.B.1b}, \ref{II.B.2a}, and \ref{II.B.3b}
are not true in general if $T \in \Ml(X)$.  

2) We shall see later (see the last remark in
Section \ref{IV.C}) that left adjointable multipliers on a general  
operator space always have a `pseudo polar decomposition'.  

3) In view of relations such as those in 
Corollaries \ref{II.B.2a},  \ref{II.B.3b} and
  \ref{II.B.3a}, it is tempting to 
believe that the kernel of a map $T \in \Al(X)$
is a right $M$-ideal of $X$.  But this is not 
true. For a counterexample, let $T = P + Q$ for
two left $M$-projections $P, Q$ on $X$.  If 
$T x = 0$, then $$\langle P x + Q x , x \rangle = 
\langle P x , P x \rangle + \langle Q x , Q x  \rangle
= 0, $$
where $\langle \cdot , \cdot \rangle$ is the 
inner product in Section \ref{II.BBB}.  
This implies that $P x = Q x = 0$. Thus 
$\ker(P + Q) = \ker(P) \cap \ker(Q)$.  However, we shall
see in Example \ref{IV.D.6} that the intersection of 
right $M$-summands need not be a right $M$-ideal.   

4) It is not true that every right $M$-ideal is the closure of the
range of a left adjointable multiplier. For example, $K(H)$ is a right $M$-ideal
of $\MIN(B(H))$, but $\Al(\MIN(B(H))) = \mathbb{C}\Id_{B(H)}$
(see Subsection \ref{III.B}).
  
\subsection{One-Sided $M$-summands and Relative
Multipliers}

A recent paper \cite{BSur} introduces
 one-sided multipliers between two different spaces.
There are several connections between this theory
and one-sided $M$-ideals.  For example, as shown there
it is natural to define such multipliers $T : X \to Y$ 
in the case that $X, Y$ are complementary  one-sided $M$-summands
in a third space $V$.  Or, one could  
attempt to generalize to the new class
of multipliers the results in the last Subsection  
and their applications to one-sided $M$-ideals.  We will
just list one result here, a sufficient condition
for a subspace to be a one-sided $M$-summand.  In fact  
the following generalizes a known characterization
of adjointable projections on $C^*$-modules:
 
\begin{lemma} \label{auor}   Suppose that $X$ is a subspace of
an operator space $Y$, and that there is a completely
contractive projection $P$ from $Y$ onto $X$.
We suppose that we have fixed
Hamana triple envelopes ${\mathcal T}(X)$ and ${\mathcal T}(Y)$
of $X$ and $Y$ respectively,
which are right $C^*$-modules over the same algebra
${\mathcal B}$.  If $P$ and the inclusion map
$\epsilon : X \to Y$ are
relative left multipliers in the sense of \cite{BSur},
then $\epsilon \circ P$ is a left $M$-projection on $Y$,
$X$ is a right $M$-summand of $Y$, and $P, \epsilon$ are
relatively adjointable, with $P^\star = \epsilon$.
\end{lemma}
 
\begin{proof}   By the definition in \cite{BSur}
there are contractive right 
${\mathcal B}$-module maps $\tilde{P} : {\mathcal T}(Y) \to
{\mathcal T}(X)$, and $\tilde{\epsilon} : {\mathcal T}(X) \to
{\mathcal T}(Y)$ extending $P$ and $\epsilon$ respectively.
For $x \in X$ and $f \in {\mathcal B}$ we have
$$\tilde{P}(\tilde{\epsilon}(x b)) \; = \; \tilde{P}(\tilde{\epsilon}(x))
b \; = \; x b .$$
Thus $\tilde{P} \circ \tilde{\epsilon}$ is the identity map
on ${\mathcal T}(X)$.  From the theory of $C^*$-modules
it follows that
$\tilde{\epsilon}$ and $\tilde{P}$ are adjointable and that
$\tilde{\epsilon}^\star = \tilde{P}$.   Thus
$\tilde{\epsilon} \circ \tilde{P}$ is an
adjointable projection on ${\mathcal T}(Y)$, and so its
restriction to $Y$ is an orthogonal projection in
$\Al(Y)$.   The other assertions are evident, except for the
fact that $P^\star = \epsilon$, which follows from
the appropriate definition in \cite{BSur}.
\end{proof}
  
\section{Examples} \label{III}

Many examples of one-sided $M$-structure of particular 
operator spaces are listed in our previous papers
\cite{BEZ,BSZ}.   We will need to restate a few main facts 
here, omitting proofs given in those papers to avoid repetition.   
Because this section is devoted
to {\em examples}, we will feel free to quote from
results proved only later in this paper.     

\subsection{Two-Dimensional Operator Spaces} \label{III.A}

In this section we consider
the one-sided $M$-structure of two-dimensional operator spaces. 
First, however,  we give a result for
finite-dimensional operator spaces.
 
\begin{proposition} \label{fd}
Let $X$ be a $n$-dimensional operator space. Then $\Al(X)$ is 
$*$-isomorphic to a unital *-subalgebra of $M_m$, for some $m \leq n$.
\end{proposition}
 
\begin{proof}
Since $\Al(X) \subset B(X)$, $\Al(X)$ is finite-dimensional, and is
therefore *-isomorphic to
$M_{n_1} \oplus M_{n_2} \oplus ... \oplus M_{n_k}$ for some positive integers $n_1, n_2, ..., n_k$.
Since $\dim(X) = n$, $\Al(X)$ has at most $n$ mutually orthogonal nonzero projections. Hence,
$n_1 + n_2 + ... + n_k \leq n$.
\end{proof}

{\bf Remark.}  Since  $\Ml(X) \subset B(X)$,
it is clear that $\Ml(X)$ is also finite-dimensional
for an $n$-dimensional operator space $X$. However,
there may exist no finite 
$m$ with $\Ml(X) \subset M_m$ completely isometrically.
This may be seen via the following well-known trick:
take an finite-dimensional  
operator space $Y \subset B(H)$, and let $X$ be the unital operator algebra
in $M_2(B(H))$ with `$Y$ as the $1-2$ corner', and `scalars on the
diagonal'.  Then $Y \subset X = \Ml(X)$
completely isometrically.  If $Y$ is not an exact operator space
\cite{Pisier} then the assertion follows. \\

From this point forward in this subsection, 
$X$ will be a two-dimensional operator space, so that the nontrivial right $M$-summands
correspond to certain unit vectors $x \in X$. 
We remark in passing that if $J$ is a right $M$-summand 
or $M$-ideal of a general
operator space $X$, if $x$ is a unit 
vector in $J$, and if $y$ is a unit vector in the complementary piece
$J^\rightthreetimes$, then by the last remark in Section \ref{II.BBB},
one may easily see that
$J_0 = \spn\{x\}$ is a right
$M$-summand in $X_0 = \spn \{x , y\}$. Thus, in
some sense, the 2-dimensional case
captures at least some of the
geometry of general right $M$-ideals. 
It would be nice to have a tidy geometric characterization 
of right $M$-summands of two-dimensional operator spaces,
however this still eludes us. To appreciate some of the
complications, note that surprisingly even $\ell^p_2$ 
(with some operator space structure) can have nontrivial
right $M$-summands (see Example \ref{exin}). We can at least
say the following, 
which rules out the existence 
of right $M$-summands in some classical spaces:  

\begin{proposition} \label{III.A.1}
If $X$ is a two-dimensional operator space with a nontrivial right
$M$-summand, then
\begin{enumerate}
\item[(i)] $X$ is completely isometric to the span in $B(H)$ of two operators $S$ and $T$ of norm one such that 
$S^*T = 0$.
\item[(ii)] There is a completely contractive surjection $\rho : C_2 \to X$
taking $e_1$ and $e_2$ to unit vectors in the 
complementary summands.
\end{enumerate}
\end{proposition}

\begin{proof}
Let $x$ be a unit vector in the nontrivial right $M$-summand $J$ of $X$, and $y$ be a unit vector in the complementary
summand. Let $P \in \Al(X)$ be the left $M$-projection onto $J$. Then there exists a complete isometry $\sigma:X \to B(H)$
such that $\sigma(Pv)^*\sigma(w) = \sigma(v)^*\sigma(Pw)$ for all $v, w \in X$.
Let $S = \sigma(x)$ and $T = \sigma(y)$.
Then
\[
        S^*T = \sigma(x)^*\sigma(y) = \sigma(Px)^*\sigma(y) = \sigma(x)^*\sigma(Py) = 0.
\]
This proves (i). Now define $\rho:C_2 \to \spn\{S, T\}$ by
\[
        \rho\left(\begin{bmatrix} a\\ b \end{bmatrix}\right) = aS + bT
\]
for all $a, b \in \mathbb{C}$. For $A, B \in M_n$, one has that
\begin{eqnarray*}
        \|A \otimes S + B \otimes T\|^2 &=& \|A^*A \otimes S^*S + A^*B \otimes S^*T + B^*A \otimes T^*S +
                B^*B \otimes T^*T\|\\
        &=& \|A^*A \otimes S^*S + B^*B \otimes T^*T\| \leq \|(A^*A + B^*B) \otimes I\|\\
        &=& \|A^*A + B^*B\| = \left\|\begin{bmatrix} A\\ B 
\end{bmatrix}\right\|^2.
\end{eqnarray*}
This proves (ii).
\end{proof}

We can say a little more about the multiplier algebras
of two-dimensional spaces.  
 
\begin{proposition} \label{III.A.2}
If $X$ is a two-dimensional operator space, then $\Al(X)$ is either
$\mathbb{C}$, $\ell_2^\infty = \mathbb{C} \oplus \mathbb{C}$, or $M_2$.
\begin{itemize}
\item [(i)] $\dim \Al(X) = 4$ if and only if $X \cong C_2$ completely
isometrically,
\item [(ii)] $\dim \Al(X) = 2$ and $\dim \Ar(X) = 2$,
if and only if $X \cong \ell^\infty_2$ completely
isometrically,
\item [(iii)] $X$ admits a nontrivial left $M$-projection $P$
such that $P^*$ is a right $M$-projection if and only if $X \cong C_2$.
\end{itemize}
\end{proposition}
 
\begin{proof}  
The first assertion follows from Proposition \ref{fd}.

(i) This follows from Theorem \ref{IV.J.3}.
 
(ii)  Suppose that $P$ is a nontrivial  left $M$-projection
and that $Q$ is a  nontrivial complete right $M$-projection on $X$.
Then $P Q = Q P$ as maps on $X$. By elementary linear
algebra it follows that $Q = P$ or $Q = \Id - P$.
Thus $P$ is both a left and a right
$M$-projection, and hence is a complete (two-sided) $M$-projection. That is,
$X \cong P(X) \oplus^\infty (\Id - P)(X) \cong \ell^\infty_2$ completely isometrically.

(iii) Suppose that $P$ is a nontrivial left $M$-projection on $X$ such that $P^*$ is a right $M$-projection
on $X^*$. Let $x, y \in X$ be unit vectors such that $Px = x$ and $Py = 0$.
By Proposition \ref{III.A.1}, the map $T:C_2 \to X:\begin{bmatrix} a\\ b \end{bmatrix} \mapsto ax + by$
is a surjective complete contraction. Now let $\phi, \psi \in \ball(X^*)$ with $\phi(x) = \psi(y) = 1$.
Set $f = P^*(\phi)$ and $g = (\Id - P^*)(\psi)$. Then $P^*f = f$, $f(x) = 1$, and $f(y) = 0$. Likewise,
$P^*g = 0$, $g(x) = 0$, and $g(y) = 1$. By the `other-handed' version of Proposition \ref{III.A.1}, the map
$S:R_2 \to X^*:\begin{bmatrix} \alpha & \beta \end{bmatrix} \mapsto \alpha f + \beta g$ is a surjective
complete contraction. Since
\[
	\langle S^*(ax + by), \begin{bmatrix} \alpha & \beta \end{bmatrix} \rangle
	= \langle ax + by, \alpha f + \beta g \rangle = a\alpha + b\beta
	= \left\langle \begin{bmatrix} a\\ b \end{bmatrix}, \begin{bmatrix} \alpha & \beta \end{bmatrix} \right\rangle,
\]
we have $S^* = T^{-1}$, and so $T$ is a complete isometry.
\end{proof}

\subsection{$\MIN$ and $\MAX$ Spaces} \label{III.B}

We recall the two most common 
methods to regard a Banach space $X$ as an operator 
space: namely $\MIN(X)$ and $\MAX(X)$.      
The following facts about a 
Banach space $X$ may be found in \cite[Section 6.1]{BEZ}.
Firstly, the left $M$-projections on $\MIN(X)$ are exactly the classical $M$-projections on $X$.
Therefore, the right $M$-summands and right
$M$-ideals of $\MIN(X)$ are just the $M$-summands and $M$-ideals of $X$, respectively. Likewise,
the right $L$-projections on $\MAX(X)$ are the $L$-projections on $X$, which implies that
the left $L$-summands of $\MAX(X)$ are the $L$-summands of $X$.
 
\subsection{Hilbertian Operator Spaces} \label{III.C}

For a Hilbertian operator space $H$, 
every left $M$-projection on $H$ is an orthogonal projection on $H$.
Consequently, if $J$ is a right $M$-summand of $H$, then its orthogonal complement and complementary right $M$-summand
are one and the same.

Recall that an operator space $X$ is \emph{homogeneous} if $B(X) = CB(X)$ isometrically. This implies that every surjective
isometry is a complete isometry, and for a Hilbertian operator space, the latter condition is equivalent to homogeneity
(see \cite{Pisier}).

\begin{theorem} \label{III.C.1}
Let $H$ be a homogeneous Hilbertian operator space, with $\dim(H) \geq 2$. Then the following are equivalent:
\begin{enumerate}
\item[(i)] $H$ is a Hilbert column space.
\item[(ii)] $H$ has a nontrivial right $M$-summand.
\item[(iii)] $H$ has a nontrivial left-adjointable multiplier (i.e. $\Al(H) \neq \mathbb{C}$).
\end{enumerate}
\end{theorem}

\begin{proof}
Clearly (i) $\Rightarrow$ (ii) $\Rightarrow$ (iii). 

(iii) $\Rightarrow$ (ii)  Since Banach space reflexivity
implies operator space reflexivity, $H$ is a dual operator space, and $\Al(H)$ is a $W^*$-algebra by Theorem
\ref{wst}. Therefore, $\Al(H)$ is the norm-closed linear span of its orthogonal projections (the left $M$-projections on
$H$). Thus, if $H$ has no nontrivial right $M$-summands, then $\Al(H) = \mathbb{C}$. 

(ii) $\Rightarrow$ (i) Suppose $J$ 
is a nontrivial right $M$-summand of $H$, with corresponding left $M$-projection $P$. Let $U$ be a unitary operator on $H$. 
Then $U$ is a completely isometric surjection, and $UPU^{-1}$ is a left $M$-projection on $H$, as discussed at the end of 
Section \ref{I.B}. Thus, every closed linear subspace of $H$ of the same (Hilbertian) dimension as $J$ is a 
right $M$-summand of $H$.
Now let $x \in J$ and $y \in J^\perp$ be unit vectors. Define $K = J \cap \{x\}^\perp$. Then $\spn\{y\} \vee K$ has the
same dimension as $\spn\{x\} \vee K = J$. Thus, $\spn\{y\} \vee K$ is a right $M$-summand of $H$. 
By Proposition \ref{IV.D.1}, $K = J \cap (\spn\{y\} \vee K)$ is a right $M$-summand of $H$, which implies that $K^\perp$ is 
as well. Again by Proposition \ref{IV.D.1}, $\spn\{x\} = J \cap K^\perp$ is a right $M$-summand of $H$. Arguing as above, 
every one-dimensional linear subspace of $H$ is a right $M$-summand of $H$. Hence, every rank-one orthogonal projection 
on $H$ is an element of $\Al(H)$. By polarization, every rank-one operator is an element of $\Al(H)$. 
Therefore, if $\{e_i: i \in I\}$ is an orthonormal basis for $H$, then 
$\{e_i \otimes \overline{e}_i: i \in I\}$ is a family of mutually 
orthogonal and equivalent left $M$-projections which add to the identity in the point-norm topology. By Theorem 
\ref{IV.J.3}, we have $H \cong C_I \cong H_c$.
\end{proof}

\begin{theorem} \label{III.C.2}
If $X$ is an operator space, then the following are equivalent:
\begin{enumerate}
\item[(i)] $X$ is completely isometric to a Hilbert column space.
\item[(ii)] $B(X) = \Al(X)$. That is, every bounded operator on $X$ is a left adjointable multiplier.
\item[(iii)] $CB(X) = \Al(X)$. That is, every completely bounded operator on $X$ is a left adjointable multiplier.
\item[(iv)] $B(X) = \Ml(X)$ isometrically.
\item[(v)] $CB(X) = \Ml(X)$ isometrically.
\item[(vi)] Every one-dimensional subspace of $X$ is a right $M$-summand.
\end{enumerate}
\end{theorem}

\begin{proof}
(i) $\Rightarrow$ (ii) See \cite{Shilov} if necessary.

(ii) $\Rightarrow$ (iii) We have the inclusions 
$\Al(X) \subset CB(X) \subset B(X) = \Al(X)$. 

(iii) $\Rightarrow$ (v) We have the isometric
inclusion $CB(X) = \Al(X) \subset \Ml(X)$ as well as the contractive inclusion $\Ml(X) \subset CB(X)$. It follows
that $CB(X) = \Ml(X)$ isometrically. 

(v) $\Rightarrow$ (vi) Let $x$ be a unit vector in $X$, and
$\psi \in X^*$ be such that $\|\psi\| = 1$ and $\psi(x) = 1$. Define $P:X \to X$ by $P(y) = \psi(y)x$ for all $y \in X$.
Then $P$ is a completely contractive projection. Since $CB(X) = \Ml(X)$ by assumption, $P$ is a left multiplier of
$X$ and $\|P\|_{\Ml(X)} \leq 1$. It follows that $P$ is a left $M$-projection on $X$, and that $\spn\{x\}$ is a right 
$M$-summand of $X$. 

(vi) $\Rightarrow$ (i) Suppose that
every one-dimensional subspace of $X$ is a 
right $M$-summand of $X$.
First suppose that $\dim(X) \geq 3$.
Then every finite-dimensional 
subspace of $X$ is a right 
$M$-summand of $X$, by Proposition \ref{IV.D.2}. In particular,
every two-dimensional subspace of $X$ is 
contractively complemented. By a result of Kakutani and Bohnenblust 
(\cite{Kakutani, Bohnenblust}), $X$ is Hilbertian. 
Arguing as in the conclusion of the proof of Theorem \ref{III.C.1},
we see that $X = C_I$ completely isometrically, where $I$ is any index set corresponding to an orthonormal basis
for $X$.   In the case that $\dim(X) = 2$ we first observe 
that (vi) implies that there exist noncommuting left $M$-projections
on $X$. Hence $\dim(\Al(X)) = 4$.  
Proposition \ref{III.A.2} now gives $X = C_2$.     

(ii) $\Rightarrow$ (iv) $\Rightarrow$ (vi)   This is
nearly identical to the proof of
(iii) $\Rightarrow$ (v) $\Rightarrow$ (vi), and is left to the reader.
\end{proof}

{\bf Remark.}  
We do not know
 whether the condition $B(X) = \Ml(X)$ (as sets) implies
that $X$ is Hilbert column space completely isometrically.
If every bounded map on $X$ is a left multiplier,
then the open mapping theorem implies that $B(X)$ is isomorphic to
the operator algebra $\Ml(X)$. It is well-known that this implies 
that $X$ is isomorphic to a Hilbert space (see \cite{Ble_cb}).
Also, it implies that $X$ is $\lambda$-homogeneous for some $\lambda$.
The latter conditions, while necessary, are not sufficient. For example,
$R_n$ ($n$-dimensional row space) is a homogeneous Hilbertian operator space,
yet $\Ml(X) = \mathbb{C} \neq B(X)$.

\subsection{$C^*$-Algebras} \label{III.D}

Let $\mathcal{A}$ be a $C^*$-algebra. As observed in \cite{BEZ}, \S 1, the left $M$-projections on 
$\mathcal{A}$ are precisely the maps
\[
	P:\mathcal{A} \to \mathcal{A}:T \mapsto ET,
\]
where $E$ is an orthogonal projection in $\M(\mathcal{A})$, the multiplier algebra of $\mathcal{A}$ 
(see \cite{Wegge}, Chapter 2, for the definition of this algebra). More generally,
$\Al(\mathcal{A}) = \M(\mathcal{A})$ and $\Ml(\mathcal{A}) = \LM(\mathcal{A})$, the left multiplier algebra
of $\mathcal{A}$ (cf. \cite{Wegge}, Exercise 2.F). 
It follows that the right $M$-ideals of $\mathcal{A}$
are simply the closed right ideals. On the other hand, $\mathcal{A}$ has trivial left $L$-structure
(\cite{BSZ}, Theorem 5.1). For further information on the one-sided $M$- and $L$-structure of $C^*$-algebras,
see \cite{BSZ}. See \cite{Effros,Prosser,Akemann70,RN} for 
 studies of the one-sided ideal structure of $C^*$-algebras, and the connection with
important projections in the second dual.
  
\subsection{Nonselfadjoint Operator Algebras} \label{III.E}

Let $\mathcal{B}$ be an operator algebra with a two-sided
contractive approximate identity.  In fact
many of the results below only require 
a one-sided contractive approximate identity, but for simplicity
we will stick to the two-sided case. 
It is shown in \cite{Shilov} that $\Ml(\mathcal{B})$ is the usual 
operator algebra left multiplier algebra
$\LM(\mathcal{B})$, and thus $\Al(\mathcal{B}) = \{ T \in \LM(\mathcal{B}) : 
T^* \in \LM(\mathcal{B}) \}$.
Consequently, the left $M$-projections 
on $\mathcal{B}$ are simply the maps of left multiplication
by an orthogonal projection $e$ in $\LM(\mathcal{B})$, 
the right $M$-summands of $\mathcal{B}$ are the 
principal right ideals $e \mathcal{B}$ for such $e$,    
and the right $M$-ideals of $\mathcal{B}$ are the
closed right ideals having a contractive left 
approximate identity (\cite{BEZ}, Proposition 6.4). 
One can find additional 
results in \cite{BSZ,Dual,OpAlg} for example.
We will content ourselves here with two new results.
They represent typical applications of the
results proved in the present paper. Note that the first application 
may be described
as approximation result, since approximate identities 
are often used to approximate elements:

\begin{corollary} \label{III.E.1}  Let $\mathcal{B}$ be an
operator algebra with a two-sided
contractive approximate identity.  The closed span
in $\mathcal{B}$ of a family of closed
right ideals, each possessing a left
contractive approximate identity for itself,
is a right ideal also possessing a left
contractive approximate identity.
\end{corollary}    

\begin{proof} This follows from remarks above together with
Corollary \ref{IV.D.2}. \end{proof}

\textbf{Remark.}  
We are grateful to George Willis for some examples that
show that the analogous result fails quite badly
for Banach algebras.  In fact it is not hard to find three
or four dimensional counterexamples.   
  
\begin{corollary} \label{III.E.2}  Let $\mathcal{B}$ be an
operator algebra with a two-sided
contractive approximate identity, and fix an element $a$ in the 
diagonal $\mathcal{B} \cap \mathcal{B}^*$ (or in the diagonal of 
$\LM(\mathcal{B})$). Let $J = a\mathcal{B} = 
\{ a b : b \in \mathcal{B} \}$, a principal right ideal.
\begin{itemize}
\item [(i)]   The norm closure of  $J$ is a right ideal possessing a  
left contractive approximate identity. 
\item [(ii)] If  $J$ is already
closed, then $J =  e \mathcal{B}$, for an orthogonal projection 
$e \in \LM(\mathcal{B})$.  In fact,
the projection $e$ may be chosen in $J$.        
\end{itemize} 
\end{corollary}  

\begin{proof} Note that $J = \ran(T)$ where 
$T$ is the map of left multiplication by $a$ on $B$.
This map is adjointable by  remarks above.

(i)  Appeal to Corollary \ref{II.B.3b}.  Alternatively, this
may be argued using ideals in the diagonal $C^*$-algebra.  

(ii) 
By Corollary \ref{II.B.1b} applied to the map $Tb = ab$, and the remarks above,
we have $J =  e \mathcal{B}$, for an orthogonal projection
$e \in \LM(\mathcal{B})$. If $\mathcal{B}$ is unital then we are done,
and note that in this case $J$ has a left identity of norm one.
If not we apply (ii) in the unital case but with $\RM(\mathcal{B})$ in
place of $\mathcal{B}$.  
Note that $J \subset a \RM(\mathcal{B})$ clearly.
For $T \in \RM(\mathcal{B})$, and $(e_t)$ a contractive approximate identity
of $\mathcal{B}$, we have 
$a T = \lim_t a e_t T \in J$.  Thus $J = a \RM(\mathcal{B})$.
We deduce from the unital case that $J$ contains a left identity 
$e$ of norm one. Thus $J = e \mathcal{B}$. \end{proof}

\subsection{Hilbert $C^*$-Modules} \label{III.F}
 
Let $X$ be a right Hilbert $C^*$-module (see e.g. \cite{La}).
Then the right $M$-ideals of $X$ are the closed right
submodules of $X$, and the right $M$-summands of $X$ are the 
orthogonally complemented right
submodules of $X$. For the details, see \cite{BEZ}, 
where we also give one or two applications of one-sided $M$-ideal 
theory to $C^*$-modules. We will not explicitly list more    
such applications in the present paper. In fact, the 
flow is mostly in the other direction: most of the results in the present
paper may be seen as {\em generalizations} of known results 
about $C^*$-modules. We invite the reader
to glance once through our results with this in mind;
the parallels are quite striking.  This perspective is 
highlighted in \cite{BSur}.

It is interesting that via Theorems \ref{BEZ} and \ref{Hermitian},
one has a purely
operator space characterization of bounded module maps/adjointable
maps  on a $C^*$-module. In contrast to 
the general operator space case, we have been able to find a similar 
characterization
for such maps $T : W \to Z$ between two {\em different} 
$C^*$-modules. Since we will not need this here,
we will present the details elsewhere.

\subsection{Banach Algebras} \label{III.G}

A useful and striking
result due to Smith and Ward \cite{SW} states that the 
classical $M$-ideals in a unital Banach algebra  $\mathcal{A}$ are necessarily 
subalgebras of $\mathcal{A}$.   This suggests that the 
one-sided $M$-ideals in a unital `completely contractive Banach algebra'
 are also subalgebras.   We have not been able to prove or disprove 
this. We do have the following
preliminary step, which we will need elsewhere:  

\begin{lemma} \label{III.G.1}
Let $\mathcal{A}$ be a unital Banach algebra which is also an operator space. If $T \in \Al(\mathcal{A})_{sa}$, then
$T(1) \in \Her(\mathcal{A})$.
\end{lemma}

\begin{proof}
Without loss of generality, we may assume that $0 \leq T \leq I$. Then we have
\[
	\|T(1) + it1\| = \|(T + itI)(1)\| \leq \|T + itI\| \leq \sqrt{1 + t^2}
\]
for all $t \in \mathbb{R}$. Thus, by I.10.10 in \cite{BonsallDuncan}, we may deduce that $T(1) \in \Her(\mathcal{A})$.
\end{proof}

\subsection{Operator Modules}

We will not define operator modules here, but simply give the 
equivalent characterization discovered in \cite{Shilov}
(see \cite{BP01} for an alternative account): they correspond
to the module actions given by the formula $a x = \pi(a)(x)$,
for $a \in {\mathcal A}, x \in X$,
 for a completely contractive unital homomorphism $\pi : {\mathcal A}
 \to \Ml(X)$.
Here ${\mathcal A}$ is an operator space and a unital algebra. 
If ${\mathcal A}$ is a $C^*$-algebra we may replace $\Ml(X)$ by $\Al(X)$ here.
Thus for example
from a result in \cite{BSZ} one sees that there are no
nontrivial (i.e.\ nonscalar) operator module actions
of a $C^*$-algebra on an $L^1$ space, or  
more generally on the predual of a von Neumann algebra.  This
is because for such spaces all left adjointable maps
are trivial (i.e. a scalar multiple of the identity).
 On the other hand such `$L^1$ spaces' $X$ have rich one-sided
Cunningham 
algebras, and thus $X$ is a nontrivial bimodule over
$\Cl(X)$ and $\Cr(X)$.     
  
\begin{proposition} \label{Msubm}  If $J$ is a left $M$-ideal in an operator
space $X$, then $T(J) \subset J$ for all $T \in \Ml(X)$.
Thus if $X$ is a left operator module over ${\mathcal A}$ then
$J$ is an ${\mathcal A}$-submodule of $X$.
\end{proposition}
 
\begin{proof}
Let $P$ be the right $M$-projection onto $J^{\perp \perp}$.
Inside $X^{**}$ we have $$T(J) = T^{**}(J) \subset T^{**}(P(X^{**}))
= P(T^{**}(X^{**})) \subset J^{\perp \perp} ,$$
using the fact from Subsection \ref{IV.C} that $T^{**} \in \Ml(X^{**})$, and the
last point in Subsection \ref{I.A}.  Thus
$T(J) \subset J^{\perp \perp} \cap X = J$, by
basic functional analysis.
\end{proof}               

From this Proposition and 
the situation in  Subsection \ref{III.F},
one might expect that every submodule of a left operator module 
$X$ is a left $M$-ideal.   However this is not true in general.
In well behaved situations like when $X$ is an  operator algebra 
with contractive approximate identity
one at least needs an extra
condition (see Subsection \ref{III.E}).  Another good example of this
is furnished by the strong Morita equivalence 
${\mathcal A}$-${\mathcal B}$-bimodules of
\cite{BMP}.  These are the analogues for nonselfadjoint operator
algebras of the $C^*$-algebraic strong Morita equivalence bimodules of
Rieffel (see \cite{La} for example),
or equivalently the `full' Hilbert $C^*$-modules 
discussed in Subsection \ref{III.F}.
By analogy with the $C^*$-module case one would expect that
if $X$ is a strong Morita equivalence 
${\mathcal A}$-${\mathcal B}$-bimodule in the sense of
\cite{BMP}, then the
right $M$-ideals of $X$ are exactly the closed right ${\mathcal B}$-submodules.
However in fact one needs to impose further conditions, and even
then the situation is quite complicated---see the discussion 
in \cite[Section 10]{BSur}.  
The  right $M$-summands of such a bimodule are easier to deal with,
they are just the ranges of completely contractive 
idempotent maps on $X$ which are right ${\mathcal B}$-module maps.
That is, they are exactly the analogue in this situation
of the `orthogonally complemented' 
submodules.  The proof of this 
right $M$-summand case (from joint work of the first author with 
Solel from '00) may be found in
\cite[Section 10]{BSur}; it proceeds by first establishing
that the `Hamana triple envelope' or `noncommutative
Shilov boundary' of $X$ is a canonical  
$C^*$-algebraic strong Morita equivalence bimodule 
over the $C^*$-envelopes of ${\mathcal A}$ and ${\mathcal B}$.
 
\subsection{Operator Systems and $M$-Projective Units} \label{Osy}
 
If $X$ is an operator system, then the left $M$-projections are just the `involutions' of
the right $M$-projections, where the `involution' of 
a map $T$ on $X$ is the map $x \mapsto T(x^*)^*$.  
Indeed, it is clear from the definition
of the multiplier algebras in terms of $C_e^*(X)$
(see e.g.\ \cite{Shilov}, or \cite[p. 505--506]{Dual})
that $\Al(X)$ and $\Ar(X)$ may be viewed as 
{\em subsets} of $X$, and in fact these two
subsets of $X$ coincide.
We are not saying that every left adjointable multiplier
is right adjointable, this would imply 
that every closed left ideal of a unital $C^*$-algebra is
a closed right ideal! In fact, the left $M$-structure 
of an operator system is related to the right
$M$-structure by the above-mentioned involution.  Since this
will not play a role for us, we leave the details to the
reader.

We remark in passing that for a certain  class of 
operator systems, multipliers were investigated 
by Kirchberg \cite{Kirchberg}.  This is related 
to his study of the sum $J + K$ of a 
left and a right ideal in a $C^*$-algebra,
which was important in his profound work on nuclear
$C^*$-algebras, etc.   

We say that an element $x$ in
an operator system $X$ is a {\em left $M$-projective unit}, 
if $x = Pe$ for a  left $M$-projection $P$ (here $e$ is the unit in $X$).
We write $X^+_1$ for the
positive part of the unit ball of $X$.
 
\begin{corollary}  \label{kw} Let $X$ be an
operator system which is also a dual operator space.
Then $X$ has a product with respect to which
it is a $W^*$-algebra (with the same underlying operator system structure), 
if and only if
every extreme point of $X^+_1$ is a left $M$-projective unit in $X$. In this case,
such a product is unique.
\end{corollary}
 
\begin{proof}  The uniqueness of the product may be deduced 
from any operator space variant of the Banach-Stone theorem
(see e.g.\ \cite[Corollary 5.2.3]{ERbook}).

($\Rightarrow$)   If $X$ is a $W^*$-algebra, 
then its left  $M$-projective units are
exactly its orthogonal projections,
by facts in Section \ref{III.D}.
Thus the result follows from the characterization
of extreme points in the positive 
part of the unit ball of a
$W^*$-algebra \cite[Proposition 1.6.2]{Sakaib}.

($\Leftarrow$)
Consider the map $\rho : \Al(X) \to X : T \mapsto Te$.   
It is explained above Theorem 1.9 in 
\cite{Dual} that $\rho$ is a complete isometry, and 
it is easy to see that $\rho$ is weak-* continuous.
Hence, by the Krein-Smulian Theorem \ref{KS}, we
see that $\rho$ is a
weak-* homeomorphism with weak-* closed range.
By hypothesis, $\ran(\rho)$ contains all 
extreme point of $X^+_1$. By the Krein-Milman
Theorem, $\ran(\rho)$ contains $X^+_1$, and hence
$\rho$ maps onto $X$.
\end{proof}                                                

This result was inspired by, and is
closely related
to, an old theorem of K. H. Werner \cite[Corollary 5.1]{kw}.

\subsection{Locally Reflexive Operator Spaces} \label{III.I}
For results on the one-sided $M$-structure of 
locally reflexive operator spaces, we refer the reader to a paper
in preparation by the second author \cite{ZarikianMpaper}.
  
\section{Constructions} \label{IV}

In this section we develop the heart of the calculus of one-sided $M$-ideals. 
It is quite interesting that there
are perfect analogs of approximately half the results 
in the classical $M$-ideal calculus. For the other half,
some additional and natural hypotheses need to be imposed.

It would take up too much space to thoroughly define and introduce 
every one of the constructions listed in this section.  Readers 
interested in a particular construction should consult the 
basic texts on operator spaces, or the original papers, for 
precise definitions and basic facts.

\subsection{Opposite and Conjugate} \label{IV.A}

If $X$ is an operator space, then $X^{\op}$ is defined to 
be the same Banach space, but with `transposed' matrix
norms: $\left\|\begin{bmatrix} x_{ij} \end{bmatrix}\right\|_{\op} \equiv 
\left\|\begin{bmatrix} x_{ji} \end{bmatrix}\right\|$. 
Then $X^{\op}$ is an operator space, as is well known and 
easy to see.  If $\sigma : X \to {\mathcal A}$ is a complete isometry into
a $C^*$-algebra, then we obtain a matching 
complete isometry $\sigma^{\op} : X^{\op} \to {\mathcal A}^{\op}$.  The
 latter space is a $C^*$-algebra, namely ${\mathcal A}$ with its reversed
multiplication.   From this it is easy but tedious to check that 
$\Ml(X^{\op}) \cong \Mr(X)$ and $\Al(X^{\op}) \cong \Ar(X)$ in a natural 
way.  
We leave the precise assertion and its proof to the reader, since we
shall not need it.  It is also easy to 
see that if $J$ is a right $M$-ideal or right $M$-summand of $X$, then $J$ is a 
left $M$-ideal or left $M$-summand 
of $X^{\op}$. Likewise, left $L$-summands of $X$ correspond to right 
$L$-summands of $X^{\op}$.

Similar assertions hold for the {\em conjugate operator space}
 $\overline{X}$, which is  defined to 
be the set $\{\overline{x}: x \in X\}$
together with the linear structure defined by
\[
	\overline{x} + \overline{y} \equiv \overline{x + y} \text{ and } \alpha\overline{x} \equiv 
		\overline{\overline{\alpha}x},
\]
and the matrix norms
\[
	\left\|\begin{bmatrix} \overline{x}_{ij} \end{bmatrix}\right\| \equiv \left\|\begin{bmatrix} x_{ij} 
		\end{bmatrix}\right\|.
\]
Clearly, $\overline{X}$ is an operator space. If $J$ is a right $M$-ideal or right $M$-summand of $X$, then
$\overline{J} \equiv \{\overline{x}: x \in J\}$ is a right $M$-ideal or right $M$-summand of $\overline{X}$.
Likewise, the left $L$-summands of $X$ correspond (via `conjugation') to the left $L$-summands of $\overline{X}$.

We say that an operator space is {\em symmetric } if $X = X^{{\rm op}}$.
  For a symmetric operator space it is easy to see that $T \in \Ml(X)$ 
if and only if $T \in \Mr(X)$,
and similarly for $\Al(X)$ and $\Ar(X)$.  It follows from 
the last fact in Subsection \ref{I.A} that $\Al(X)$ and 
$\Ar(X)$ are commutative.  In fact it follows from Subsection
\ref{S71} that $\Al(X) = \Ar(X) = Z(X)$, with the notation
from that Subsection.  Thus from that Subsection
we see that the left $M$-projections (resp. right 
$M$-summands, right
$M$-ideals) on $X$
are precisely the `complete $M$-projections' (resp.
complete $M$-summands, complete $M$-ideals).
   
An example of this is the 
$O\ell^p$ of Pisier \cite{Pisier}.
  These are formed by interpolating between the
symmetric spaces $O\ell^1 = c_0^*$ and $\ell^\infty$.
Thus they are symmetric.  In this case there are 
no nontrivial classical $M$-ideals or $M$-summands if $p < \infty$,
hence no nontrivial  complete $M$-ideals or complete $M$-projections,
and hence no nontrivial one-sided $M$-ideals or projections.
Since in this case  $\Al(O\ell^p)$ is
a $W^*$-algebra and therefore generated by its projections,
we must have that $\Al(O\ell^p)$ is one dimensional.

\subsection{Subspace and Quotient} \label{IV.B}

To facilitate our development, we introduce two operators associated with an operator and 
an invariant subspace:

\begin{definition} \label{IV.B.1}
Let $X$ be an operator space, $Y$ be a closed linear subspace of $X$, and $T:X \to X$ be a linear map such that 
$T(Y) \subset Y$. By $T|_Y:Y \to Y$ we denote the restriction of $T$ to $Y$, and by $T/Y:X/Y \to X/Y$ we denote 
the (well-defined) linear map
\[
        (T/Y)(x + Y) = Tx + Y.
\]
\end{definition}

We observe that if $T$ is completely bounded, then so are $T|_Y$ and $T/Y$, and 
$\|T|_Y\|_{cb}, \|T/Y\|_{cb} \leq \|T\|_{cb}$. The latter inequality is a consequence of the calculation
\[
        \|Tx + Y\| = \inf\{\|Tx + y\|: y \in Y\} \leq \inf\{\|T(x + y)\|: y \in Y\} 
		\leq \|T\|\inf\{\|x + y\|: y \in Y\} = \|T\|\|x + Y\|,
\]
as well as the fact that under the isometric isomorphism $M_n(X/Y) = M_n(X)/M_n(Y)$, $(T/Y)_n$ is identified
with $T_n/M_n(Y)$.

\begin{proposition}[\cite{BEZ}, Proposition 5.8] \label{IV.B.2}
Let $X$ be an operator space, $Y$ be a closed linear subspace of $X$, and $T:X \to X$ be a linear map such that
$T(Y) \subset Y$.
\begin{enumerate}
\item[(i)] If $T \in \Ml(X)$, with $\|T\|_{\Ml(X)} \leq 1$, then $T|_Y \in \Ml(Y)$ and $T/Y \in \Ml(X/Y)$,
with $\|T|_Y\|_{\Ml(Y)} \leq 1$ and $\|T/Y\|_{\Ml(X/Y)} \leq 1$.
\item[(ii)] If $T \in \Al(X)$ and $T^\star(Y) \subset Y$, then $T|_Y \in \Al(Y)$ and $T/Y \in \Al(X/Y)$. 
Furthermore, $(T|_Y)^\star = T^\star|_Y$ and $(T/Y)^\star = T^\star/Y$.
\end{enumerate}
\end{proposition}

We begin by examining subspaces. Our first result is a
hereditary property of right $M$-summands (resp. right $M$-ideals):

\begin{proposition} \label{IV.B.3}
Let $X$ be an operator space, $J$ be a closed linear subspace of $X$, and $Y$ be a closed linear subspace 
of $J$.
\begin{enumerate}
\item[(i)] If $Y$ is a right $M$-summand of $X$, then $Y$ is a right $M$-summand of $J$.
\item[(ii)] If $Y$ is a right $M$-ideal of $X$, then $Y$ is a right $M$-ideal of $J$.
\end{enumerate}
\end{proposition}

\begin{proof}
(i) Let $P \in \Al(X)_{sa}$ be the left $M$-projection onto $Y$. Then $P(J) = Y \subset J$. 
Therefore, by Proposition \ref{IV.B.2}, $P|_J \in \Al(J)_{sa}$. 
Clearly, $P|_J$ is an idempotent with 
range $Y$. Thus, $Y$ is a right $M$-summand of $J$. 

(ii) By (i), $Y^{\perp\perp}$ is a right $M$-summand 
of $J^{\perp\perp}$. It follows that $Y^{\perp \perp}$ is a right $M$-summand of $J^{**}$, which says that 
$Y$ is a right $M$-ideal of $J$.
\end{proof}

On the other hand, if $Y$ is a right $M$-ideal of $J$ and $J$ is a 
right $M$-ideal of $X$, 
then it need not be the case that $Y$ is a right $M$-ideal of $X$ (compare 
with \cite{HWW}, Proposition I.1.17).
The following  example of this phenomenon is due to K. 
Kirkpatrick (an R.E.U. student advised
by the first author):

\begin{example} \label{IV.B.4}
Let
\[
	X = \left\{\begin{bmatrix} a & 0\\ b & 0\\ 0 & \frac{a + b}{2}\\ 0 & c \end{bmatrix}:
		a, b, c \in \mathbb{C}\right\} \subset M_{4,2}.
\]
Define
\[
	J = \left\{\begin{bmatrix} a & 0\\ b & 0\\ 0 & \frac{a + b}{2}\\ 0 & 0 \end{bmatrix}:
		a, b \in \mathbb{C}\right\} \text{ and }
		Y = \left\{\begin{bmatrix} a & 0\\ 0 & 0\\ 0 & \frac{a}{2}\\ 0 & 0 
		\end{bmatrix}: a \in \mathbb{C}\right\}.
\]
It is easy to see that $J$ is a right $M$-summand of $X$, and that the `natural' projection $P$ of
$X$ onto $J$ is the corresponding left $M$-projection. Likewise, $Y$ is a right $M$-summand of $J$, and
the natural projection $Q$ of $J$ onto $Y$ is the corresponding left $M$-projection. This
follows from the fact that the map
\[
	C_2 \to J:\begin{bmatrix} a\\ b \end{bmatrix}
		\mapsto \begin{bmatrix} a & 0\\ b & 0\\ 0 & \frac{a + b}{2}\\ 0 & 0
		\end{bmatrix}
\]
is a completely isometric isomorphism. For future reference, let $\tilde{J}$ 
be the complementary
right $M$-summand of $J$ in $X$ (namely, $\tilde{J}$ 
is the copy of $\mathbb{C}$ in the $4-2$ corner of $X$), and let
$\tilde{Y}$ be the complementary right $M$-summand of $Y$ in $J$.
Now suppose by way of contradiction that $Y$ is 
a right $M$-summand of $X$. Since $E = Q \circ P:X \to X$ is a contractive projection
with range $Y$, it must be the left $M$-projection of $X$ onto $Y$ (cf.\ the fourth
`bullet' 
in Section \ref{I.C}). But 
$\nu_E^c:X \to C_2(X)$ is not isometric, as can be seen by applying it to
\[
	x_\lambda = \begin{bmatrix} 1 & 0\\ 1 & 0\\ 0 & 1\\ 0 & \lambda \end{bmatrix}
\]
for $\lambda$ sufficiently large. Indeed, 
$\|\nu_E^c(x_\lambda)\| = \max\left\{\sqrt{2},\sqrt{\frac{1}{2} + \lambda^2}\right\}$ whereas 
$\|x_\lambda\| = \max\{\sqrt{2},\sqrt{1 + \lambda^2}\}$. This is a contradiction.
\end{example}

We remark that by symmetry, $\tilde{Y}$ is also a right $M$-summand of $J$ which fails to be a right
$M$-summand of $X$. 

Before moving on to quotients, we identify the `culprit' behind Example \ref{IV.B.4}. Namely, we prove in
Proposition \ref{herind} below that if $J$ is a right $M$-summand of a dual operator space $X$, then all
right $M$-summands of $J$ are right $M$-summands of $X$ if and only if $\Al(J) \cong P\Al(X)P$, where
$P$ is the left $M$-projection of $X$ onto $J$. In the case of Example \ref{IV.B.4}, $\Al(J) = M_2$ whereas
$P\Al(X)P = \ell_2^\infty$. We begin with an important lemma.

\begin{lemma} \label{IV.B.4b}   Let $X$ be an operator space,  
let $P$ be a left $M$-projection on $X$,
and let $\rho$ be the map $T \mapsto T|_{P(X)}$ from 
$P \Al(X) P$ into $\Al(P(X))$.
\begin{itemize} \item [(1)]  $\rho$
identifies the $C^*$-algebra $P \Al(X) P$ $*$-isomorphically with 
a (possibly proper) $C^*$-subalgebra of $\Al(P(X))$.    
\item [(2)]  If $X$ is a dual operator space, then $\rho$ is a
weak-* homeomorphic $*$-isomorphism from $P \Al(X) P$ onto a 
$W^*$-subalgebra of $\Al(P(X))$.
\end{itemize}
 \end{lemma} 

\begin{proof}   
Let $\rho$ be the above map, which is clearly a one-to-one homomorphism
into $B(P(X))$.  By Proposition \ref{IV.B.2} we see that
$\rho$  does indeed map into $\Al(P(X))$, 
and is a $*$-homomorphism.  Hence $\rho$  is
also isometric and has closed range.   The remark before the Lemma
shows that this range may be a proper $C^*$-subalgebra.   This
gives (1).  

If $X$ is a dual operator space then $P(X)$ is
weak-* closed, by a remark after Theorem \ref{left_M}. 
Then (2) follows from (1), the Krein-Smulian theorem
\ref{KS}, and Theorem \ref{wst}.   \end{proof}

\begin{proposition} \label{herind}
Let  $X$ be a dual operator space with a right $M$-summand
$J$ and corresponding left $M$-projection
$P$.  Then every right $M$-summand of $J$ is
a right $M$-summand of $X$
if and only if
$\Al(J) \cong P \Al(X) P$.
\end{proposition}

\begin{proof}
$(\Leftarrow)$  This implication does not need
$X$ to be a dual space. Suppose that $Q$ is a left $M$-projection 
on $J$. By assumption, $Q = \rho(T)$ for some $T \in P \Al(X) P \subset \Al(X)$.
Since $\rho$ is an injective *-homomorphism,
$T$ is a left $M$-projection on $X$. We have
$T(X) = T P(X) = Q P(X) = \ran(Q)$. Hence
$\ran(Q)$ is a right $M$-summand of $X$.

$(\Rightarrow)$
Because the orthogonal
projections generate any von Neumann algebra,
we need only prove that every left $M$-projection
$Q$ on $J$ equals $P R P|_J$ for a left
$M$-projection $R$ on $X$.  By hypothesis,
there exists a left $M$-projection $R$ on $X$
such that $\ran(R) = \ran(Q)$. Evidently,
$P R P|_J$ is a left $M$-projection on $J$ with range
$\ran(Q)$. It follows from the uniqueness of
left $M$-projections with a given range that $Q = P R P|_J$.
\end{proof}                            

Perhaps the most frustrating  problem in the theory of one-sided
$M$-structure is that the equivalent conditions in the last Proposition 
may fail (see Example \ref{IV.B.4}). 
In any case, this `problem' suggests a definition:

\begin{definition} \label{IV.B.4c}  A 
right $M$-summand $J$ of
an operator space $X$ is called {\em hereditary} if
every right $M$-summand of $J$  is also a
right $M$-summand of $X$.   
\end{definition}

One might then ask for a
reasonable sufficient condition for a right $M$-summand to be
hereditary,
and a first guess at might be that all $M$-projections involved commute.
However a simple
modification of Example \ref{IV.B.4}, with the `$a$-$b$
column' replaced
by $diag \{a,b\}$, shows that such `commutativity' does not help.
(We remark that for this modified example it is quickest to
compute the multiplier algebras explicitly from their
formulation in terms of the `triple envelope' (see \cite{Hamana,Shilov});
in this case the triple envelope is easily seen to be
$\ell^\infty_2 \oplus C_2$.)   

The following is a rather strong sufficient condition.
We say that 
 a closed subspace
$J$ of an operator space $X$ is a {\em Shilov subspace} if
there is a triple envelope $Z$ of $X$ 
such that the triple subsystem of $Z$ generated by $J$
(that is, the smallest closed subspace of $Z$ containing
$J$ and closed under the `triple product')  is a triple 
envelope of $J$.  By the universal property of the triple envelope 
(see e.g.\ \cite{Hamana}), one can show that
this notion does not depend on the
particular triple envelope $Z$ above.

\begin{proposition} \label{shsub}  If a right $M$-summand $J$ of
an operator space $X$ is a Shilov subspace of $X$, then 
$J$ is hereditary.
\end{proposition} 

\begin{proof}   Let $P$ be the left $M$-projection of $X$ onto $J$.
If $J$ is a Shilov subspace, then the
restriction of the `Shilov inner product' on $X$ (see Section
\ref{II.BBB}), to $J$, is a `Shilov inner product' for
$J$.   If $Q$ is a left $M$-projection on $J$ then 
$Q$ is a projection in $\Al(J)$, and so by (\ref{adeq}) we have
$$\langle Qx , y \rangle = \langle x , Q y \rangle , \qquad
x, y \in J .$$
Since $P Q = Q$, we have for $x, y \in X$ that  
 $$\langle Q P x , y \rangle = \langle Q P x , P y
\rangle = \langle P x , Q P y  \rangle = \langle x , Q P y\rangle.$$
Thus $Q P$ is a left adjointable contractive projection on $X$
with range $J$, which establishes the asserted statement.
\end{proof}

Turning to quotients, we have the following positive result, 
which tells us that right $M$-summands and right
$M$-ideals are stable under the taking of quotients.

\begin{proposition} \label{IV.B.5}
Let $X$ be an operator space, $J$ be a closed linear subspace of $X$, and $Y$ be a closed linear subspace of $J$.
\begin{enumerate}
\item[(i)] If $J$ is a right $M$-summand of $X$, then $J/Y$ is a right $M$-summand of $X/Y$.
\item[(ii)] If $J$ is a right $M$-ideal of $X$, then $J/Y$ is a right $M$-ideal of $X/Y$.
\end{enumerate}
\end{proposition}
 
\begin{proof}
(i) Let $P \in \Al(X)_{sa}$ be the left $M$-projection onto $J$. Then $P(Y) = Y$. Therefore, by
Proposition \ref{IV.B.5}, $P/Y \in \Al(X/Y)_{sa}$. Clearly, $P/Y$ is an idempotent with range $J/Y$. Thus, 
$J/Y$ is a right $M$-summand of $X/Y$. 

(ii) By (i), $J^{\perp\perp}/Y^{\perp\perp}$ is a right $M$-summand of 
$X^{**}/Y^{\perp\perp}$. It follows that $(J/Y)^{\perp\perp}$ is a right $M$-summand of $(X/Y)^{**}$, which says 
that $J/Y$ is a right $M$-ideal of $X/Y$.
\end{proof}

On the other hand, not every right $M$-ideal of a quotient of $X$ 
arises as the 
image of a right $M$-ideal of $X$ under the quotient map 
(compare again with \cite{HWW}, Proposition I.1.17). 
An instance of this phenomenon is furnished by the 
following example:

\begin{example} \label{IV.B.6}   
Let $X$, $J$, $\tilde{J}$, $Y$, and $\tilde{Y}$ have the same meanings as in Example \ref{IV.B.4}. 
We will show that $(Y + \tilde{J})/\tilde{J}$ is a right $M$-summand of $X/\tilde{J}$, but that there does not exist
a right $M$-summand $K$ of $X$
such that $(Y + \tilde{J})/\tilde{J} = K/\tilde{J}$.
To prove the first assertion, note that under the canonical completely isometric isomorphism $X/\tilde{J} = J$,
$(Y + \tilde{J})/\tilde{J}$ is identified with $Y$. Since $Y$ is a right $M$-summand of $J$, 
$(Y + \tilde{J})/\tilde{J}$ is a right $M$-summand of $X/\tilde{J}$. To prove the second assertion, we
proceed by contradiction. So let $K$ be a right $M$-summand of $X$ 
such that
$K/\tilde{J} = (Y + \tilde{J})/\tilde{J}$.   Replacing
$K$ by $K + \tilde{J}$, which is still a right $M$-summand of $X$
by Proposition \ref{II.A.3} for example, we may suppose that
$K$ contains $\tilde{J}$.
It is easy to see that $K$ must equal $Y + \tilde{J}$. Let $F$ be
the left $M$-projection 
of $X$ onto $Y + \tilde{J}$. Since $\tilde{J} \subset Y + \tilde{J}$, we 
have
 $P^\perp F = FP^\perp = P^\perp$. Therefore, if
\[
	F\left(\begin{bmatrix} 0 & 0\\ 1 & 0\\ 0 & \frac{1}{2}\\ 0 & 0 \end{bmatrix}\right) =
		\begin{bmatrix} a & 0\\ 0 & 0\\ 0 & \frac{a}{2}\\ 0 & c \end{bmatrix},
\]
then we may conclude that $c = 0$. On the other hand, the fact that $\nu_F^c:X \to C_2(X)$ is (completely) 
isometric allows us to conclude that $a = 0$. That is, $F$ is the `natural' projection of $X$ onto $Y + \tilde{J}$. 
But then $F^\perp$ is the natural projection of $X$ onto $\tilde{Y}$, and so $\tilde{Y}$ is a right $M$-summand of 
$X$. This contradicts the remark following Example \ref{IV.B.4}.
\end{example}

The behavior of one-sided $L$-structure with respect to subspaces and quotients is similar. This is a consequence
of the following lemma.

\begin{lemma} \label{IV.B.7}
Let $X$ be an operator space, $Y$ be a closed linear subspace of $X$, and $T:X \to X$ be a linear map such that
$T(Y) \subset Y$. Then
\begin{enumerate}
\item[(i)] Under the completely isometric isomorphism $Y^* = X^*/Y^\perp$, $(T|_Y)^*:Y^* \to Y^*$ is identified
with $T^*/Y^\perp:X^*/Y^\perp \to X^*/Y^\perp$.
\item[(ii)] Under the completely isometric isomorphism $(X/Y)^* = Y^\perp$, $(T/Y)^*:(X/Y)^* \to (X/Y)^*$ is
identified with $T^*|_{Y^\perp}:Y^\perp \to Y^\perp$.
\end{enumerate}
\end{lemma}

The proof of this lemma is an easy exercise which we omit.

\begin{corollary} \label{IV.B.8}
Let $X$ be an operator space, $Y$ be a closed linear subspace, and $T \in \Cr(X)$. If $T(Y) \subset Y$ and
$T^\star(Y) \subset Y$, then $T|_Y \in \Cr(Y)$ and $T/Y \in \Cr(X/Y)$. Furthermore, 
$(T|_Y)^\star = T^\star|_Y$ and $(T/Y)^\star = T^\star/Y$.
\end{corollary}

\begin{proof}
Recall from Section \ref{I.D} that for any operator space $V$, the linear isometry $CB(V) \to CB(V^*):T \mapsto T^*$ 
restricts to a *-anti-isomorphism $\Cr(V) \to \Al(V^*)$ which is also a weak-* homeomorphism. Therefore, to
show that $T|_Y \in \Cr(Y)$, it suffices to show that $(T|_Y)^* \in \Al(Y^*)$. Likewise, to show that
$T/Y \in \Cr(X/Y)$, it suffices to show that $(T/Y)^* \in \Al((X/Y)^*)$. By Lemma \ref{IV.B.7}, it is equivalent
to show that $T^*/Y^\perp \in \Al(X^*/Y^\perp)$ and $T^*|_{Y^\perp} \in \Al(Y^\perp)$, respectively. But this 
follows from Proposition \ref{IV.B.2}, as do the formulas appearing in the statement of the corollary.
\end{proof}

\begin{corollary} \label{IV.B.9}
Let $X$ be an operator space, $J$ be a closed linear subspace of $X$, and $Y$ be a closed linear subspace of
$J$.
\begin{enumerate}
\item[(i)] If $Y$ is a left $L$-summand of $X$, then $Y$ is a left $L$-summand of $J$.
\item[(ii)] If $J$ is a left $L$-summand of $X$, then $J/Y$ is a left $L$-summand of $X/Y$.
\end{enumerate}
\end{corollary}

To prove this corollary one mimics the proofs of Propositions \ref{IV.B.3} and \ref{IV.B.5}, with the role
of Proposition \ref{IV.B.2} being played by Corollary \ref{IV.B.8}. 
We leave the details to the interested
reader.

\subsection{Dual and Bidual} \label{IV.C}

As we saw in Section \ref{I.D}, the linear isometry $CB(X) \to CB(X^*):T \mapsto T^*$ restricts to
a *-anti-isomorphism $\Cr(X) \to \Al(X^*)$ which is also a weak-* homeomorphism. We claim that it
also restricts to a *-anti-homomorphism $\Al(X) \to \Cr(X^*)$. Our argument will be indirect: we will show that
the linear isometry $\Phi:CB(X) \to CB(X^{**}):T \to T^{**}$ restricts to a *-homomorphism $\Al(X) \to \Al(X^{**})$,
which then can be composed with the inverse of the *-anti-isomorphism $\Cr(X^*) \to \Al(X^{**})$. It is interesting
that a direct argument does not present itself. Indeed, since the generic element of $\Al(X)$ is not a norm
limit of linear combinations of left $M$-projections on $X$, it is not clear a priori that the Banach space 
adjoint of an element of $\Al(X)$ is an element of $\Cr(X^*)$. Nonetheless, it is true. Obviously, by composing
the *-anti-isomorphism $\Cr(X) \to \Al(X^*)$ with the *-anti-homomorphism $\Al(X^*) \to \Cr(X^{**})$, we
obtain a *-homomorphism $\Cr(X) \to \Cr(X^{**}):T \mapsto T^{**}$. There is no reason to expect that this
map is weak-* continuous. We turn now to the details.\\

\begin{proposition} \label{IV.C.0}   Let $T:X \to X$ be a linear map
on an operator space. 
Then $T$ is an element of the closed unit ball of $\Ml(X)$ if and only if 
$T^{**}$ is an element of the closed unit ball of $\Ml(X^{**})$. 
\end{proposition}

\begin{proof}  By two applications of Theorem \ref{BEZ}
we have that
\[
	\|T\|_{\Ml(X)} \leq 1 \Leftrightarrow \|\tau_T^c\|_{cb} \leq 1 \Leftrightarrow 
		\|(\tau_T^c)^{**}\|_{cb} \leq 1 \Leftrightarrow \|\tau_{T^{**}}^c\|_{cb} \leq 1 
		\Leftrightarrow \|T^{**}\|_{\Ml(X^{**})} \leq 1.
\]
(A direct proof without using Theorem \ref{BEZ} may also be given).
\end{proof}

Thus the map
$\Phi$ above
restricts to a unital isometric homomorphism $\Ml(X) \to \Ml(X^{**})$. 
It is easy to check that a unital contraction between unital
Banach algebras takes Hermitian elements  to Hermitian elements.
For a unital subalgebra ${\mathcal A}$ of a 
unital $C^*$-algebra ${\mathcal B}$ it is clear that the 
Hermitian elements of ${\mathcal A}$ are the 
selfadjoint elements of ${\mathcal B}$ which are in ${\mathcal A}$,
or equivalently the selfadjoint elements of the diagonal 
of ${\mathcal A}$.  Since $\Al(X)$ is the
diagonal of $\Ml(X)$ (cf. Section \ref{I.C}),  
it follows from the last few facts that $\Phi$ 
takes selfadjoint elements 
in  $\Al(X)$ to selfadjoint elements in $\Al(X^{**})$.
Hence $\Phi$ restricts to an isometric 
 *-homomorphism from $\Al(X)$  into $\Al(X^{**})$.
Consequently,
$\Al(X)$ may be regarded as a 
unital $C^*$-subalgebra of $\Al(X^{**})$. Typically, the inclusion
$\Al(X) \subset \Al(X^{**})$ is proper.

A simple consequence of the preceding discussion is the following result, which we will
make use of shortly (see the remarks at the end of the section).

\begin{lemma} \label{ab2s}
Let $X$ be an operator space, and let
$T \in \Al(X)$.  Then $|T|^{**} = |T^{**}|$.
\end{lemma}

\begin{proof} We have observed that 
the map $\Al(X) \to \Al(X^{**}): T \mapsto T^{**}$ is a $*$-homomorphism.
Thus $|T|^{**}$ is positive in the $C^*$-algebra
$\Al(X^{**})$.  We also deduce that $$(|T|^{**})^2 =
(T^\star T)^{**} = (T^{**})^\star T^{**}.$$
Hence $|T|^{**} = |T^{**}|$.  \end{proof}

We have seen that that $T \in \Ml(X)$ if and only if 
$T^{**} \in \Ml(X^{**})$. Likewise, we saw earlier that 
if $T \in \Al(X)$, then $T^{**} \in \Al(X^{**})$. 
It is natural to expect therefore that if $T^{**} \in \Al(X^{**})$,
then  $T \in \Al(X)$.  In fact this is not true in general.
Indeed, if $\mathcal{A}$ is a $C^*$-algebra and 
$A \in \mathcal{LM(A)} \backslash \mathcal{M(A)}$, then the map $L_A:\mathcal{A} \to \mathcal{A}$
of left multiplication by $A$ is a left multiplier of $\mathcal{A}$ which is not left adjointable. On the
other hand, $L_{A}^{**}:\mathcal{A}^{**} \to \mathcal{A}^{**}$ is again the map of left multiplication by $A$,
which is a left adjointable multiplier of $\mathcal{A}^{**}$ (since $A \in \mathcal{A}^{**}$). The following
proposition clarifies the situation.

\begin{proposition} \label{IV.C.1}
Let $X$ be an operator space and $T:X \to X$ be linear. 
\begin{itemize} \item [(i)]  $T \in \Al(X)_{sa}$ if and only if
$T^{**} \in \Al(X^{**})_{sa}$. 
 \item [(ii)]  $T \in \Al(X)$ if and only if $T^{**} \in \Al(X^{**})$
and $T = T_1 + iT_2$, where $T_1$ and $T_2$ are
Hermitian elements of $B(X)$.
\end{itemize} 
\end{proposition}

\begin{proof}
(i)  We have proven above that if $T \in \Al(X)_{sa}$, then 
$T^{**} \in \Al(X^{**})_{sa}$. Conversely, suppose that
$T^{**} \in \Al(X^{**})_{sa}$. Since $T^{**}(X) = T(X) \subset X$, 
$T = T^{**}|_X \in \Al(X)_{sa}$ by Proposition \ref{IV.B.2}. 

(ii) ($\Rightarrow$) We proved above that if 
$T \in \Al(X)$ then $T^{**} \in \Al(X^{**})$.  The other assertion 
follows since $\Al(X)$ is a $C^*$-algebra, and is a unital 
subalgebra of $B(X)$.   

($\Leftarrow$)   Suppose that $T = T_1 + iT_2$, where $T_1$ and 
$T_2$ are Hermitian elements of $B(X)$. Then $T^{**} = T_1^{**} + iT_2^{**}$, where $T_1^{**}$ and $T_2^{**}$ are
Hermitian elements of $B(X^{**})$. On the other hand, 
assuming that $T^{**} \in \Al(X^{**})$, 
one has that $T^{**} = R + iS$, where $R, S \in \Al(X^{**})_{sa}$. In 
particular, $R$ and $S$ are Hermitian elements of $B(X^{**})$. Thus, by Lemma \ref{App.A.8}, 
$T_1^{**}, T_2^{**} \in \Al(X^{**})_{sa}$. By what has already been shown, $T_1, T_2 \in \Al(X)_{sa}$ and so
$T \in \Al(X)$. 
\end{proof}
 
{\bf Remark.}  From this last result 
one may give a solution to a problem raised in
\cite{Behrends84}.
To describe the problem we first
list some notation.  Behrends defines a map $R$ on 
a Banach space $X$ to
`possess an adjoint' if there exists
an operator $S$ on $X$ (necessarily unique) such that  both $\frac{1}{2}(R+S)$ and $
\frac{1}{2i}(R-S)$ are Hermitian in $B(X)$.  We shall 
call such an $S$ a  `Behrends adjoint of
$R$', and write $S = R^{\flat}$.
Behrends asks the following question: if $T$ is a map on a Banach
space, such that $T^* \in B(X^*)$ possesses an adjoint
in  his sense, then is this adjoint necessarily
weak-* continuous?  This is equivalent to asking whether
$T$ `possesses a Behrends  adjoint' if and only if $T^*$ does.
In fact, any linear map $T$ on an operator 
space $X$ with $T \notin \Al(X)$ but $T^{**} \in  \Al(X^{**})$
(for example a left multiplier
$T$ on a $C^*$-algebra which is not a two-sided multiplier)
gives rise to a counterexample to the question.
Indeed, given such a $T$, let $V$ be the involution in $\Al(X^{**})$ of $T^{**}$.
We have that 
$\frac{1}{2}(T^{**} + V)$ and $\frac{1}{2i}(T^{**} - V)$ are selfadjoint in
$\Al(X^{**})$, 
and consequently Hermitian also in $B(X^{**})$.
We know that $V$ is weak-* continuous by a fact listed 
after Theorem  
\ref{wst}.  Hence $V = S^*$ for some $S \in B(X^*)$.
Now $\frac{1}{2}(T^{*} + S)$ and  $\frac{1}{2i}(T^{*} - S)$ are 
both Hermitian in $B(X^{*})$,
since their Banach space adjoints are Hermitian.   
So $T^*$ `possesses a Behrends  adjoint':
$(T^*)^\flat = S$. However  $S$ is
not weak*-continuous.  For if it was, and $S = R^*$ say, then 
$\frac{1}{2}(T + R)$ and $\frac{1}{2i}(T - R)$
would be Hermitian. Thus by the last Proposition,
$T$ would be left adjointable, which is a contradiction. \\

In general, there is no way to produce from an element of $\Al(X^{**})$ a corresponding element of $\Al(X)$.
If, however, $X$ is a dual space, then there is. More precisely, we have the following result:

\begin{theorem} \label{IV.C.2}
Let $X$ be an operator space. Then there 
exists a normal (i.e. weak-* continuous) conditional expectation $E$ of 
$\Al(X^{***})$ onto $\Al(X^*)$. 
 The weak-* topologies on these
spaces are the canonical ones with respect to which they 
are $W^*$-algebras.
\end{theorem}

\begin{proof}
For any operator space $Y$, denote by $\iota_Y:Y \to Y^{**}$ the canonical completely isometric inclusion. We 
define $E:CB(X^{***}) \to CB(X^*)$ by the formula
\[
	E(T) = \iota_X^* \circ T \circ \iota_{X^*}.
\]
Obviously, $E$ is unital and contractive. Now suppose that
$T \in \Ml(X^{***})$, with $\|T\|_{\Ml(X^{***})} \leq 1$.
Then (claim) $E(T) \in \Ml(X^*)$, with $\|E(T)\|_{\Ml(X^*)} \leq 1$. Indeed, since $\|\tau_T^c\|_{cb} \leq 1$, one 
has that
\[
	\left\|\begin{bmatrix} E(T)(f)\\ g \end{bmatrix}\right\| = \left\|
		\begin{bmatrix} \iota_X^*(T(\iota_{X^*}(f)))\\ \iota_X^*(\iota_{X^*}(g))
		\end{bmatrix}\right\| \leq \left\|\begin{bmatrix} T(\iota_{X^*}(f))\\
		\iota_{X^*}(g) \end{bmatrix}\right\| \leq \left\|\begin{bmatrix}
		\iota_{X^*}(f)\\ \iota_{X^*}(g) \end{bmatrix}\right\| =
		\left\|\begin{bmatrix} f\\ g \end{bmatrix}\right\|
\]
for all $f, g \in X^*$, which says that $\|\tau_{E(T)}^c\| \leq 1$. Repeating this calculation for matrices over 
$X^*$ shows that in fact $\|\tau_{E(T)}^c\|_{cb} \leq 1$, which in turn implies that $E(T) \in \Ml(X^*)$, with 
$\|E(T)\|_{\Ml(X^*)} \leq 1$. Since $E:\Ml(X^{***}) \to \Ml(X^*)$ is unital and contractive, it restricts to 
a positive map from $\Al(X^{***})$ to $\Al(X^*)$. The map $E:\Al(X^{***}) \to \Al(X^*)$ is surjective, since
\[
	E(T^{**}) = \iota_X^* \circ T^{**} \circ \iota_{X^*} = \iota_X^* \circ \iota_{X^*} \circ T = T
\]
for all $T \in \Al(X^*)$. Similarly, we have 
\[
	E(S \circ T^{**}) = \iota_X^* \circ S \circ T^{**} \circ \iota_{X^*} =
		\iota_X^* \circ S \circ \iota_{X^*} \circ T = E(S) \circ T
\]
for all $S \in \Al(X^{***})$ and all $T \in \Al(X^*)$. Taking involutions we also have 
$E(T^{**} \circ S) = T \circ E(S)$, so that $E$ is a conditional expectation. Finally, we claim that it is normal. 
Indeed, if $(T_i)$ is a bounded net in $\Al(X^{***})$, $T \in \Al(X^{***})$, and $T_i \to T$ weak-*, then for any 
$f \in X^*$, we have that
\[
	T_i(\iota_{X^*}(f)) \to T(\iota_{X^*}(f))
\]
with respect to $\sigma(X^{***},X^{**})$.  This implies that
\[
	E(T_i)(f) = \iota_X^*(T_i(\iota_{X^*}(f))) \to \iota_X^*(T(\iota_{X^*}(f)))
		= E(T)(f)
\]
with respect to $\sigma(X^*,X)$. Since the choice of $f$ was arbitrary, $E(T_i) \to E(T)$ 
in the weak-* topology of $\Al(X^*)$.
\end{proof}

\begin{corollary} \label{IV.C.3}
Let $X$ be an operator space and $J$ be a right $M$-ideal of $X^*$. Then 
$\overline{J}^{\wks}$, the closure of $J$ in the 
weak-* topology of $X^*$, is a right
$M$-summand of $X^*$. 
\end{corollary}

\begin{proof}
Let $P \in \Al(X^{***})$ be the left $M$-projection onto $J^{\perp\perp}$. Set $Q = E(P) \in \Al(X^*)$. We claim 
that $Q$ is a left $M$-projection with range $\overline{J}^{\wks}$. 
Clearly, $Q$ is a contraction in 
$\Al(X^*)$, and 
it is positive as an element of that 
$C^*$-algebra since $E$ is a positive map.    
Now let $f \in J$. Then $\iota_{X^*}(f) \in J^{\perp\perp}$, which implies that
\[
	Q(f) = \iota_X^*(P(\iota_{X^*}(f))) = \iota_X^*(\iota_{X^*}(f)) = f.
\]
Since $Q$ is weak-* continuous, as noted 
after Theorem \ref{wst},
we have  $Q(f) = f$ for all $f \in \overline{J}^{\wks}$. Since
\[
	Q(X^*) = \iota_X^*(P(\iota_{X^*}(X^*))) \subset \iota_X^*(J^{\perp\perp}) = 
		\iota_X^*(\overline{\iota_{X^*}(J)}^{\wks}) \subset 
		\overline{\iota_X^*(\iota_{X^*}(J))}^{\wks} = \overline{J}^{\wks},
\]
we deduce that  $\ran(Q) = \overline{J}^{\wks}$ and 
$Q^2 = Q$, which proves the claim.
\end{proof}

\begin{corollary} \label{IV.C.4}
Let $X$ be an operator space and $J$ be a right $M$-ideal of $X^*$. Then $\phantom{}^\perp J$ is a left $L$-summand
of $X$.
\end{corollary}

\begin{proof}
$\phantom{}^\perp J$ is a closed linear subspace of $X$ whose annihilator is a right $M$-summand of $X^*$:
\[
	(\phantom{}^\perp J)^\perp = \overline{J}^{\wks}.
\]
\end{proof}

\begin{corollary} \label{IV.C.4b}  A dual operator space with no nontrivial
right $M$-summands has only weak-* dense, or trivial, right $M$-ideals.
\end{corollary}

\begin{corollary} \label{IV.C.5}
Let $X$ be an operator space and $J$ be a right $M$-ideal of $X^*$. If $x \in \overline{J}^{\wks}$, then there exists
a net $(x_i)$ in $J$, with $\|x_i\| \leq \|x\|$ for all $i$, such that $x_i \to x$ weak-*.
\end{corollary}

\begin{proof}
Let $P \in \Al(X^{***})$ be the left $M$-projection onto $J^{\perp\perp}$, and $Q = E(P) \in \Al(X^*)$ be
the left $M$-projection onto $\overline{J}^{\wks}$. Set $\hat{x} = \iota_{X^*}(x) \in X^{***}$. Then 
$P\hat{x} \in J^{\perp\perp} \cong J^{**}$. Thus, there exists a net $(x_i) \in J$, with 
$\|x_i\| \leq \|P\hat{x}\| \leq \|\hat{x}\| = \|x\|$ for all $i$, such that
$\iota_{X^*}(x_i) \to P\hat{x}$ weak-*. But then $x_i = Qx_i = \iota_X^*(P(\iota_{X^*}(x_i))) \to \iota_X^*(P(\hat{x})) =
Qx = x$ weak-*.
\end{proof}
                 
{\bf Remarks.}  1)  The classical analogs of the last four
 results are true too, and follow as a special case of the above.
We have not seen these results in the literature, though.

2) A good illustration of Corollary \ref{IV.C.4b}
is $H^\infty(D)$, the classical 
$H^\infty$ space of the unit disk.  The $M$-projections on $H^\infty(D)$ 
(which are the same as the one-sided $M$-projections by
facts in Section \ref{III.B}) are multiplications
by characteristic functions in $H^\infty(D)$, and hence are trivial.
Thus by  Corollary \ref{IV.C.4b}    
the $M$-ideals (which coincide with the one-sided $M$-ideals
by \ref{III.B}), and which are known to be the closed ideals 
of $H^\infty(D)$ possessing a bounded approximate identity
(see e.g.\ \cite{HWW}), are either trivial or 
weak-* dense, by the above. In this
 spirit, note that the disk algebra $A(D)$ also has no nontrivial
$M$-projections, but has plenty of $M$-ideals.  

Another interesting example comes from work of Powers \cite{Powers}.
Namely, there exists a $C^*$-algebra with no nontrivial $M$-ideals,
but with uncountably many pairwise non-isometric $M$-summands in the
second dual \cite[Proposition V.4.6]{HWW}. It is easy to 
show using basic spectral theory, however, that every nontrivial
$C^*$-algebra has nontrivial closed left ideals. On the
 other hand,
there are 
plenty of nonselfadjoint unital operator algebras with no nontrivial
left $M$-ideals.  In fact it is easy to find two or three dimensional 
`upper triangular' matrix algebras with this property.   

3) Obviously, there is a normal conditional expectation $\tilde{E}:\Cr(X^{**}) \to \Cr(X)$ which makes the following
diagram commute:
\[
\xymatrix{
	\Al(X^{***}) \ar[r]^E & \Al(X^*)\\
	\Cr(X^{**}) \ar[u]^{T \mapsto T^*}\ar[r]_{\tilde{E}} & \Cr(X) \ar[u]_{T \mapsto T^*}
}
\]

Despite this, and in contrast to
Corollary \ref{IV.C.3}, it is not true that the weak-* closure of
a left $L$-summand in a dual space is also a 
left $L$-summand.    To see this we will use the 
fact that the intersection $J \cap K$ of two 
right $M$-summands $J$ and $K$ need not be a 
right $M$-ideal (see Example \ref{IV.D.6}).  Hence
$(J \cap K)^\perp$, which by Appendix \ref{App.A.2} (ii)
equals $\overline{J^\perp + K^\perp}^{\wks}$, 
is not a left $L$-summand.  But the last set is the 
weak-* closure of $\overline{J^\perp + K^\perp}$,
which is a left $L$-summand by Corollary \ref{II.A.4}.    

4)  One may use some of the techniques above
to construct a kind of polar decomposition for 
a left adjointable multiplier $T$ in a general operator space $X$. 
The idea is to apply
Theorem \ref{II.B.2} together with Lemma \ref{ab2s}, to the 
map $T^{**} \in  \Al(X^{**})$.   We obtain a partial 
isometry $W \in \Al(X^{**})$ such that 
$T^{**} = W |T|^{**}$ and $|T|^{**} = W^\star T^{**}$.
Restricting to $X$ we obtain $T = V |T|$ and $|T| = U T$,
where $V = W|_X$ and $U = W^\star|_X$.   Unfortunately 
$V$ (resp. $U$) maps into $X^{**}$, although on 
$\ran(|T|)$  (resp. $\ran(T)$) it maps 
completely isometrically into $X$.

5) Corollary \ref{IV.C.5} is a `Kaplansky density theorem' for right $M$-ideals in a dual operator space.

\subsection{Sum and Intersection} \label{IV.D}

In this section we describe how one-sided $M$-summands, $M$-ideals, and $L$-summands behave with respect to
sums and intersections. Generally speaking, the behavior with respect to sums is better than the
behavior with respect to intersections.

Our first result follows immediately from Proposition \ref{II.A.5}.

\begin{proposition} \label{IV.D.1}
\begin{enumerate}
\item[(i)] If $X$ is a dual operator space and $\{J_a: a \in A\}$ is a family of right $M$-summands
of $X$, then $\cap_{a \in A} J_a$ and $\overline{\spn}^{\wks}\left\{\cup_{a \in A} J_a\right\}$ are
right $M$-summands of $X$.
\item[(ii)] If $X$ is an operator space and $\{J_a: a \in A\}$ is a family of left $L$-summands of $X$,
then $\cap_{a \in A} J_a$ and $\overline{\spn}\left\{\cup_{a \in A} J_a\right\}$ are left $L$-summands
of $X$.
\end{enumerate}
\end{proposition}

\begin{corollary} \label{IV.D.2}
Let $X$ be an operator space and $\{J_a: a \in A\}$ be a family of right $M$-ideals of $X$. Then
$\overline{\spn}\left\{\cup_{a \in A} J_a\right\}$ is a right $M$-ideal of $X$.
\end{corollary}

\begin{proof}
By Lemma \ref{App.A.2} and Proposition \ref{IV.D.1}, $\overline{\spn}\left\{\cup_{a \in A} J_a\right\}$ 
is a closed linear subspace of $X$ whose annihilator, $\cap_{a \in A} J_a^\perp$, is a left $L$-summand 
of $X^*$.
\end{proof}

As the following two examples show, the `closure' is necessary in both Proposition \ref{IV.D.1} and
Corollary \ref{IV.D.2}.

\begin{example} \label{IV.D.3}
Let $H$ be a Hilbert space, and $J_1$ and $J_2$ be closed linear subspaces of $H$ such that $J_1 + J_2$
is not closed (cf. \cite{Halmos}, \S 15). Let $X = H_c$, the column Hilbert space corresponding to
$H$. Then $J_1$ and $J_2$ are right $M$-summands of $X$ (cf. \cite{BEZ}, Proposition 6.10), but
$J_1 + J_2$ is not a right $M$-ideal of $X$. Consequently, the norm closure is necessary in 
Corollary \ref{IV.D.2}.
\end{example}

\begin{example} \label{IV.D.4}
There exists a von Neumann algebra $\mathcal{M}$ and projections $P, Q \in \mathcal{M}$ such that
$\overline{P\mathcal{M} + Q\mathcal{M}}$ (norm closure) is not equal to $R\mathcal{M}$ for any projection
$R \in \mathcal{M}$. Thus, there exist right $M$-summands $J_1$ and $J_2$ of a dual operator space $X$
such that $\overline{J_1 + J_2}$ is not a right $M$-summand of $X$. Consequently, the weak-* closure is
necessary in Proposition \ref{IV.D.1}. 
\end{example}

\begin{proof}
Let $H$ be a Hilbert space, and $H_1$ and $H_2$ be closed linear subspaces of $H$ such that
$H_1 + H_2$ is not closed. We may assume that $H = \overline{H_1 + H_2}$. Let $P \in B(H)$ be the projection
onto $H_1$ and let 
$Q \in B(H)$ be the projection onto $H_2$. Let $\mathcal{M} \subset B(H)$ be any
von Neumann algebra containing $P$ and $Q$. Suppose that 
$\overline{P\mathcal{M} + Q\mathcal{M}} = R\mathcal{M}$ for some projection $R \in \mathcal{M}$. Then
$RP = P$ and $RQ = Q$, so that $R = I$. Thus for any $0 < \eps < 1$, there exist $S, T \in \mathcal{M}$
such that $\|I - PS - QT\| \leq \eps$. But then $PS + QT$ is invertible, which implies that
$H = \ran(PS + QT) \subset \ran(P) + \ran(Q) = H_1 + H_2$, a contradiction.
\end{proof}

Of course, in the presence of commutativity, the 
need for the closure disappears. Namely,

\begin{proposition} \label{IV.D.5}
Let $X$ be an operator space.
\begin{enumerate}
\item[(i)] If $J$ and $K$ are right $M$-summands of $X$ (resp. left $L$-summands of $X$) whose corresponding
left $M$-projections (resp. right $L$-projections) commute, then $J + K$ is a right $M$-summand of $X$
(resp. left $L$-summand of $X$). In particular, $J + K$ is closed.
\item[(ii)] If $J$ and $K$ are right $M$-ideals of $X$, then $J + K$ is a right $M$-ideal of $X$ if and only
if $J + K$ is closed. This happens, for example, when the left $M$-projections corresponding
to $J^{\perp\perp}$ and $K^{\perp\perp}$ (resp. the right $L$-projections corresponding to
$J^\perp$ and $K^\perp$) commute.
\end{enumerate}
\end{proposition}

\begin{proof}
(i) Suppose $P$ corresponds to $J$ and $Q$ corresponds to $K$. Then $P + Q - PQ$ corresponds to
$J + K$. 

(ii) The first assertion follows from Corollary \ref{IV.D.2}. The second assertion follows from (i)
above, with the help of Lemma \ref{App.A.4}.
\end{proof}

As the next example shows, the omission in Corollary \ref{IV.D.2} of a statement concerning the intersection
of right $M$-ideals was not an oversight---no general statement can be made.
The lack of closure under arbitrary intersections of $M$-ideals 
is a feature of the classical $M$-ideal
theory (cf. \cite{HWW}, Example II.5.5).  However, 
the lack of closure under finite intersections is a new feature of the
one-sided theory.

\begin{example} \label{IV.D.6}
There exists a unital operator algebra $\mathcal{B}$ and projections $\tilde{P}, \tilde{Q} \in \mathcal{B}$
such that $\tilde{P}\mathcal{B} \cap \tilde{Q}\mathcal{B}$ has no contractive left approximate identity.
Thus, by facts mentioned in
Section \ref{III.E}, the intersection of two right $M$-summands of an operator space $X$ need not be 
a right $M$-ideal of $X$, let alone a right $M$-summand. Consequently, the intersection of two right
$M$-ideals of an operator space need not be a right $M$-ideal.
\end{example}

\begin{proof}
Let $\mathcal{A} \subset B(H)$ 
be the unital $C^*$-algebra constructed in Example 7.1 of \cite{BSZ}, the important
features of which we now recall:
\begin{itemize}
\item $\mathcal{A}$ contains projections $P$ and $Q$ such that $\ran(P \wedge Q)$ is separable with orthonormal basis 
$\{x_1, x_2, ...\}$, and $P \wedge Q \notin \mathcal{A}$.
\item The one-dimensional projection $E_n$ with range $\spn\{x_n\}$ is an element of $\mathcal{A}$ for all 
$n \in \mathbb{N}$.
\end{itemize}
Set $\mathcal{B} = R_\infty^w(C_\infty(\mathcal{A}))$. Then $\mathcal{B}$ is a unital operator algebra
(cf. Appendix \ref{App.B}). Set $\tilde{P} = P \otimes I_\infty$ and $\tilde{Q} = Q \otimes I_\infty$. 
Clearly, $\tilde{P}$ and $\tilde{Q}$ are projections in $\mathcal{B}$. Suppose $\{F_\lambda: \lambda \in \Lambda\}$ 
is a contractive left approximate identity for
\[
	J = \tilde{P}\mathcal{B} \cap \tilde{Q}\mathcal{B} = 
		R_\infty^w(C_\infty(P\mathcal{A})) \cap 
		R_\infty^w(C_\infty(Q\mathcal{A})) = R_\infty^w(C_\infty(P\mathcal{A} \cap 
		Q\mathcal{A})).
\]
Then
\[
	F_\lambda
	\begin{bmatrix}
		E_1 & E_2 & E_3 & \dots\\
		0 & 0 & 0 & \dots\\
		0 & 0 & 0 & \dots\\
		\vdots & \vdots & \vdots & \ddots
	\end{bmatrix}
	\to
	\begin{bmatrix}
		E_1 & E_2 & E_3 & \dots\\
		0 & 0 & 0 & \dots\\
		0 & 0 & 0 & \dots\\
		\vdots & \vdots & \vdots & \ddots
	\end{bmatrix}.
\]
We may view this as a convergent net in $M_\infty^w(B(H))$. Multiplying on the right by
\[
	\begin{bmatrix}
		E_1 & 0 & 0 & \dots\\
		E_2 & 0 & 0 & \dots\\
		E_3 & 0 & 0 & \dots\\
		\vdots & \vdots & \vdots & \ddots
	\end{bmatrix},
\]
we see that
\[
	F_\lambda
	\begin{bmatrix}
		P \wedge Q & 0 & 0 & \dots\\
		0 & 0 & 0 & \dots\\
		0 & 0 & 0 & \dots\\
		\vdots & \vdots & \vdots & \ddots
	\end{bmatrix}
	\to
	\begin{bmatrix}
		P \wedge Q & 0 & 0 & \dots\\
		0 & 0 & 0 & \dots\\
		0 & 0 & 0 & \dots\\
		\vdots & \vdots & \vdots & \ddots
	\end{bmatrix}.
\]
From this it follows that $F_\lambda(1,1)(P \wedge Q) \to (P \wedge Q)$. Using Lemma 2.9 in \cite{OpAlg}, we deduce 
that\\ $F_\lambda(1,1)^* = F_\lambda(1,1)^* (P \wedge Q) \to (P \wedge Q)$. Thus $P \wedge Q \in \mathcal{A}$, which 
is a contradiction.
\end{proof}

There are a number of ways to alleviate this lack of closure under 
finite intersections.   They all hinge on the 
following lemma.

\begin{lemma} \label{IV.D.7}
Let $X$ be an operator space, and $J$ and $K$ be right $M$-ideals of $X$. If $J + K$ is norm closed (equivalently
if $J^\perp + K^\perp$ is norm closed or weak-* closed), then $J \cap K$ is a right $M$-ideal of $X$.
\end{lemma}

\begin{proof}
By Lemma \ref{App.A.4}, the fact that $J + K$ is norm closed is equivalent to the fact that $J^\perp + K^\perp$
is norm closed, which in turn is equivalent to the fact that $J^\perp + K^\perp$ is weak-* closed. Now suppose that
 $J + K$ is norm closed. Then by Lemma \ref{App.A.2},
\[
	(J \cap K)^\perp = \overline{J^\perp + K^\perp}^{\wks} = J^\perp + K^\perp = \overline{J^\perp + K^\perp},
\]
which is a left $L$-summand of $X^*$ by Proposition \ref{IV.D.1}.
\end{proof}

The following results are stated for two right $M$-ideals;
 we leave the case of a finite number of right $M$-ideals as an exercise.

One situation when we can apply Lemma \ref{IV.D.7} is when there is `commutativity'. Namely,

\begin{proposition} \label{IV.D.8}
Let $X$ be an operator space.
\begin{enumerate}
\item[(i)] If $J$ and $K$ are right $M$-summands of $X$ whose corresponding left $M$-projections commute, then 
$J \cap K$ is a right $M$-summand of $X$.
\item[(ii)] If $J$ and $K$ are right $M$-ideals of $X$ such that the left $M$-projections corresponding
to $J^{\perp\perp}$ and $K^{\perp\perp}$ (resp. the right $L$-projections corresponding to
$J^\perp$ and $K^\perp$) commute, then $J \cap K$ is a right $M$-ideal of $X$.
\end{enumerate}
\end{proposition}

\begin{proof}
(i) Let $P$ and $Q$ be the left $M$-projections of $X$ corresponding to $J$ and $K$, respectively. Then
$PQ$ is the left $M$-projection of $X$ corresponding to $J \cap K$. 

(ii) 
This follows immediately from
Proposition \ref{IV.D.5} and Lemma \ref{IV.D.7}.
\end{proof}

A second situation when we can apply 
Lemma \ref{IV.D.7} is when there is `finite-dimensionality'. Namely,

\begin{corollary} \label{IV.D.9}
Let $X$ be an operator space, and $J$ and $K$ be right $M$-ideals of $X$. If $J$ is finite-dimensional or
finite-codimensional, then $J \cap K$ and $J + K$ are right $M$-ideals of $X$.
\end{corollary}

\begin{proof}
Suppose $J$ is finite-dimensional. Then $J + K$ is closed (Lemma \ref{App.A.7}). 
Thus $J \cap K$ is a right
$M$-ideal of $X$ by Lemma \ref{IV.D.7},
and $J + K$ is a right
$M$-ideal by Proposition \ref{IV.D.5} (ii). Indeed 
since $J \cap K$ is finite-dimensional, it
 is a right $M$-summand of $X$ (by a 
fact stated towards the end of
\S \ref{I.C}). Now suppose that
 $J$ is finite-codimensional. Then $J^\perp \cong (X/J)^*$ is finite-dimensional,
so that $J^\perp + K^\perp$ is closed. Again we appeal to Lemma \ref{IV.D.7}
and Proposition \ref{IV.D.5}.
\end{proof}

A third situation when we can apply Lemma \ref{IV.D.7} is when there is 
`perpendicularity':

\begin{corollary} \label{IV.D.9b}  Let $X$ be an operator space, and let
$J$ and $K$ be right $M$-ideals of $X$.  If $J \rightthreetimes K$ then 
$J + K$ and $J \cap K$ are right $M$-ideals of $X$. Of course, $J \cap K = \{0\}$ in this case.
\end{corollary} 

\begin{proof}  There exists a
completely isometric embedding $\sigma$ of $X$ into a
$C^*$-algebra satisfying (\ref{pere}) for all $x \in J, y \in K$.
Then $\sigma^{**}$ is a completely isometric embedding of $X^{**}$
into a
$W^*$-algebra.  Since $\sigma^{**}$ is weak-* continuous,
it is easy to check that $\sigma^{**}$  satisfies the 
analog of (\ref{pere}) for all $x \in J^{\perp \perp},
y \in K^{\perp \perp}$.  Let 
$P$ and $Q$ be the left $M$-projections
onto $J^{\perp \perp}$ and $K^{\perp \perp}$, respectively.
If 
$\langle \cdot , \cdot \rangle$ is the `left Shilov inner product' 
discussed in 
subsection \ref{II.BBB}, then we have by an observation there that
$$\langle Q P x , y \rangle = \langle P x , Q y \rangle = 0,$$ 
for all $x, y \in X^{**}$. Thus $Q P = 0 = P Q$, and we may appeal to 
Propositions \ref{IV.D.5} and  \ref{IV.D.8}.  \end{proof}   

A final situation when we can apply Lemma \ref{IV.D.7} is when there is 
`obliqueness'. Namely,

\begin{proposition} \label{IV.D.10}
Let $X$ be an operator space, and let
$J$ and $K$ be right $M$-ideals of $X$. Let $P$ and $Q$ be the
left $M$-projections corresponding to $J^{\perp\perp}$ and $K^{\perp\perp}$ (resp. the right $L$-projections
corresponding to $J^\perp$ and $K^\perp$). If $\|P(Q - P \wedge Q)\| < 1$, then 
$J \cap K$ and $J + K$ are right $M$-ideals of $X$.
\end{proposition}

Before proceeding with the proof, we interpret the condition appearing in the proposition. Consider orthogonal
projections $P$ and $Q$ on a Hilbert space $H$. Let $J$ and $K$ be the ranges of $P$ and $Q$, respectively.
Then one can easily show that $\|PQ\| < 1$ if and only if there exists $0 \leq \alpha < 1$ such that 
$|\langle x, y \rangle| \leq \alpha$ for all $x \in \ball(J)$ and all $y \in \ball(K)$. In particular, 
$J \cap K = \{0\}$. More generally, $\|P(Q - P \wedge Q)\| < 1$ if and only if there exists $0 \leq \alpha < 1$ such 
that $|\langle x, y \rangle| \leq \alpha$ for all $x \in \ball(J \cap (J \cap K)^\perp)$ and all
$y \in \ball(K \cap (J \cap K)^\perp)$. So as long as neither $J$ nor $K$ is contained in the other, the
condition $\|P(Q - P \wedge Q)\| < 1$ expresses the fact that the `angle' between $J$ and $K$ is positive (in
the Hilbert space case, at least).

Now we come to the proof of the proposition.

\begin{proof}
By Lemma \ref{IV.D.7}, it suffices to prove that $\ran(P) + \ran(Q)$ is closed. By Lemma \ref{App.A.9},
$\ran(P) + \ran(Q - P \wedge Q)$ is closed. But (claim) $\ran(P) + \ran(Q) = \ran(P) + \ran(Q - P \wedge Q)$.
Indeed,
\[
	Px + Qy = P(x + (P \wedge Q)y) + (Q - P \wedge Q)y
\]
and
\[
	Px + (Q - P \wedge Q)y = Px + Q(y - (P \wedge Q)y).
\]
\end{proof}

{\bf Remarks.} 1) Proposition \ref{IV.D.10} is a 
more general form of Theorem II.7 of \cite{Akemann69}, and 
indeed our result was inspired by
that Theorem.  We thank L. G. Brown
for drawing our attention to this result of Akemann.

2) Proposition \ref{IV.D.10} contains 
Proposition \ref{IV.D.8} and Corollary \ref{IV.D.9} as special cases, a fact which the 
energetic reader can confirm. \\
  
One question remains here: 
Is the intersection of two right $M$-ideals of $X^*$ again a right $M$-ideal of $X^*$?  We doubt it, but no counter-example comes to mind.

\subsection{Algebraic Direct Sum} \label{IV.E}

We turn to another item in the classical $M$-ideal `calculus'.
 Namely, if $X$ is a Banach space, and 
$J$ and $K$ are $M$-ideals of $X$ such that $X = J \oplus K$ (internal algebraic direct sum), 
then $J$ and $K$ are in fact 
complementary $M$-summands, i.e. $X = J \oplus_{\infty} K$ (see
\cite{Behrends}, Proposition 2.8). Furthermore, an
$M$-ideal $J$ of a Banach space $X$ is an $M$-summand of $X$ if and only if there exists an $M$-ideal $K$ of $X$
such that $X = J \oplus K$ (same reference). The corresponding statements are not true for operator spaces and 
one-sided $M$-ideals.

\begin{example} \label{IV.E.1}
Let $\mathcal{A} = C[0,1]$ and $\mathcal{I} = \{h \in C[0,1]: h(0) = 0\}$. Define $X = \mathcal{A} \oplus \mathcal{I}$
(external direct sum). Then $X$ is a left Hilbert
$C^*$-module over  $\mathcal{A}$. Set
\[
	J = \{(h,h): h \in \mathcal{I}\} \text{ and } K = \{(f,0): f \in \mathcal{A}\}.
\]
Then $J$ and $K$ are closed left submodules (i.e. left $M$-ideals) of $X$ such that $X = J \oplus K$ 
(internal direct sum). However, $J$ is not orthogonally complemented (i.e. not a left $M$-summand). Indeed, 
if $\langle (f,g),(h,h) \rangle = 0$ for all $h \in \mathcal{I}$, then $(f + g)\overline{h} = 0$ for all 
$h \in \mathcal{I}$. This yields $f + g = 0$ on $(0,1]$, which in turn yields that $f + g = 0$. 
Thus, $f \in \mathcal{I}$ and $(f,g) \in \mathcal{I} \oplus \mathcal{I}$. It follows that not every element of $X$ can 
be written as the sum of an element of $J$ and an element orthogonal to $J$.
\end{example}

It is precisely `commutativity' that is the missing ingredient (cf. Proposition \ref{IV.E.3} below). We begin with a 
simple lemma.

\begin{lemma} \label{IV.E.2}
Let $X$ be an operator space, and $J$ and $K$ be closed linear subspaces of $X$. Then:
\begin{enumerate}
\item[(i)] $J$ and $K$ are complementary right $M$-summands of $X$ if and only if $J^\perp$ and $K^\perp$ are
complementary left $L$-summands of $X^*$.
\item[(ii)] $J$ and $K$ are complementary left $L$-summands of $X$ if and only if $J^\perp$ and $K^\perp$ are
complementary right $M$-summands of $X^*$.
\end{enumerate}
\end{lemma}

\begin{proof}
(i) Suppose $J$ and $K$ are complementary right $M$-summands of $X$. Then there exists a left $M$-projection
$P$ on $X$ such that $\ran(P) = J$ and $\ker(P) = K$. Then $P^*$ is a right $L$-projection on $X^*$ such that
$\ker(P^*) = \ran(P)^\perp = J^\perp$ and $\ran(P^*) = \ker(I - P^*) = \ran(I - P)^\perp = \ker(P)^\perp = K^\perp$.
Thus, $J^\perp$ and $K^\perp$ are complementary left $L$-summands of $X^*$. Conversely, suppose that $J^\perp$ and
$K^\perp$ complementary left $L$-summands of $X^*$. Then there exists a right $L$-projection $Q$ on $X^*$ such
that $\ran(Q) = J^\perp$ and $\ker(Q) = K^\perp$. By Lemma \ref{App.A.6}, $Q$ is weak-* continuous. Thus, there
exists a bounded projection $P$ on $X$ such that $P^* = Q$. Clearly, $P$ is a left $M$-projection. Furthermore,
$\ran(P) = \ker(I - P) = \phantom{}^\perp\ran(I - Q) = \phantom{}^\perp\ker(Q) = \phantom{}^\perp(K^\perp) =
\overline{K} = K$ and $\ker(P) = \phantom{}^\perp\ran(Q) = \phantom{}^\perp(J^\perp) = \overline{J} = J$. Thus,
$J$ and $K$ are complementary right $M$-summands of $X$. 

(ii) The proof of this assertion is similar to that of
the previous one (easier, in fact).
\end{proof}

\begin{proposition} \label{IV.E.3}
Let $X$ be an operator space, and $J$ and $K$ be right $M$-ideals of $X$ such that $X = J \oplus K$. Then $J$ and
$K$ are complementary right $M$-summands of $X$ if and only if the left $M$-projections corresponding to 
$J^{\perp\perp}$ and $K^{\perp\perp}$ (resp. the right $L$-projections corresponding to $J^\perp$ and $K^\perp$) 
commute.
\end{proposition}

\begin{proof}
Let $P$ and $Q$ be the right $L$-projections corresponding to $J^\perp$ and $K^\perp$, respectively. Suppose that
$J$ and $K$ are complementary right $M$-summands of $X$. Then by Lemma \ref{IV.E.2}, $J^\perp$ and
$K^\perp$ are complementary left $L$-summands of $X^*$. Thus, $Q = I - P$, so that $PQ = QP = 0$. Conversely,
suppose that $PQ = QP$. Since
 $X = J \oplus K$, $X^* = J^\perp \oplus K^\perp$ (Lemma \ref{App.A.5}). 
But then
$PQ = QP$ has range $J^\perp \cap K^\perp = \{ 0 \}$.
Consequently, $P + Q = P + Q - PQ = I$.
Thus $J^\perp$ and $K^\perp$ are complementary left
$L$-summands of $X^*$. By Lemma \ref{IV.E.2}, $J$ and $K$ are complementary right $M$-summands of $X$.
\end{proof}

\begin{corollary} \label{IV.E.4}
Let $X$ be an operator space, and $J$ be a right $M$-ideal of $X$. Then $J$ is a right $M$-summand of $X$ if and only
if there exists a right $M$-ideal $K$ of $X$ such that $X = J \oplus K$ and the left $M$-projections corresponding
to $J^{\perp\perp}$ and $K^{\perp\perp}$ (resp. the right $L$-projections corresponding to $J^\perp$ and $K^\perp$)
commute.
\end{corollary}

\subsection{One-Sided $M$-Summands in Tensor Products} \label{IV.G0}

We give a method to verify that if $P$ is a 
left $M$-projection on an operator space $X$, 
then $P \otimes \Id_Z$ is  a
left $M$-projection on $X \otimes_\beta Z$, and that
$P(X) \otimes_\beta Z$ is the corresponding right $M$-summand,
for a fairly general class of
`tensor products' $\otimes_\beta$.
Indeed, the argument sketched in the next paragraph applies to any 
tensor product $\otimes_\beta$ with the following three
properties:
\begin{itemize}
\item [(1)]  $- \otimes_\beta \Id_Z$ is {\em functorial}.  That is,
  if $T : X_1 \to X_2$ is completely contractive,
then $T \otimes \Id_Z : X_1 \otimes_\beta Z \to X_2 \otimes_\beta Z$ is completely contractive.
\item [(2)]  The canonical map $C_2(X) \otimes Z \to C_2(X \otimes Z)$ extends to a completely
isometric isomorphism $C_2(X) \otimes_\beta Z \cong C_2(X \otimes_\beta Z)$.  
\item [(3)]  The span of
elementary tensors $x \otimes z$ for $x \in X, z \in Z$ is 
dense in $X \otimes_\beta Z$.
\end{itemize} 
Properties (1)--(3) certainly hold, for example, for
the minimal (i.e. spatial),  Haagerup, or extended  Haagerup
tensor products of operator spaces.  The argument also applies to
the module Haagerup tensor product over an
operator
algebra $\mathcal{A}$, provided $P$ is also a right ${\mathcal A}$-module
map.   

Suppose then that $P$ is a left $M$-projection on
$X$.  By Theorem \ref{left_M} (v), this is equivalent to saying 
that the maps $\nu_P^c$ and $\mu_P^c$ mentioned there give a
sequence of complete contractions
$$X \; \overset{\nu_P^c}{\longrightarrow} C_2(X) \;  \overset{\mu_P^c}{\longrightarrow} X. $$ 
Applying hypothesis (1) above,
 we obtain an induced sequence of complete contractions
$$X \otimes_\beta Z \; \overset{\nu_P^c \otimes \Id_Z}{\longrightarrow}
\; C_2(X) \otimes_\beta Z \; \overset{\mu_P^c \otimes \Id_Z}{\longrightarrow}
\; X \otimes_\beta Z .$$
From (2) 
we obtain a sequence 
$$X \otimes_\beta Z \; \longrightarrow \; C_2(X \otimes_\beta Z)
\; \longrightarrow X \otimes_\beta Z.$$
From (3) it is easy to
see that the first map in the last sequence equals $\nu^c_{P \otimes \Id_Z}$,
while the second equals $\mu^c_{P \otimes \Id_Z}$.   Thus these are 
complete contractions. Also, by (3), $P \otimes \Id_Z$ is idempotent. Thus,
by Theorem \ref{left_M} (v), we deduce that
$P \otimes \Id_Z$ is a  left $M$-projection on $X \otimes_\beta Z$.
By (1) again the canonical maps
$$P(X) \otimes_\beta Z \; \longrightarrow \; 
X \otimes_\beta Z \; \overset{P \otimes \Id_Z}{\longrightarrow} \;
P(X) \otimes_\beta Z$$
are complete contractions, which by (3) compose to the identity map.
Hence the first map in the last sequence is a complete isometry.
We conclude that $P(X) \otimes_\beta Z$ may be identified 
completely isometrically with a closed linear subspace of $X \otimes_\beta Z$.
It is easy to see that 
this subspace is exactly
 the image of the left $M$-projection $P \otimes \Id_Z$.
   
In  certain cases, a modification of the argument above
will show that right $M$-ideals in $X$ give 
rise to right $M$-ideals in $X \otimes_\beta Z$.  The following
proof is a good illustration of this technique.

\begin{theorem} \label{IV.F.f}  If $J$ is a right $M$-ideal in an
operator space $X$, then $J \otimes_h Y$ is a right $M$-ideal in
$X \otimes_h Y$ for any operator space $Y$.
\end{theorem} 

\begin{proof} \ (Sketch) \   Suppose that $P$ is the corresponding left 
$M$-projection from $X^{**}$ onto $J^{\perp \perp}$.   
We need to show that $(J \otimes_h Y)^{\perp \perp}$ is a
right $M$-summand in $(X \otimes_h Y)^{**}$.
The latter space may be identified with the {\em normal 
Haagerup tensor product} $X^{**} \otimes_{\sigma h} Y^{**}$.
See \cite{ERoc} for the definition and basic properties of this
tensor product.  In particular it satisfies
analogues of properties (1)--(3) above;  namely
\begin{itemize}
\item [(1)]  $- \otimes_{\sigma h} \Id_{Y^{**}}$ 
is {\em functorial}.  That is,
  if $T : X_1 \to X_2$ is weak-* continuous and completely contractive,
then $T \otimes \Id_{Y^{**}} : X_1 \otimes_{\sigma h} Y^{**}
 \to X_2 \otimes_{\sigma h} Y^{**}$ is 
weak-* continuous and  completely contractive.
\item [(2)]  $C_2(X^{**}) \otimes_{\sigma h} Y^{**} \cong
C_2(X^{**} \otimes_{\sigma h} Y^{**})$ completely
isometrically and weak-* homeomorphically.
\item [(3)]  The span of
elementary tensors $x \otimes y$ for $x \in X, y \in Y$ is
weak-* dense in $X^{**} \otimes_{\sigma h} Y^{**}$.
\end{itemize}   
To see (2) for example,
note that by basic properties of the Haagerup 
tensor product, $$C_2(X^{**} \otimes_{\sigma h} Y^{**})
\cong C_2((X \otimes_h Y)^{**}) \cong C_2(X \otimes_h Y)^{**}
\cong C_2(X^{**}) \otimes_{\sigma h} Y^{**},$$ 
and one can verify that each `$\cong$' here is
a weak-* homeomorphism.

We may now follow the proof above Theorem \ref{IV.F.f}, 
checking that all maps
are weak-* continuous.  We deduce that $P \otimes \Id_{Y^{**}}$
is a left $M$-projection on $X^{**} \otimes_{\sigma h} Y^{**}$,
and its range is completely isometric and weak-* homeomorphic
to $J^{\perp \perp} \otimes_{\sigma h} Y^{**}$.  
Since $P \otimes \Id_{Y^{**}}$ is
the identity map when restricted to
$J \otimes_h Y$, it is easy to see that
the range of  $P \otimes \Id_{Y^{**}}$ is 
$(J \otimes_h Y)^{\perp \perp}
= \overline{J \otimes_h Y}^{\wks}$.
 \end{proof} 
 
We now turn to different methods to analyze tensor products.    

\subsection{Minimal Tensor Product} \label{IV.G} 

We write $X \itimes Y$ for the minimal or 
spatial tensor product of operator spaces
(see the standard operator space texts for more details).

A similar analysis to that in the last proof shows that in cases
that $(X \itimes Y)^{**} \cong
X^{**} \bar{\otimes} Y^{**}$ (see e.g.\ 
\cite[Chapter 8]{ERbook} for this notation),
if $J$ is a right $M$-ideal in
 $X$ then $J \itimes Y$ is a right $M$-ideal in
$X \itimes Y$.  We leave it to the reader to make 
this precise.
Note however that $(X \itimes Y)^{**} \neq 
X^{**} \bar{\otimes} Y^{**}$ in many cases, even if 
$X$ is finite dimensional---this is related to the topic
of local reflexivity.
                 
\begin{proposition} \label{IV.G.1}
Let $X$ and $Y$ be operator spaces, $S \in \Ml(X)$, and $T \in \Ml(Y)$. Then $S \otimes T \in \Ml(X \itimes Y)$, and
$\|S \otimes T\|_{\Ml(X \itimes Y)} \leq \|S\|_{\Ml(X)}\|T\|_{\Ml(Y)}$. If $S = \Id_X$ or $T = \Id_Y$, then in fact we
have equality in the previous inequality. If $S \in \Al(X)$ and $T \in \Al(Y)$, then $S \otimes T \in \Al(X \itimes Y)$
and $(S \otimes T)^\star = S^\star \otimes T^\star$.
\end{proposition}

\begin{proof}
Let $(\sigma,A)$ and $(\rho,B)$ be implementing pairs for $S$ and $T$, respectively. Then $\sigma:X \to B(H)$ and
$\rho:Y \to B(K)$ are complete isometries. Thus, by the injectivity of the minimal tensor product,
$\sigma \otimes \rho:X \itimes Y \to B(H) \itimes B(K) \subset B(H \otimes K)$ is a complete isometry. We claim that
$(\sigma \otimes \rho,A \otimes B)$ is an implementing pair for $S \otimes T$. Indeed, for all $x \in X$ and
$y \in Y$,
\begin{eqnarray*}
	(\sigma \otimes \rho)((S \otimes T)(x \otimes y)) &=& (\sigma \otimes \rho)(Sx \otimes Ty) =
		\sigma(Sx) \otimes \rho(Ty)\\
	&=& A\sigma(x) \otimes B\rho(y) = (A \otimes B)(\sigma(x) \otimes \rho(y))\\
	&=& (A \otimes B)(\sigma \otimes \rho)(x \otimes y).
\end{eqnarray*}
It follows that $S \otimes T \in \Ml(X \itimes Y)$ and 
$\|S \otimes T\|_{\Ml(X \itimes Y)} \leq \|A \otimes B\| = \|A\|\|B\|$. Since the choice of implementing pairs was 
arbitrary, $\|S \otimes T\|_{\Ml(X \itimes Y)} \leq \|S\|_{\Ml(X)}\|T\|_{\Ml(Y)}$. Furthermore, if
$A^*\sigma(X) \subset \sigma(X)$ and $B^*\rho(Y) \subset \rho(Y)$, then 
$(A^* \otimes B^*)(\sigma \otimes \rho)(X \itimes Y) \subset (\sigma \otimes \rho)(X \itimes Y)$. Thus, if
$S \in \Al(X)$ and $T \in \Al(Y)$, then $S \otimes T \in \Al(X \itimes Y)$ and 
$(S \otimes T)^\star = S^\star \otimes T^\star$. Finally, suppose that
$T = \Id_Y$. Let $(\theta,C)$ be an implementing
pair for $S \otimes T$, so that $\theta:X \itimes Y \to B(L)$ is a complete isometry. Fix a unit vector $y \in Y$
and define $\theta_y:X \to B(L)$ by $\theta_y(x) = \theta(x \otimes y)$ for all $x \in X$. Then $\theta_y$ is a
complete isometry and
\[
	\theta_y(Sx) = \theta(Sx \otimes y) = \theta(Sx \otimes Ty) = C\theta(x \otimes y) = C\theta_y(x)
\]
for all $x \in X$. In other words, $(\theta_y,C)$ is an implementing pair for $S$. Thus, 
$\|S\|_{\Ml(X)}\|T\|_{\Ml(Y)} = \|S\|_{\Ml(X)} \leq \|C\|$. Since the choice of the implementing pair was
arbitrary, $\|S\|_{\Ml(X)}\|T\|_{\Ml(Y)} \leq \|S \otimes T\|_{\Ml(X \itimes Y)}$. The same type of argument works
if $S = \Id_X$.
\end{proof}

In particular, a left multiplier (resp. left adjointable multiplier) on an operator space $X$ amplifies to
a left multiplier (resp. left adjointable multiplier) on $M_I(X)$ or $K_I(X)$, for any index set $I$. The
amplification has the same multiplier norm as the original map. 
In fact, the same is true for the amplification
to $M_I^w(X)$, although this doesn't follow from the above proposition 
(it follows from a modification of the proof of
 Proposition \ref{IV.I.1} below; see also Corollary \ref{IV.I.f}).

\subsection{Haagerup Tensor Product} \label{IV.H}

\begin{proposition} \label{IV.H.1}  
If $T$ is a left multiplier on an operator space $X$, and if $Y$ is any other operator space, then $T \otimes \Id_Y$
is a left multiplier on $X \htimes Y$, where $\htimes$ is the Haagerup tensor product. Furthermore, the multiplier
norms of $T$ and $T \otimes \Id_Y$ are the same. If $T$ is a left adjointable multiplier on $X$, then
$T \otimes \Id_Y$ is a left adjointable multiplier of $X \htimes Y$, and $(T \otimes \Id_Y)^\star =
T^\star \otimes \Id_Y$.
\end{proposition}

\begin{proof}
Suppose $T \in \Ml(X)$ and $\|T\|_{\Ml(X)} \leq 1$. Under the completely isometric identification 
$C_2(X \htimes Y) \cong C_2(X) \htimes Y$, $\tau_{T \otimes \Id_Y}^c$ corresponds to $\tau_T^c \otimes \Id_Y$. 
Therefore,
\[
	\|\tau_{T \otimes \Id_Y}^c\|_{cb} = \|\tau_T^c \otimes \Id_Y\|_{cb} \leq \|\tau_T^c\|_{cb}\|\Id_Y\|_{cb} 
		\leq 1,
\]
which says that $T \otimes \Id_Y \in \Ml(X \htimes Y)$ and $\|T \otimes \Id_Y\|_{\Ml(X \htimes Y)} \leq 1$. Now 
suppose that
$T \otimes \Id_Y \in \Ml(X \htimes Y)$. Let $(\sigma,A)$ be an implementing pair for $T \otimes \Id_Y$. Fix 
a unit vector $y \in Y$ arbitrarily and define $\sigma_y:X \to B(H)$ by
\[
	\sigma_y(x) = \sigma(x \otimes y)
\]
for all $x \in X$. Then $\sigma_y$ is a complete isometry and
\[
	\sigma_y(Tx) = \sigma(Tx \otimes y) = \sigma((T \otimes \Id_Y)(x \otimes y)) = A\sigma(x \otimes y) = 
		A \sigma_y(x)
\]
for all $x \in X$. In other words, $(\sigma_y,A)$ is an implementing pair for $T$. Thus, $\|T\|_{\Ml(X)} \leq \|A\|$. 
Since the choice of $(\sigma,A)$ was arbitrary, $\|T\|_{\Ml(X)} \leq \|T \otimes \Id_Y\|_{\Ml(X \htimes Y)}$. Now 
suppose that
 $T \in \Al(X)_{sa}$. Then for all $t \in \mathbb{R}$,
\[
	\|\exp(it\tau_{T \otimes \Id_Y}^c)\|_{cb} = \|\exp(it(\tau_T^c \otimes \Id_Y))\|_{cb} = 
		\|\exp(it\tau_T^c) \otimes \Id_Y\|_{cb} \leq \|\exp(it\tau_T^c)\|_{cb} = 1.
\]
Thus, $T \otimes \Id_Y \in \Al(X \htimes Y)_{sa}$. Finally, if $T \in \Al(X)$, then $T = T_1 + iT_2$, where
$T_1, T_2 \in \Al(X)_{sa}$, so that
\[
	T \otimes \Id_Y = (T_1 \otimes \Id_Y) + i(T_2 \otimes \Id_Y) \in \Al(X \htimes Y).
\]
Furthermore,
\[
	(T \otimes \Id_Y)^\star = (T_1 \otimes \Id_Y) - i(T_2 \otimes \Id_Y) = T^\star \otimes \Id_Y.
\]
\end{proof}

{\bf Remarks.}  1) We observe that due to the lack of `commutativity' in the Haagerup tensor product, $\Id_X \otimes T$ need not be a 
left multiplier (resp. left adjointable multiplier) of $X \htimes Y$ when $T$ is a left multiplier (resp. 
left adjointable multiplier) of $Y$. The correct statement is that $\Id_X \otimes T$ is a right multiplier 
(resp. right adjointable multiplier) of $X \htimes Y$ when $T$ is a right multiplier (resp. right adjointable 
multiplier) of $Y$.

2)  We also note that results similar to Proposition \ref{IV.H.1} hold for the extended and module Haagerup tensor 
products.\\

It follows immediately from Proposition \ref{IV.H.1} that the map $T \mapsto T \otimes \Id_Y$ is an isometric 
homomorphism from $\Ml(X)$ into $\Ml(X \htimes Y)$, and that this map restricts to an injective $*$-homomorphism 
from $\Al(X)$ into $\Al(X \htimes Y)$. Thus $\Al(X)$ may be viewed as a $*$-subalgebra of $\Al(X \htimes Y)$. 
Typically this inclusion is proper. For example, $\Al(C_n) = M_n$, but $\Al(C_n\htimes C_n) = \Al(C_{n^2}) = M_{n^2}$. 
If $\mathcal{A}$ and $\mathcal{B}$  are nontrivial unital $C^*$-algebras, however, we have that 
$\Al(\mathcal{A}) = \Al(\mathcal{A} \htimes \mathcal{B})$ (see Theorem \ref{IV.H.3} below). 
First we recall  a well-known fact
 about the Haagerup tensor product,
for which we include a proof for completeness.

\begin{lemma} \label{IV.H.2}
Let $V$ and $W$ be operator spaces and $x \in V \htimes W$. If $(\phi \otimes \Id_W)(x) = 0$ for all $\phi \in V^*$, then $x = 0$. Likewise, if $(\Id_V \otimes \psi)(x) = 0$ for all $\psi \in W^*$, then $x = 0$.
\end{lemma}

\begin{proof}
We only prove the first assertion. By \cite{ERbook}, Proposition 9.2.5 and Theorem 9.4.7, we have the completely 
isometric inclusions
\[
	V \htimes W \to V^{**} \htimes W^{**} \to (V^* \htimes W^*)^*.
\]
Let $\Phi_x \in (V^* \htimes W^*)^*$ be the image of $x$ under these inclusions. 
We have
\[
	\Phi_x(\phi \otimes \psi) = (\phi \otimes \psi)(x) = \psi((\phi \otimes \Id_W)(x)) = 0
\]
for all $\phi \in V^*$ and all $\psi \in W^*$. It follows that $\Phi_x = 0$, which implies that $x = 0$.
\end{proof}

\begin{theorem} \label{IV.H.3}
Let $\mathcal{A}$ and $\mathcal{B}$ be unital $C^*$-algebras, 
with neither $\mathcal{A}$ nor $\mathcal{B}$ equal to 
$\mathbb{C}$. Then $\mathcal{A} \cong \Al(\mathcal{A} \htimes \mathcal{B})$ and 
$\mathcal{B} \cong \Ar(\mathcal{A} \htimes \mathcal{B})$.
\end{theorem}  

\begin{proof}
Let $\lambda:\mathcal{A} \to \Al(\mathcal{A} \htimes \mathcal{B})$ and 
$\rho:\mathcal{B} \to \Ar(\mathcal{A} \htimes 
\mathcal{B})$ be the canonical injective *-homomorphisms (see the discussion
above Lemma \ref{IV.H.2}). We aim to show 
that $\lambda$ is surjective (the proof that $\rho$ is surjective being similar). 
It suffices to show that any self-adjoint $T \in \Al(\mathcal{A} 
\htimes \mathcal{B})$, with 
$\|T\|_{\Ml(\mathcal{A} \htimes \mathcal{B})} \leq 1$, is in the range 
of $\lambda$. 
We will use R. R. Smith's
observation that $\mathcal{A} \htimes \mathcal{B}$ is a unital Banach algebra with product
\[
        (a \otimes b)(a' \otimes b') = aa' \otimes bb'
\]
(see e.g.  Proposition 2 in \cite{Bgeo}). 
From Lemma \ref{III.G.1} we know that 
$T(1 \otimes 1) \in \Her(\mathcal{A} \htimes \mathcal{B})$.
Thus by (1) on p. 126 of \cite{Bgeo} we have that 
\begin{equation} \label{ther}
	T(1 \otimes 1)  = h \otimes 1 + 1 \otimes k
\end{equation}
for some $h \in \mathcal{A}_{sa}$, $k \in \mathcal{B}_{sa}$. Since left and right multipliers of an operator space 
automatically commute, $\rho(\mathcal{B})$ commutes with $T$. Thus
\begin{equation} \label{tism}
	T(a \otimes b) = T(\rho(b)(a \otimes 1)) = \rho(b)(T(a \otimes 1)) = T(a \otimes 1) (1 \otimes b)
\end{equation}
for $a \in \mathcal{A}$, $b \in \mathcal{B}$. We next will prove the identity
\begin{equation} \label{tis}
	(\Id_{\mathcal{A}} \otimes \psi)(T(a \otimes w)) = (\Id_{\mathcal{A}} \otimes \psi)(T(1 \otimes w)) a  
\end{equation}
for all $\psi \in \mathcal{B}^*$ and $w \in \mathcal{B}$. First suppose that $w$ is a unitary in $\mathcal{B}$, and 
that $\psi \in \mathcal{B}^*$ satisfies $\psi(w) = 1 = \|\psi\|$. Consider the operator 
$u(a) = (\Id_{\mathcal{A}} \otimes \psi)(T(a \otimes w))$ on $\mathcal{A}$. We have for any $a' \in \mathcal{A}$ that 
\[
	\left\|\begin{bmatrix} u(a)\\ a' \end{bmatrix}\right\| =
	\left\|\begin{bmatrix} (\Id_{\mathcal{A}} \otimes \psi)(T(a \otimes w))\\ 
		(\Id_{\mathcal{A}} \otimes \psi)(a' \otimes w) \end{bmatrix}\right\| \leq
	\left\|\begin{bmatrix} T(a \otimes w)\\ a' \otimes w \end{bmatrix}\right\| \leq
	\left\|\begin{bmatrix} a \otimes w\\ a' \otimes w \end{bmatrix}\right\|,
\]
the last inequality by Theorem \ref{BEZ}. Since we clearly have
\[
	\left\|\begin{bmatrix} a \otimes w\\ a' \otimes w \end{bmatrix}\right\| = \left\|\begin{bmatrix}
		a\\ a' \end{bmatrix}\right\|,
\]
we see using Lemma 4.1 in \cite{BSZ}
that there exists an $a_{w,\psi} \in \mathcal{A}$ such that 
\[
	(\Id_{\mathcal{A}} \otimes \psi)(T(a \otimes w)) = a_{w,\psi} a
\]
for all $a \in \mathcal{A}$. Setting $a = 1$ gives $a_{w,\psi} = (\Id_{\mathcal{A}} \otimes \psi)(T(1 \otimes w))$, 
and this establishes (\ref{tis}) in this case.

If $g \in \mathcal{B}^*$, then $g_{w} = g(\cdot w) \in \mathcal{B}^*$. Since $\mathcal{B}^*$ is the span of the 
states on $\mathcal{B}$, we may write $g_{w} = \sum_{k=1}^4 \alpha_k f^k$ for $\alpha_k$ scalars and $f^k$ states on 
$\mathcal{B}$. Then $g = \sum_{k=1}^4 \alpha_k f^k_{w^*}$. Setting $\psi = f^k_{w^*}$ in (\ref{tis}), and using the 
fact that Equation (\ref{tis}) is linear in $\psi$, we now have (\ref{tis}) with $w$ unitary. By the well-known fact 
that the unitary elements span a $C^*$-algebra, and the linearity of Equation (\ref{tis}) in $w$, we obtain 
(\ref{tis}) for any $w \in \mathcal{B}$. Thus we have proved (\ref{tis}) in general.

Combining (\ref{tis}) and (\ref{tism}) we obtain
\[
	(\Id_{\mathcal{A}} \otimes \psi)(T(a \otimes b)) = (\Id_A \otimes \psi)(T(1 \otimes 1) (1 \otimes b)) a.
\]
Writing $T(1 \otimes 1)$ as in (\ref{ther}), we have that
\[
	(\Id_{\mathcal{A}} \otimes \psi)(T(a \otimes b)) = \psi(b)ha + \psi(kb)a = 
		(\Id_{\mathcal{A}} \otimes \psi)((h \otimes 1 + 1 \otimes k)(a \otimes b)).
\]
Since $\psi$ is arbitrary, we deduce by Lemma \ref{IV.H.3} that
\[
	T(a \otimes b) = (h \otimes 1 + 1 \otimes k) (a \otimes b)
\]
for all $a \in \mathcal{A}$, $b \in \mathcal{B}$. 

For any $x \in \mathcal{A} \htimes \mathcal{B}$, denote by $L_x$ the operator of left multiplication by $x$ on
$\mathcal{A} \htimes \mathcal{B}$. Since $L_{h \otimes 1} = \lambda(h)$, 
$L_{h \otimes 1} \in \Al(\mathcal{A} \htimes \mathcal{B})_{sa}$. Since 
$T \in \Al(\mathcal{A} \htimes \mathcal{B})_{sa}$ too, we deduce  that 
$L_{1 \otimes k} \in \Al(\mathcal{A} \htimes \mathcal{B})_{sa}$. We claim that $k$ is a scalar multiple of $1$. 
Suppose not. Let $c \in \mathcal{A}_{sa}$. Then 
$L_{c \otimes 1} = \lambda(c) \in \Al(\mathcal{A} \htimes \mathcal{B})_{sa}$ and $L_{c \otimes 1}$ commutes with 
$L_{1 \otimes k}$. It follows that 
$L_{c \otimes k} = L_{c \otimes 1}L_{1 \otimes k} \in \Al(\mathcal{A} \htimes \mathcal{B})_{sa}$. By Lemma 
\ref{III.G.1},
\[
	c \otimes k = L_{c \otimes k}(1 \otimes 1) \in \Her(\mathcal{A} \htimes \mathcal{B}),
\]
and so by (1) on p. 126 of \cite{Bgeo},
\[
	c \otimes k = h_c \otimes 1 + 1 \otimes k_c
\]
for some $h_c \in \mathcal{A}_{sa}$, $k_c \in \mathcal{B}_{sa}$. Now let $\phi \in \mathcal{A}^*$ be a state. Then
\[
	\phi(c)k = (\phi \otimes \Id_{\mathcal{B}})(c \otimes k) = 
		(\phi \otimes \Id_{\mathcal{B}})(h_c \otimes 1 + 1 \otimes k_c) = \phi(h_c) 1 + k_c,
\]
so that
\[
	k_c = \phi(c)k - \phi(h_c) 1.
\]
But then
\[
	c \otimes k = h_c \otimes 1 + 1 \otimes \phi(c)k - 1 \otimes \phi(h_c) 1, 
\]
which implies that
\[
	(c - \phi(c)1) \otimes k = (h_c - \phi(h_c)1) \otimes 1.
\]
Since $k$ and $1$ are linearly independent, $c = \phi(c)1$. 
Since the choice of $c$ was arbitrary, 
$\mathcal{A} = \mathbb{C}$, a contradiction.
\end{proof}

\begin{corollary} \label{IV.H.4}
Given any two 
unital $C^*$-algebras $\mathcal{A}$ and 
$\mathcal{B}$, there exists an operator space $X$ 
with $\Al(X) \cong \mathcal{A}$ and $\Ar(X) \cong \mathcal{B}$. Given any two lattices $L_1$ and $L_2$ 
which are the projection lattices for two von Neumann algebras, there exists an operator space $X$ whose 
lattice of left $M$-summands equals $L_1$ and whose lattice of right $M$-summands equals $L_2$.
\end{corollary}

\begin{proof}
If $\mathcal{A}$ and
$\mathcal{B}$ are both nontrivial then this follows from the 
last Theorem immediately.  So suppose that  $\mathcal{A} = 
\mathbb{C}$ (we leave the remaining cases to the reader).
Let $X = R_n \htimes \mathcal{B}$.    If $T \in 
\Al(X)$, then by Proposition \ref{IV.G.1} we have that
$\Id_{C_n} \otimes T \in \Al(C_n \itimes X)$.
By basic operator space theory we have
 $$C_n \itimes X \cong C_n \htimes R_n \htimes \mathcal{B}
\cong M_n \htimes \mathcal{B}.$$ 
By the last Theorem it follows that there exists a matrix 
$A \in M_n$, with $$x \otimes T(u) = Ax \otimes u$$
for all $x \in C_n, u \in X$. By linear algebra
we may deduce that $T \in \mathbb{C} \Id_X$. Since the choice of $T$ was arbitrary,
$\Al(X) = \mathbb{C}$. On the other hand, by Proposition
\ref{IV.H.1},
\[
	\mathcal{B} \subset \Ar(R_n \htimes \mathcal{B}) \subset \Ar(M_n \htimes \mathcal{B}) = \mathcal{B}.
\]
\end{proof} 

\subsection{Interpolation} \label{IV.F}
 
In this section we shall see that the class of one-sided $M$-summands
is closed under interpolation, and we will use this
to construct some `exotic' examples.
We will use the interpolation formula
$M_n(X_\theta) = M_n(X)_\theta$ and Proposition 2.1
from \cite{OH}.  Suppose that $(E_0,E_1)$ is a compatible couple of
operator spaces in the sense of interpolation theory.
If $T : E_0 + E_1 \to E_0 + E_1$ restricts to
a map $E_0 \to E_0$ which is a contractive left multiplier on $E_0$,
and restricts to
a map $E_1  \to E_1$ which is a contractive left multiplier on $E_1$,
then we claim that $T$ induces a left multiplier on $E_\theta$,
for all $\theta \in [0,1]$.   To see this
note that $\tau_T^c : C_2(E_0) + C_2(E_1) \to C_2(E_0) + C_2(E_1)$
restricts to complete contractions  $\tau_T^c : C_2(E_0)
\to C_2(E_0)$ and  $\tau_T^c : C_2(E_1)
\to C_2(E_1)$.  By \cite{OH} Proposition 2.1,
$\tau_T^c : C_2(E_\theta) \to C_2(E_\theta)$ is a complete
contraction.  Thus $T$ has induced a canonical
contractive left multiplier on each $E_\theta$.
Furthermore, if $T$ is idempotent, then we have obtained
a left $M$-projection on each $E_\theta$, and have
therefore interpolated between the left $M$-summands $T(E_0)$ and
$T(E_1)$.
It should be remarked that $T$ need not be a left multiplier
on $E_0 + E_1$. Indeed, $\tau_T^c$ is completely contractive on $C_2(E_0) + C_2(E_1)$,
but not necessarily on $C_2(E_0 + E_1)$.
 
\begin{example} \label{exin} Let $E_0 = \ell^\infty$, and $E_1 = \ell^2$ with
its column Hilbert space structure.   View $\ell^2 \subset \ell^\infty$.
Let $I$ be a nontrivial subset
of the natural numbers, and let  $F_0$ and $F_1$ be the subspaces
of $E_0 $ and $E_1$ supported on $I$.  These are
right $M$-summands in $E_0 $ and $E_1$ respectively.  Interpolating
between these spaces yields the somewhat surprising existence
of right $M$-summands in $\ell^p$ (the latter Banach space equipped
with some operator space structure) for
$2 < p < \infty$.    On the other hand
it follows from (4) in \cite{BEZ},
in conjunction with Theorem \ref{left_M} (iv),
that for $1 \leq p < 2$, there is no
operator space structure on
$\ell^p$ such that the sequences supported on $I$ form a
right $M$-summand.
\end{example}
 
\begin{example} \label{exin2}   A similar technique works
for $L^p$-spaces and their noncommutative variant.
In the commutative case note that $L^\infty(\Omega)$ acts
as left multipliers both on $L^2(\Omega)_c$ and on
$L^\infty(\Omega)$.   Interpolating, we see that
left multiplication by a function $f$ in $L^\infty(\Omega)$
is a left multiplier on $L^p(\Omega)$ (the latter Banach space equipped
with some operator space structure) for $2 < p < \infty$.
Taking $f$ to be a characteristic function yields left $M$-projections.
 
In the noncommutative case we do the same thing, after observing that a
(finite) von Neumann algebra $\mathcal{M}$  acts
as left multipliers both on the Hilbert space
$L^2(\mathcal{M})_c$, and on $L^\infty(\mathcal{M}) = \mathcal{M}$.
Interpolating we obtain exotic right $M$-summands and
left multipliers on the Banach space
$L^p(\mathcal{M})$ (with some operator space structure) for $2 < p < \infty$.
 \end{example}
 
The operator space structure on $L^p(M)$ in the last examples
is not quite the standard one of \cite{Pisier}.  However
it is not a very obscure structure.  Marius Junge has
informed us that he has encountered it in some of his work.
  In certain cases this structure has a simple description.
For example if $M = M_n$ then we are interpolating between
$M_n$ with its Hilbert-Schmidt norm and column
Hilbert space structure, and $M_n$.
That is,
in operator space notation,
 we are interpolating between $C_n \otimes_h C_n$ and
$C_n \otimes_h R_n$.
By a formula in \cite{OH}, we have
$$(C_n \otimes_h C_n,C_n \otimes_h R_n)_\theta \cong
C_n \otimes_h (C_n,R_n)_\theta \cong C_n(O\ell^p_n)
\cong CB(R_n,O\ell^p_n),$$
where $O\ell^p_n$ is the standard operator space structure on $\ell^p_n$
from \cite{Pisier}.  Thus the operator space structure which
we are putting on $L^p(M_n)$ is simply described as
$C_n(O\ell^p_n)$.  Equivalently, we are viewing
the Schatten class operators on $\ell_n^2$ as endowed with the
$CB(R_n,O\ell^p_n)$ operator space structure.
From this fact,
as pointed out to us by Junge, the left adjointable multiplier algebra
is easily computable.  By fact (\ref{I.A.4}) and the discussion at
the end of
Subsection \ref{IV.A}, we have
$\Al(C_n(O\ell^p_n)) \cong M_n(\Al(O\ell^p_n)) \cong M_n$,
if $p < \infty$.
From this and facts in the previous section,
it is easy to see that the only right $M$-ideals in
our interpolated space are of the form $J \otimes_h O\ell^p_n$,
for a subspace $J$ of $C_n$.

\subsection{Infinite Matrices and Multipliers} \label{IV.I}

We have the isometric identification $M_n(\Ml(X)) = \Ml(C_n(X))$
(see formula (\ref{I.A.4}) from 
Section \ref{I.A}).  For an infinite index set however,  $M_I^w(\Ml(X))$ 
 usually differs from $\Ml(C_I^w(X))$, even in the case that $X$ is a 
unital $C^*$-algebra (see Appendix \ref{App.B} for the relevant definitions).
For dual spaces the situation is much better: 
 
\begin{theorem} \label{IV.I.0}  Let  $X$ be a dual 
operator space, and let $I, J$ be index sets. Then
\begin{itemize} \item [(i)] $\Al(C^w_I(X)) \cong M^w_I(\Al(X)) \cong 
\Al(M^w_I(X))$ as von Neumann
algebras.
\item [(ii)] $\Al(M^w_{I,J}(X)) \cong M^w_{I}(\Al(X))$  as von Neumann
algebras.
\item [(iii)]  $\Ml(M^w_{I,J}(X)) \cong M^w_I(\Ml(X))$    as dual
operator algebras (that is, via a
completely isometric homomorphism which is a weak-* homeomorphism).  
\end{itemize}   \end{theorem} 

Note that this result is only new if $I, J$ are infinite.
In this case, the precise isomorphisms here are  
essentially the same as in the case that  $I, J$ are finite.

\begin{proof}  (i)  Suppose that $T \in \Al(C^w_I(X))$.
Fix $i , j \in I$, and let $Q_i$ and $U_{ij}$ be respectively the
projection onto the $i$th entry, and the 
permutation of the $i$ and $j$ entries. We have
$Q_i U_{ij} = U_{ij} Q_j$.  
  It is easy to see that $Q_i, U_{ij}$ are in $\Al(C^w_I(X))$,
and hence so is $S_{ij} = Q_j U_{ij} T Q_j$.
 Writing $\epsilon_i$ for the
inclusion of $X$ into $C^w_I(X)$ as the $i$th entry,
we have $S_{ij}(\epsilon_j(X)) \subset \epsilon_j(X)$,
and $S_{ij}^\star(\epsilon_j(X)) \subset \epsilon_j(X)$.
Thus $S_{ij}$ may be regarded as an element $T_{ij}$ of
$\Al(X)$.   Indeed $T_{ij} = \pi_i \circ T \circ \epsilon_j$,
where $\pi_i : C^w_I(X) \to X$ is the $i$th coordinate 
function.  We let $\theta(T) = \begin{bmatrix} T_{ij} \end{bmatrix}$, and we will proceed
to show that this is in $M^w_I(\Al(X))$. 

We introduce some notation: for a finite subset 
$\Delta$ of $I$ we define $P_\Delta = \sum_{i \in \Delta} Q_i$,
the projection of $C^w_I(X)$ onto the space of columns 
`supported on $\Delta$'.  For $T \in \Al(C^w_I(X))$,
let $T^\Delta = P_\Delta T P_\Delta$.
Also, given any $I \times I$ matrix $x$, 
write $x^\Delta$ for the same matrix but with entries 
$x_{ij}$ replaced by $0$ if $i$ or $j$ is not in $\Delta$.
It is then easy to see that 
$$\theta(T^\Delta) = \theta(T)^\Delta.$$
Let $A_\Delta \in M_{\Delta}(\Al(X))$ be the 
`$\Delta$-submatrix' of $\begin{bmatrix} T_{ij} \end{bmatrix}$.  That is,
$A_\Delta$ is the `possibly nonzero part' of $\theta(T)^\Delta.$
We identify $C_{\Delta}(X)$ with a
subspace of $C^w_I(X)$, and notice that this 
subspace is invariant under
$T^\Delta$ and its adjoint.  We write $T^\Delta|_{C_{\Delta}(X)}$
for the restriction of $T^\Delta$ to this subspace;
note that this restriction is in $\Al(C_{\Delta}(X))$ by
Lemma \ref{IV.B.4b}, and has the same norm there.    
 Since $\Al(C_n(X)) \cong M_n(\Al(X))$ we have
$$\Vert A_\Delta \Vert_{M_{\Delta}(\Al(X))}
= \Vert T^\Delta|_{C_{\Delta}(X)} \Vert_{\Al(C_{\Delta}(X))}
= \Vert T^\Delta 
\Vert_{\Al(C^w_I(X))} \leq \Vert T \Vert_{\Al(C^w_I(X))}.$$
From this we deduce that $\theta(T) = \begin{bmatrix} T_{ij} \end{bmatrix} \in 
M^w_I(\Al(X))$.  Also we now see that $\theta$ is
a unital linear contraction which is isometric 
on $P_\Delta \Al(C^w_I(X)) P_\Delta$.

Claim:  $\theta$ is weak-* continuous.  To see this
take a bounded net $(T_\lambda)$ converging to
$T$ with respect to the weak-* topology on $\Al(C^w_I(X))$.  By Theorem
\ref{wst}, this is equivalent to the fact that $T_\lambda(x) \to T(x)$ weak-*
for all $x \in C_I^w(X)$. By the last remark in Appendix \ref{App.B}, this
is the same as saying that $\pi_i(T_\lambda(x)) \to \pi_i(T(x))$ weak-* for all
$x \in C_I^w(X)$ and all $i \in I$. But then for any $x \in X$,
$(\pi_i \circ T_\lambda \circ \epsilon_j)(x) \to (\pi_i \circ T \circ \epsilon_j)(x)$
in the weak-* topology on $X$. From Theorem \ref{wst} again, we conclude that 
$\pi_i \circ T_\lambda \circ \epsilon_j \to \pi_i \circ T  \circ \epsilon_j$ with respect to the 
weak-* topology on $\Al(X)$. Since the $i$-$j$ entry of $\theta(T_\lambda)$ is 
$\pi_i \circ T_\lambda \circ \epsilon_j$, we see that $\theta(T_\lambda) \to \theta(T)$ 
weak-* in $M^w_I(\Al(X))$ (Appendix \ref{App.B}). By the Krein-Smulian theorem
\ref{KS}, $\theta$ is
weak-* continuous.

To see that $\theta$ is an isometry we suppose that
$\Vert T \Vert = 1$ but $\Vert \theta(T) \Vert \leq 
\alpha < 1$.  Then $$\Vert T^\Delta \Vert =
\Vert \theta(T^\Delta) \Vert 
= \Vert \theta(T)^\Delta \Vert \leq
\alpha < 1.$$    On the other hand, for any increasing net of projections
$\{ p_t \}$ in a von Neumann algebra $M$ with supremum $1$, it is
well known that $p_t \to 1$ strongly on the underlying 
Hilbert space.  Hence it is  clear
that the bounded net
$p_t x p_t \to x$ in the weak-* topology
for any $x \in M$.    We conclude that 
$T^\Delta \to T$ weak-* in $\Al(C^w_I(X))$, and it follows from 
the last centered equation that $\Vert T \Vert < 1$.  This
contradiction shows that $\theta$ is an isometry.  

We may now deduce
from the Krein-Smulian theorem that  $\theta$ has 
weak-* closed range.  However, $\theta$ has weak-* dense range,
 by the last fact from Appendix \ref{App.B}.  Thus, $\theta$ is surjective.   
Finally, it is easy to see that $\theta$ restricted to 
$P_\Delta \Al(C^w_I(X)) P_\Delta$ is a *-homomorphism.  By a routine 
weak-* approximation argument, we see that $\theta$ is a *-homomorphism.
Thus, $\theta$ is a $*$-isomorphism of von Neumann algebras.      

We have now proven that $\Al(C^w_I(X)) \cong M^w_I(\Al(X))$.
We next claim that $\Al(R^w_J(X)) \cong \Al(X)$ as von Neumann
algebras.
One way to see this is to first observe that there is a 
weak-* continuous 1-1 $*$-homomorphism 
$\pi : \Al(X) \to \Al(R^w_J(X))$ (we leave this as an exercise).  
Conversely, if $T \in \Al(R^w_J(X))_{sa}$ then 
$T$ commutes with every $S \in \M_r(R^w_J(X))$.  Hence 
$T$ commutes with `projection onto the $i$th entry',
and so $T$ is necessarily of the form
$[x_i] \mapsto [T_i(x_i)]$ for maps
$T_i$ which are easily seen to be in $\Al(X)$ (using Proposition
 \ref{IV.B.2}). 
If $\Delta$ is as above then $T$ restricts to a 
left multiplier of $R_\Delta(X)$, which forces
$T_i = T_j$ for $i, j \in \Delta$.
Thus $T_i = T_j$ for $i, j \in J$, and so $\pi$ above
is a surjection.

Finally, $$\Al(M^w_{I,J}(X)) \cong \Al(R^w_J(C^w_I(X)))
 \cong \Al(C^w_I(X)) \cong M^w_{I}(\Al(X)) , $$  
establishing the final isomorphism in (i) and
also (ii). 

(iii)   This is very similar to the proof above, and we leave the
details to the reader. The main obstacle one encounters is that one 
may no longer appeal to Lemma \ref{IV.B.4b} to see that the 
map  from $P_\Delta \Ml(C^w_I(X)) P_\Delta \to \Ml(C_\Delta(X))$ is isometric.
We present one possible alternative. Let $T \in \Ml(C_\Delta(X))$, with $\|T\|_{\Ml(C_\Delta(X))} \leq 1$.
By Proposition 1.6.4 of \cite{Zarikian}
(or the proof of \cite[Lemma 4.8]{Shilov}), 
there exists a complete isometry $\sigma:X \to B(H)$ and
an $A \in M_\Delta(B(H))$, with $\|A\| \leq 1$, such that
\[
	\sigma_{\Delta,1}(Tx) = A\sigma_{\Delta,1}(x)
\]
for all $x \in C_\Delta(X)$. Here $\sigma_{\Delta,1}:C_\Delta(X) \to C_\Delta(B(H))$ is the natural amplification
of $\sigma$. Now let $\tilde{A} \in M_I^w(B(H))$ be the $I \times I$ matrix whose
$\Delta \times \Delta$ submatrix equals $A$, and all of whose other entries equal $0$. Clearly,
$\tilde{A}\sigma_{I,1}(C_I^w(X)) \subset \sigma_{I,1}(C_I^w(X))$. Thus, there exists a
$\tilde{T} \in \Ml(C_I^w(X))$ such that $\sigma_{I,1}(\tilde{T}x) = \tilde{A}\sigma_{I,1}(x)$
for all $x \in C_I^w(X))$. One has that $\|\tilde{T}\|_{\Ml(C_I^w(X))} \leq \|\tilde{A}\| = \|A\| \leq 1$.
Also, $P_\Delta \tilde{T} P_\Delta|_{C_\Delta(X)} = T$. Since the choice of $T$ was arbitrary, we
conclude that the map $P_\Delta \Ml(C^w_I(X)) P_\Delta \to \Ml(C_\Delta(X))$, which we knew was
a one-to-one contraction, is also surjective, and its inverse is contractive. Hence it is isometric.
As a final note, we observe that once one has the isometry in (iii), the complete isometry follows from 
the relation $M_n(\Ml(X)) \cong
\Ml(M_n(X) \cong \Ml(C_n(X))$. \end{proof} 

{\bf Remark.}   For a general operator space $X$ there
is at least a complete isometry $\Ml(C_I(X)) \to M^w_{I}(\Ml(X))$,
and similarly for $\Al(\cdot)$ in place of $\Ml(\cdot)$.
To see the latter for example, we use the fact that
$C_I(X)^{**} \cong C_I^w(X^{**})$ (see Appendix \ref{App.B}).
Thus by remarks below Proposition \ref{IV.C.0}, 
$\Al(C_I(X))$ may be regarded as
a $C^*$-subalgebra of $\Al(C_I^w(X^{**}))$.  Consider the restriction 
$\rho$ of the 
$*$-isomorphism $\Al(C^w_I(X^{**})) \to M^w_{I}(\Al(X^{**}))$ to 
the subalgebra $\Al(C_I(X))$ of $\Al(C^w_I(X^{**}))$.
This is certainly a complete isometry. 
The second adjoint of the maps  $\epsilon_i$ and $\pi_i$ 
in the proof of Theorem \ref{IV.I.0} are the corresponding 
inclusion and projection
maps between $C_I^w(X^{**})$ and its summands.
Using this and the formula  $T_{ij} = \pi_i \circ T \circ \epsilon_j$
from the proof of Theorem \ref{IV.I.0}, it now easily 
follows that $\rho$ maps into $M^w_{I}(\Al(X))$.  
\\

The following proposition tells us that for a general operator space 
$X$, the diagonal elements of $M_I^w(\Ml(X))$
are contained isometrically in $\Ml(C_I^w(X))$. Although the
result may be deduced from the last Theorem by going to the 
second dual, we give a simpler direct proof.

\begin{proposition} \label{IV.I.1}
Let $X$ be an operator space, $I$ be an index set, and $\{T_i: i \in I\}$ be a uniformly bounded collection of elements 
of $\Ml(X)$.
\begin{enumerate}
\item[(i)] The map
\[
	\bigoplus_{i \in I} T_i:C_I^w(X) \to C_I^w(X):\begin{bmatrix} \vdots\\ x_i\\ \vdots \end{bmatrix}_{i \in I} 
		\mapsto \begin{bmatrix} \vdots\\ T_ix_i\\ \vdots \end{bmatrix}_{i \in I}
\]
is an element of $\Ml(C_I^w(X))$, of multiplier norm $\sup_{i \in I} \|T_i\|_{\Ml(X)}$.
\item[(ii)] $\bigoplus_{i \in I} T_i$ maps $C_I(X)$ into itself, and so restricts to an element of $\Ml(C_I(X))$. 
The restriction process is norm-preserving.
\item[(iii)] If each $T_i \in \Al(X)$, then $\bigoplus_{i \in I} T_i \in \Al(C_I^w(X))$, and 
$\left(\bigoplus_{i \in I} T_i\right)^\star = \bigoplus_{i \in I} T_i^\star$.
\end{enumerate} 
\end{proposition}

\begin{proof}
(i) Without loss of generality, we may assume that $\sup_{i \in I} \|T_i\|_{\Ml(X)} = 1$. It is not hard to see
that $T = \bigoplus_{i \in I} T_i$ is a well-defined linear contraction on $C_I^w(X)$.
Let $\sigma:X \to B(H)$ be a ``Shilov embedding'' 
of $X$ (cf. \cite{BEZ}, \S 1). Then for each $i \in I$, there exists an 
$A_i \in B(H)$ such that $(\sigma,A_i)$ is an implementing pair for $T_i$ and $\|A_i\| = \|T_i\|_{\Ml(X)}$. 
It is easy to see that the $I$-fold `column-amplification' of $\sigma$, namely the
map 
$\sigma^I:C_I^w(X) \to C_I^w(B(H)) \cong B(H,H^{(I)})$ taking $[x_i]$ to $[\sigma(x_i)]$, is a complete isometry. 
 Let $A = \bigoplus_{i \in I} A_i \in M_I^w(B(H)) \cong B(H^{(I)})$.
Then $(\sigma^I,A)$ is an implementing pair for $T$. Thus, $T \in \Ml(C_I^w(X))$ and 
$\|T\|_{\Ml(C_I^w(X))} \leq \|A\| = \sup_{i \in I} \|A_i\| 
= \sup_{i \in I} \|T_i\|_{\Ml(X)} = 1$. Now fix $i \in I$. 
Let $X_i$ denote the closed linear subspace of $C_I^w(X)$ consisting of columns supported in the $i$th position. It is 
clear that $T$ maps $X_i$ into itself. Thus, $T|_{X_i} \in \Ml(X_i)$ and 
$\|T|_{X_i}\|_{\Ml(X_i)} \leq \|T\|_{\Ml(C_I^w(X))}$. But under the completely isometric identification $X_i = X$, 
$T|_{X_i}$ is identified with $T_i$. Thus, $\|T_i\|_{\Ml(X)} \leq \|T\|_{\Ml(C_I^w(X))}$. Since the choice of $i$ was 
arbitrary, $\|T\|_{\Ml(C_I^w(X))} = 1$. 

(ii) This statement is clear. 

(iii) Suppose each $T_i \in \Al(X)_{sa}$.
Then each $A_i \in B(H)_{sa}$, which implies that the matrix
$A$ above is in $M_I^w(B(H))_{sa}$, which in turn implies that 
$T \in \Al(C_I^w(X))_{sa}$. The result now follows.
\end{proof}

\begin{corollary} \label{IV.I.e} Let $\{ J_i : i \in I \}$ be a collection of 
right $M$-ideals (resp. right $M$-summands)
 of an operator space $X$.  Then 
$\{ [x_i] \in C_I(X) : x_i \in J_i \}$ is a right $M$-ideal
of $C_I(X)$ (resp. 
$\{ [x_i] \in C_I^w(X) : x_i \in J_i \}$ is a right $M$-summand of 
$C_I^w(X)$ and $\{ [x_i] \in C_I(X) : x_i \in J_i \}$ is a right $M$-summand 
of 
$C_I(X)$).
\end{corollary}

\begin{proof}  First suppose that
the $J_i$ are right $M$-summands.
If $\{ P_i  : i \in I \}$ is the matching collection
of left $M$-projections, then by the previous proposition
$P = \oplus_{i \in I} P_i $ is in $\Al(C_I^w(X))_{sa}$, and it is
easy to see that $P$ is idempotent, and restricts to an 
idempotent $P' \in \Al(C_I(X))_{sa}$.  Thus $P, P'$ are 
left $M$-projections, and the result is now clear.

Next suppose that
the $J_i$ are right $M$-ideals, so that $J_i^{\perp \perp}$
 is a right $M$-summand
of $X^{**}$ for all $i$. By the first part, $K = \{ [\eta_i] \in C_I^w(X^{**}) : 
\eta_i \in J_i^{\perp \perp} \}$ is a right $M$-summand of 
$C_I^w(X^{**})$.  In particular it is weak-* closed. We will use the 
fact that $C_I(X)^{**} = C^w_I(X^{**})$.
Let $J = \{ [x_i] \in C_I(X) : x_i \in J_i \}$.  Since 
$J \subset K$ we have $J^{\perp \perp} \subset K$.  On the 
other hand, if $z \in K$, then the net of `finitely supported submatrices'
of $z$ is a bounded net converging in the weak-* topology to $z$.
Now any `finitely supported' column in $K$ is weak-* approximable by
a bounded net of `finitely supported' columns in $J$.
Hence $J$ is weak-* dense in $K$, so that $J^{\perp \perp} = K$.
Thus $J$ is a right $M$-ideal of $C_I(X)$.
\end{proof}    

Analogues of the last two results hold
with $C^w_I(X)$ replaced by $M^w_{I,J}(X)$ for 
another index $J$, and with $C_I(X)$ replaced by $K_{I,J}(X)$.
The proofs are almost identical, and are omitted.    Thus
for example we have:

\begin{corollary} \label{IV.I.f} Let $J$ be a
right $M$-ideal (resp. right $M$-summand)
 of an operator space $X$, and let $I$ be an index set.  Then 
$K_I(J)$ 
is a right $M$-ideal
of $K_I(X)$ (resp.
$M^w_I(J)$ is a right $M$-summand of
$M_I^w(X)$, $M_I(J)$
is a right $M$-summand of $M_I(X)$, and 
$K_I(J)$ is a right $M$-summand
of
$K_I(X)$).
\end{corollary}

\subsection{Mutually Orthogonal and Equivalent One-Sided $M$-Projections} \label{IV.J}

\begin{lemma} \label{IV.J.1}
Let $X$ be an operator space, $P_1, P_2, ..., P_n \in \Al(X)$ be mutually orthogonal left $M$-projections, and 
$J_1, J_2, ..., J_n \subset X$ be the corresponding right $M$-summands. Then
\begin{enumerate}
\item[(i)] $P_1 + P_2 + ... + P_n$ is the left $M$-projection onto $J_1 + J_2 + ... + J_n$,
and the complementary left $M$-projection is $P_1^\perp P_2^\perp ... P_n^\perp$. 
In particular, 
$J_1 + J_2 + ... + J_n$ is closed.
\item[(ii)] The map
\[
	X \to C_{n+1}(X):x \mapsto \begin{bmatrix} P_1x\\ P_2x\\ \vdots\\ P_nx\\ P_1^\perp P_2^\perp ... P_n^\perp x
		\end{bmatrix}
\]
is a complete isometry.
\item[(iii)] The map
\[
	J_1 + J_2 + ... + J_n \to C_n(X):x \mapsto \begin{bmatrix} P_1x\\ P_2x\\ \vdots\\ P_nx \end{bmatrix}
\]
is a complete isometry.
\end{enumerate}
\end{lemma}

\begin{proof}
(i) Let $P = P_1 + P_2 + ... + P_n \in \Al(X)$. It is easy to verify that $P$ is a left $M$-projection 
(i.e. $P^2 = P$ and $P^\star = P$). Likewise, it is easy to see that the range of $P$ is $J_1 + J_2 + ... + J_n$. 

(ii) 
We proceed by induction on $n$. The case $n = 1$ follows immediately from Theorem \ref{left_M}. Now suppose that the 
statement is true for $n \geq 1$. That is, suppose that the map
\[
	X \to C_{n+1}(X):x \mapsto \begin{bmatrix} P_1x\\ P_2x\\ \vdots\\ P_nx\\ P_1^\perp P_2^\perp ... P_n^\perp x
		\end{bmatrix}
\]
is a complete isometry. Since, by Theorem \ref{left_M} again, the map
\[
	\mu_{P_{n+1}}^c:X \to C_2(X):x \mapsto \begin{bmatrix} P_{n+1}x\\ P_{n+1}^\perp x \end{bmatrix}
\]
is a complete isometry, we conclude that the map
\[
	X \mapsto C_{2n+2}(X):x \mapsto 
		\begin{bmatrix} P_{n+1}P_1x\\ P_{n+1}^\perp P_1x\\ P_{n+1}P_2x\\ P_{n+1}^\perp P_2x\\
		\vdots\\ P_{n+1}P_nx\\ P_{n+1}^\perp P_nx\\ P_{n+1}P_1^\perp P_2^\perp ... P_n^\perp x\\
		P_{n+1}^\perp P_1^\perp P_2^\perp ... P_n^\perp x \end{bmatrix} = \begin{bmatrix} 0\\ P_1x\\ 0\\ P_2x\\
		\vdots\\ 0\\ P_nx\\ P_{n+1}x\\ P_1^\perp P_2^\perp ... P_n^\perp P_{n+1}^\perp x \end{bmatrix}
\]
is a complete isometry. It follows that the statement is true for $n+1$. 

(iii) This is an immediate consequence of (ii).
\end{proof}

\begin{corollary} \label{IV.J.2}
Let $X$ be an operator space and let
$\{P_i: i \in I\} \subset \Al(X)$ be a family of mutually orthogonal left 
$M$-projections.
\begin{enumerate}
\item[(i)] If $\sum_{i \in I} P_i = \Id_X$ in the point-norm topology, then the map
\[
	X \to C_I(X):x \mapsto \begin{bmatrix} \vdots\\ P_ix\\ \vdots \end{bmatrix}_{i \in I}
\]
is a complete isometry.
\item[(ii)] If $X$ is a dual operator space and $\sum_{i \in I} P_i = \Id_X$ in the weak-* topology, then the map
\[
	X \to C_I^w(X):x \mapsto \begin{bmatrix} \vdots\\ P_ix\\ \vdots \end{bmatrix}_{i \in I}
\]
is a weak-* continuous complete isometry.
\end{enumerate}
\end{corollary}

\begin{proof}
There is a quick proof of (ii) using \ref{IV.I.0}.  Instead 
we give another proof  of (ii) which adapts easily to give (i).
 
(ii) Fix $x \in X$. Then $P_Fx \to x$ weak-*, where for any finite subset 
$F$ of $I$, $P_F = \sum_{i \in F} P_i$. It follows that $\|x\| \leq \sup_F \|P_Fx\|$. Now suppose that
$F = \{i_1, i_2, ..., i_n\} \subset I$. Then by Lemma \ref{IV.J.1},
\[
	\|P_Fx\| = \left\|\begin{bmatrix} P_{i_1}x\\ P_{i_2}x\\ \vdots\\ P_{i_n}x \end{bmatrix}\right\|
		\leq \left\|\begin{bmatrix} \vdots\\ P_ix\\ \vdots \end{bmatrix}_{i \in I}\right\|.
\]
Thus $\|x\| \leq \|[P_i(x)] \|$. On the other hand, for any $i_1, i_2, ..., i_n \in I$, 
one has by Lemma \ref{IV.J.1} that
\[
	\left\|\begin{bmatrix} P_{i_1}x\\ P_{i_2}x\\ \vdots\\ P_{i_n}x \end{bmatrix}\right\| \leq
		\left\|\begin{bmatrix} P_{i_1}x\\ P_{i_2}x\\ \vdots\\ P_{i_n}x\\ P_{i_1}^\perp P_{i_2}^\perp ...
		P_{i_n}^\perp x \end{bmatrix}\right\| = \|x\|,
\]
which implies that $\|[P_i(x)] \|  \leq \|x\|$. So the given map $X \to C_I^w(X)$ is an isometry. 
Since the $n$-fold amplification of this map may be identified with the map
\[
	M_n(X) \to C_I^w(M_n(X)):x \mapsto \begin{bmatrix} \vdots\\ (P_i)_nx\\ \vdots \end{bmatrix}_{i \in I},
\]
and since $\{(P_i)_n: i \in I\}$ is a family of mutually orthogonal left $M$-projections on the dual operator space 
$M_n(X)$ which add to the identity in the weak-* topology, the given map $X \to C_I^w(X)$ is in fact a complete
isometry. That this map is weak-* continuous follows easily from 
the Krein-Smulian Theorem \ref{KS} and the last
remark in Appendix \ref{App.B}, together with the fact that each $P_i$ is
weak-* continuous.  

(i)  This is
similar to (ii).   The only difference is 
that for fixed $x \in X$, we now have that $P_Fx \to x$ in norm. This is a 
favorable difference---not only does it tell 
us that the given map is a 
complete isometry from $X$ to $C_I^w(X)$, but also that the 
range lies in $C_I(X)$. 
\end{proof}

\begin{theorem} \label{IV.J.3}
Let $X$ be an operator space and $I$ be an index set.
\begin{enumerate}
\item[(i)] If $I$ is finite, then $X$ is completely isometric
to $C_I(X_0)$ 
for some operator space $X_0$ if and only if there exists a family 
$\{P_i: i \in I\}$ of mutually orthogonal and equivalent left $M$-projections on $X$ such that 
$\sum_{i \in I} P_i = \Id_X$.  
\item[(ii)] If $I$ is infinite, then $X$ is completely isometric
to $C_I(X_0)$ for some operator space $X_0$ if and only if there exists a family
$\{P_i: i \in I\}$ of mutually orthogonal and equivalent left $M$-projections on $X$ such that 
$\sum_{i \in I} P_i = \Id_X$ point-norm.
\item[(iii)] If $I$ is infinite and $X$ is a dual operator space, then $X$
is completely isometric and weak-* homeomorphic to $C_I^w(X_0)$ for some dual operator space 
$X_0$ if and only if there exists a family $\{P_i: i \in I\}$ of mutually orthogonal and equivalent left 
$M$-projections on $X$ such that $\sum_{i \in I} P_i = \Id_X$ weak-*.
\end{enumerate}

We have $\Al(X_0) \cong P_i \Al(X) P_i$ for each $i \in I$
and $X_0$ as above.
If the equivalent conditions in (i) or (iii) hold,
then we also have $\Al(X) \cong M_I^w(\Al(X_0))$.
\end{theorem}

\begin{proof}
We will only prove (iii) and the statement 
following it. The proof of 
(ii) is similar, and (i) follows immediately from (ii). 

($\Rightarrow$) Suppose
$X = C_I^w(X_0)$ for some dual operator space $X_0$. For each $i \in I$, let $P_i:C_I^w(X_0) \to C_I^w(X_0)$ be the 
projection onto the $i$th coordinate. Then $P_i$ is a left $M$-projection. 
Obviously, the $P_i$'s are mutually orthogonal. We claim that they add to the identity in the weak-* topology. 
To see this, fix $x \in C_I^w(X)$ and $j \in I$. If $F \subset I$ is a finite set containing $j$ and 
$P_F = \sum_{i \in F} P_i$, then $(P_Fx)_j = x_j$. Consequently, $(P_Fx)_j \to x_j$ weak-*. Since the choice of 
$j$ was arbitrary, $P_Fx \to x$ weak-* (see Appendix \ref{App.B}). Since the choice of $x$ was arbitrary, 
the claim is proven. To see that the $P_i$'s are mutually equivalent, fix $i, j \in I$ with $i \neq j$. 
Let $U_{ij}:C_I^w(X_0) \to C_I^w(X_0)$ be the map which swaps the $i$ and $j$ coordinates. Then $U_{ij}$ is
a unitary element of $\Al(C_I^w(X_0))$. Now define $V = U_{ij}P_i$. Then $V \in \Al(C_I^w(X_0))$, $V^\star V = P_i$, 
and $VV^\star = P_j$. Thus, $P_i \sim P_j$. 

($\Leftarrow$) Suppose there exists a family $\{P_i: i \in I\}$ 
of mutually orthogonal and equivalent left $M$-projections on $X$ such that $\sum_{i \in I} P_i = \Id_X$ in the 
weak-* topology. Fix $i \in I$. For each $j \in I$, let $V_j \in \Al(X)$ be such that $V_j^\star V_j = P_j$ and 
$V_j V_j^\star = P_i$. For the sake of definiteness, let $V_i = P_i$. By Proposition \ref{IV.I.1},
$V \equiv \bigoplus_{j \in I} V_j$ is a partial isometry in $\Al(C_I^w(X))$ with range $C_I^w(P_i(X))$. 
Note that $V$ is weak-* continuous by a fact listed after Theorem \ref{wst}. Also, $V$ is completely isometric on 
the range of the map
\[
	\theta:X \to C_I^w(X):x \mapsto \begin{bmatrix} \vdots\\ P_jx\\ \vdots \end{bmatrix}_{j \in I}.
\]
Since, by Corollary \ref{IV.J.2}, $\theta$ itself is completely isometric
and weak-* continuous, the composition $R = V \circ \theta$ is a 
completely isometric weak-* continuous isomorphism of $X$ with $C_I^w(P_i(X))$. 
However $P_i(X)$ is a right $M$-summand of a dual operator space.
Thus by a basic fact mentioned after Theorem \ref{left_M}, we have that 
$P_i(X)$ is weak-* closed, and hence is a dual operator space. By the 
Krein-Smulian theorem \ref{KS}
we see that $R$ is a homeomorphism for the weak-* topologies.

Assuming that the equivalent conditions in (iii)
hold, the fact that $\Al(X) \cong M_I^w(\Al(X_0))$ follows
from Theorem \ref{IV.I.0}. To express $\Al(X_0)$ in 
terms of $\Al(X)$, note that 
by Lemma \ref{IV.B.4b} we need only show that the `inclusion' $P_i \Al(X) P_i \to \Al(X_0)$ is 
surjective.   Since $X \cong C_I^w(X_0)$,  one can check that it suffices
to show that the `inclusion' $Q_i \Al(C_I^w(X_0)) P_i \to \Al(X_0)$ is
surjective, where $Q_i$ is as in the proof of Theorem \ref{IV.I.0}.
However this is a simple consequence of 
Proposition \ref{IV.I.1} (iii) for example.
\end{proof}

{\bf Remark.}  If, in addition to the equivalent conditions
in each part of the last theorem,
we assume that each $P_i$ is equivalent to the identity
map on $X$, then one can deduce that $X$ is `column stable',
in the sense that $X \cong C_I(X)$ 
completely isometrically in (i) or (ii), and 
that $X \cong C_I^w(X)$ completely isometrically via a
weak-* homeomorphism in (iii).  The argument for this is quite 
straightforward: Recall that the $X_0$ mentioned in (i)--(iii)
is $P_i(X)$ for one fixed $i \in I$. Since $P_i$ is
equivalent to $\Id_X$, the desired result follows from facts
towards the end of Section \ref{II.BB}.  

\begin{corollary}  \label{cot1} Let $X$ be a
dual operator space. If $X$
has no complete right $M$-summand which is
 completely 
isometric to $C_2(X_0)$ for an operator space $X_0$,
then all left $M$-projections on $X$ commute.
In particular, if $X$ does not contain $C_2$
completely isometrically, then
all left $M$-projections on $X$ commute. \end{corollary}

\begin{proof} If $\Al(X)$ is noncommutative,
then by the last remark in Section \ref{II.BB},
there exist two equivalent mutually
orthogonal left $M$-projections $P, Q$.  Suppose that
$V V^\star = P, V^\star V = Q$.  Let $R = P + Q$,
another left $M$-projection, and consider the
von Neumann algebra $R \Al(X) R$, which may be viewed
as a von Neumann subalgebra of $\Al(R(X))$,
using Lemma \ref{IV.B.4b}. Then $V R = V Q R = V Q = V$
and $R V = R P V = P V = V$.  Thus $V = R V R \in R \Al(X) R$,
and hence $P \sim Q$ in $\Al(R(X))$.
 Appealing to Theorem \ref{IV.J.3} (i) we have
$R(X) \cong C_2(X_0)$.
\end{proof}                      

\textbf{Remark.} The converse is false. $X = \ell^\infty \htimes C_2$ is a dual operator space which contains
$C_2$ completely isometrically, but $\Al(X) = \ell^\infty$ (see the proof of Corollary \ref{IV.H.4}), 
which is commutative.

\subsection{Multiple $\oplus_{\rM}$-sums} \label{nfold}

Suppose that $\{ P_i : i \in I \}$ is a collection 
of mutually orthogonal left $M$-projections
on an operator space (resp. dual operator space) $X$
which add up to $\Id_X$ in the point-norm (resp. weak-*)
topology. Then we say that $X$ is an
`$\oplus_{\rM}$-sum' of the $J_i$, and write
$X = \oplus_{\rM} J_i$ (resp. $X = \bar{\oplus}_{\rM} J_i$). 
We shall need the `weak-*' case in Section
\ref{Type}. By Corollary \ref{IV.J.2} it follows that if $X = \oplus_{\rM} J_i$ (resp. $X = \bar{\oplus}_{\rM} J_i$),
then the map 
$x \mapsto [P_i(x)]$ from $X$ into $C_I(X)$ (resp. $C_I^w(X)$) is a complete isometry
(resp. weak-* continuous complete isometry). 
We now establish the converse to this. For simplicity we will first 
restrict our attention to the case that $I$ is finite; and then we simply write
$X = J_1 \oplus_{\rM} J_2 \oplus_{\rM} \cdots \oplus_{\rM} J_n$.

\begin{proposition}  \label{glm}  
Let $X$ be an operator space, and suppose that $P_1, P_2, ..., P_n$ are 
linear projections on $X$. Then the following are equivalent:
\begin{itemize}
\item [(i)] $P_1, P_2, ..., P_n$ are left $M$-projections with sum $\Id_X$.
\item [(ii)] $P_iP_j = 0$ for all $i \neq j$, and the map $x \mapsto [P_i(x)]$ from $X$ into $C_n(X)$ 
is a complete isometry.
\item [(iii)] $P_i(X) \rightthreetimes P_j(X)$ if $i \neq j$, and $\sum_i P_i = \Id_X$. 
\item [(iv)] There exists a complete isometry $\sigma:X \to B(H)$ and orthogonal projections
$E_1, E_2, ..., E_n \in B(H)$ such that $(\sigma,E_i)$ is an implementing pair for $P_i$, $1 \leq i \leq n$, and
$E_1 + E_2 + ... + E_n = I$.
\end{itemize}  
\end{proposition}

\begin{proof} (iv) $\Rightarrow$ (iii) The condition implies that the $E_i$ are mutually orthogonal.
Thus $\sigma(P_i(x))^* \sigma(P_j(y)) = 0$ if $i \neq j$. Note also that  
$\sigma(\sum_i P_i(x)) = \sum_i E_i \sigma(x) = \sigma(x)$ for all $x \in X$.

(iii) $\Rightarrow$ (i) If $\langle \cdot , \cdot 
\rangle$ is the `left Shilov inner product' 
(see Subsection \ref{II.BBB}), then 
$\langle P_i(x) , P_j(y) \rangle = 0$ if $i \neq j$.
Hence $$\langle P_i(x) , y \rangle = 
\langle P_i(x) , \sum_j P_j(y) \rangle 
= \langle P_i(x) , P_i(y) \rangle.$$
Likewise, 
$$\langle x , P_i(y) \rangle = \langle P_i(x) , P_i(y) \rangle.$$
Thus $P_i \in \Al(X)_{sa}$. This gives (i).

(i) $\Rightarrow$ (ii) The assumption gives 
 us that $P_iP_j = 0$ if $i \neq j$. Now apply Lemma \ref{IV.J.1}.

(ii) $\Rightarrow$ (iv) Let $\nu^c$ be the map in (ii). Suppose that 
$\rho : X \to B(H)$ is any complete isometry,
and define $\sigma : X \to M_n(B(H))$ to 
be the map taking $x$ to the matrix with first 
column $[\rho(P_i(x))]$ and other columns $0$.
Let $E_i$ be the projection in $M_n(B(H))$
which is $0$ except for an $I$ in the $i$-$i$ entry.
Then $E_i \sigma(x) = \sigma(P_i(x))$. Clearly $\sum_i E_i = I$.
\end{proof}

{\bf Remark.}  A variant of this result holds with
an almost identical proof for a collection $\{ P_i : i \in I \}$ 
of linear projections on $X$.  In this case one must 
interpret the convergence of
the sums in (i) and (iii) in the `point-norm' topology,
and in (iv) one insists that $\sum_i E_i \sigma(x)$ converges
in norm to $\sigma(x)$ for all $x$.

In the `weak-* case' (which we will need in 
Section \ref{Type}) one has the following: 
 
\begin{proposition}  \label{glmd}  Let $X$ be a
dual  operator space, and suppose that $\{ P_i : i \in I \}$ is a collection of
linear projections on $X$. Then the following are equivalent:
\begin{itemize}
\item [(i)] The $P_i$ are left $M$-projections with $\sum_i P_i = \Id_X$ in the weak-* topology of $\Al(X)$.
\item [(ii)] $P_iP_j = 0$ if $i \neq j$, and the map $x \mapsto [P_i(x)]$ from
$X$ into $C^w_I(X)$ is a weak-* continuous complete isometry.
\item [(iii)] $\sum_i P_i = \Id_X$ in the point-weak-* topology, and there is a
weak-* continuous complete isometry $\sigma$ from $X$ into a von Neumann algebra with 
$\sigma(P_i(X))^* \sigma(P_j(X)) = 0$ whenever $i \neq j$.
\item [(iv)] There exists a weak-* continuous complete isometry $\sigma:X \to B(H)$ and orthogonal projections
$E_i \in B(H)$, $i \in I$, such that $(\sigma,E_i)$ is an implementing pair for $P_i$, $i \in I$, and
$\sum_i E_i = I$ (in the weak operator topology).
\end{itemize}
\end{proposition}                              

\begin{proof}
The only change to the proof that (iv) implies (iii) is
an appeal to the Krein-Smulian theorem \ref{KS} to see that
$\sigma$ is a homeomorphism for the weak-* topologies, 
from which it follows that $\sum_i P_i(x) = x$ 
in the weak-* topology for each $x \in X$.
In the proof that (iii) implies (i)
we replace the inner product by $\sigma(\cdot)^* \sigma(\cdot)$.
Also note that if $P_i \in \Al(X)$, then saying that
$\sum_i P_i = \Id_X$ in the point-weak-* topology
is the same as saying that 
$\sum_i P_i = \Id_X$ in the weak-* topology of $\Al(X)$.
In the proof that (ii) implies (iv) we take $\rho$ to also
be weak-* continuous. \end{proof}

The $\oplus_{\rM}$-sum is {\em associative}.
To see this, suppose for specificity 
that $X$ is an operator space with
$X = Y \oplus_{\rM} W \oplus_{\rM} Z$,
for subspaces $Y, W, Z$ of $X$. Then clearly
$Y, W, Z$ are all right $M$-summands of $X$. Also $Y + W, Y + Z,
W + Z$ , are all right $M$-summands of $X$.
Indeed $Y + W$ is the complementary right $M$-summand
of $Z$ in $X$, for example. These statements
are all easiest seen using obvious properties of projections in $C^*$-algebras.

Also note that if $Y, J$ are
right $M$-summands of $X$, with $Y \subset J$, then $Y$ is a right $M$-summand of $J$,
and if $W$ is the complementary right $M$-summand of $Y$
in $J$, then it is easy to see that
$W$ is a right $M$-summand of $X$, and we have $X = Y \oplus_{\rM} W \oplus_{\rM} Z$,
where $Z$ is the  complementary right $M$-summand of
$J$ in $X$. \\
  
{\bf Remark.}  It is easy to check,
by a slight variation of the argument in 
Section \ref{IV.G0}, that (finite or `point-norm')
$\oplus_{rM}$-sums in an ambient operator space
$X$ distribute over a 
Haagerup or spatial tensor product (for example) with a fixed
operator space $V$.  

\subsection{Diagonal Sums} \label{Diag}
 
We write $\oplus_{i \in I}^\infty X_i$ 
for the usual `$\ell^\infty$ direct sum'
of operator spaces. We write $\oplus_i^0 X_i$ for the usual `$c_0$ 
direct sum', that is, the closure in $\oplus_i^\infty X_i$ of the 
finitely supported tuples. Of course if $I$ is finite,
then these direct sums coincide. It is well known that
$(\oplus_i^1 X_i)^* = \oplus_i^\infty X_i^*$ completely
isometrically where $\oplus_i^1 X_i$ is the `$\ell^1$ direct sum'
of operator spaces, and that $(\oplus_i^0 X_i)^* = \oplus_i^1  X_i^*$ 
completely isometrically. Hence $(\oplus_i^0 X_i)^{**} = \oplus_i^\infty X_i^{**}$. 

\begin{proposition} \label{insu}
If $T_i$ are left multipliers (resp. left adjointable multipliers) on operator spaces
$X_i$, with the multiplier norms of the
$T_i$ uniformly bounded above,
then $\oplus_i^\infty T_i$ is a left multiplier
(resp. left adjointable multiplier) on $\oplus_i^\infty X_i$.
Moreover, the multiplier norm of $\oplus_i^\infty T_i$ is
the supremum of the multiplier norms of the $T_i$. In the adjointable case,
$(\oplus_i^\infty T_i)^\star = \oplus_i^\infty T_i^\star$.
\end{proposition}
 
\begin{proof}
This follows from an idea we have used repeatedly. Namely, suppose
that  $(\sigma_i,S_i)$ is an implementing pair for
$T_i$, for each $i \in I$.  Let
$S = \oplus_i^\infty  S_i \in B(\oplus_i H_i)$,
let $\sigma = \oplus_i^\infty \sigma_i :
\oplus_i^\infty X_i \rightarrow B(\oplus_i H_i)$,
and proceed as in the proof of Proposition
\ref{IV.I.1}.
\end{proof}

\begin{corollary} \label{sumell} Let 
$\{ X_i : i \in I \}$ be a collection of operator spaces,
and for each $i \in I$, suppose that $J_i$ 
is a right $M$-ideal (resp. right $M$-summand) of $X_i$.
Then $\{ (x_i) \in \oplus_{i \in I}^0 X_i : x_i \in J_i \}$ 
is a right $M$-ideal of $\oplus_{i \in I}^0 X_i$ (resp.
$\{ (x_i) \in \oplus_{i \in I}^\infty X_i : x_i \in J_i \}$ 
is a right $M$-summand of
$\oplus_{i \in I}^\infty X_i$ and $\{ (x_i) \in \oplus_{i \in I}^0 X_i : x_i \in J_i \}$
is a right $M$-summand of $\oplus_{i \in I}^0 X_i$).
\end{corollary}

\begin{proof}  First suppose that
the $J_i$ are right $M$-summands.
If $\{ P_i  : i \in I \}$ is the matching collection
of left $M$-projections, then by the previous proposition
$P = \oplus_{i \in I} P_i $ is 
in $\Al(\oplus_{i \in I}^\infty X_i)_{sa}$, and it is
easy to see that $P$ is idempotent, and restricts to an
idempotent $P' \in \Al(\oplus_{i \in I}^0 X_i)_{sa}$.  Thus $P, P'$ are
left $M$-projections, and the result is now clear.

Next suppose that
the $J_i$ are right $M$-ideals, so that the $J_i^{\perp \perp}$
are right $M$-summands of $X^{**}$. Then by the first part, 
$K = \{ (\eta_i) \in \oplus_{i \in I}^\infty X_i^{**} : \eta_i \in J_i^{\perp \perp} \}$ 
is a right $M$-summand in $\oplus_{i \in I}^\infty X_i^{**}$. 
The rest of the argument is almost identical
to that of Corollary \ref{IV.I.e}, except that
we use the fact that 
$(\oplus_{i \in I}^0 X_i)^{**} = \oplus_{i \in I}^\infty X_i^{**}$.
\end{proof}
 
\section{One-Sided Type Decompositions and Morita Equivalence for Dual Operator Spaces}

\subsection{One-Sided Type Decompositions} \label{Type}
 
In this section $X$ is a dual operator space, so that 
$\Al(X)$ is a $W^*$-algebra and therefore has a 
type decomposition (see any book on von Neumann algebras).   
We say that a right $M$-summand
$J$ of $X$ is finite, infinite, or properly infinite, according to whether 
the associated
left $M$-projection $P$ has this property as a projection in
$\Al(X)$ (see von Neumann algebra texts for
the definitions). Strictly speaking one should perhaps
say `left finite', `left infinite', or `left properly infinite',
but we suppress the `left' in the following to avoid 
excessive verbiage.
We will use the fact that if $\Id_X$ is not finite,
then it is a unique sum of a central finite projection $P_f$ (possibly $0$)
and a central properly infinite projection $P_{pi}$. Thus if $X$ is not finite,
we can write canonically $X = X_{f} \oplus_{\rM} X_{pi}$, where 
$X_{f} = P_f(X), X_{pi} = P_{pi}(X)$. Similarly, we may decompose the 
identity $\Id_X$ as a sum of central projections corresponding 
to the type decomposition of $\Al(X)$.  We say 
for example, if $P$ is the central projection
in this decomposition such that $\Al(X) P$ is type $\II_1$,
that the summand $P(X)$ is the {\em type $\II_1$ summand} of $X$.  
Thus we can write uniquely $X$ as a $\oplus_{\rM}$-direct sum (see section \ref{nfold})
of a finite type $\I$ summand (or (0)), a purely infinite 
type $\I$ summand (or (0)), a type $\II_1$ summand
(or (0)), a type $\II_\infty$ summand (or (0)), and a type $\III$ summand
(or (0)). The type $\I$ summand can be further 
written as a $\bar{\oplus}_{rM}$-sum of type $\I_n$ 
summands of $X$, for appropriate cardinals $n$.
We say for example that $X$ is of (left) type $\II_1$
if $\Al(X)$ is a type $\II_1$ $W^*$-algebra.

In this section we present a few simple consequences of
this type decomposition.

In regard to (iii) of the next theorem,
we recall that to say `$X$ contains a completely contractively 
complemented copy of $Y$' means that there is a completely contractive
projection of $X$ onto a 
subspace of $X$ which is completely isometrically isomorphic to $Y$.

\begin{theorem} \label{tdec} Let $X$ be a dual operator space.  
\begin{itemize}
\item [(i)] If $X$ is of type $\I_n$ for 
a cardinal $n$, then $X$ is completely isometric
and weak-* homeomorphic to $C^w_n(X_0)$, for a dual operator space $X_0$
which has the property that all its left $M$-projections commute.
\item [(ii)] If $X$ contains a nontrivial type $\I_n$ summand in its 
type decomposition ($n$ possibly infinite), then this summand  is
completely isometric and weak-* homeomorphic to $C^w_n(X_0)$, for some dual operator space $X_0$.
Similarly for a nontrivial type $\II_1$ summand, except that now $n$ can be taken to be any positive
integer.
\item [(iii)]  $X$ is properly infinite if and only if
 $X$ is `column stable' in the sense that $X \cong C_\infty^w(X)$
completely isometrically and weak-* homeomorphically. Indeed if $X$ is not finite, then for 
the properly infinite summand $W$ of $X$ we have $W \cong C_\infty^w(W)$. In this case, $X$
contains a completely contractively complemented infinite dimensional Hilbert column space. 
\item [(iv)] If $X$ is not of type $\I$, then $X$ contains completely contractively 
complemented copies of Hilbert column space of arbitrarily large finite dimension.
Hence, any dual operator space $X$ for which there is a finite $n$ such that $X$ does not contain 
a copy of $C_n$ is of type $\I$.
\item [(v)] If there are no left $M$-projections $Q$ on $X$ with $Q\Al(X)Q$ abelian (i.e.
if $X$ has no type $\I$ summand), then for every finite integer $n$ we have that $X$ is completely 
isometric and weak-* homeomorphic to $C_n(X_0)$, for some dual operator space $X_0$.
\end{itemize}
 
The $X_0$ in (i), (ii), and (v) may be taken to be a
right $M$-summand of $X$.  \end{theorem}

\begin{proof}  
We will treat (ii) in detail; and it may be used as a model for any omitted details 
in other parts.  

(ii) Suppose that $P$ is the central projection in $\Al(X)$ corresponding to the type $\II_1$
summand, and set $W = P(X)$. By Lemma \ref{IV.B.4b}, we may regard $\Al(X)P$ as a unital 
$C^*$-subalgebra of $\Al(W)$. By \cite[Lemma 6.5.6]{KR_II}, for any finite $n$, there are $n$ mutually
orthogonal and Murray-von Neumann equivalent projections $P_i \in \Al(X) P$ adding up to $P$. It is 
clear that the restrictions of these projections give $n$ mutually
orthogonal and Murray-von Neumann equivalent projections in $\Al(W)$ adding up to $\Id_W$.
Thus by Theorem \ref{IV.J.3} and the remark after it, we have that $W$ is completely isometric and
weak-* homeomorphic to $C_n(X_0)$, for a dual operator space $X_0$. In fact, $X_0$ may be taken to be the 
range of $P_i|_W$ (any $i$). But this is the same as the range of $P_i$, which is a right $M$-summand of $X$.

Similarly, let $Y_n$ be the type $\I_n$ piece of $X$, and consider the central
left $M$-projection $Q$ onto $Y_n$. Then $Q$ is the sum of $n$ equivalent 
mutually orthogonal (abelian) projections in $\Al(Y_n)Q$. This implies, as in the last paragraph, 
that $Y_n \cong C^w_n(X_0)$ completely isometrically and weak-* homeomorphically, for a subspace $X_0$ of 
$Y_n$ which is a right $M$-summand of $X$. 

(iii)
 If $X$ is properly infinite then
 we have by \cite[Proposition 2.2.4]{Sakaib} that
$\Id_X = \sum_i P_i$ for a countably infinite family of mutually orthogonal 
left $M$-projections, each equivalent to $\Id_X$. 
Now apply Theorem \ref{IV.J.3} (iii),
and the remark after that result. 
If $X$ is `column stable' then it 
is easy to find a sequence of mutually orthogonal  projections
in $\Al(X)$, each  Murray-von Neumann equivalent to $Id_X$,
and  adding up to $Id_X$.
This is easily seen to imply that $X$ is properly infinite. 

The second part of (iii) follows from the first part; 
we leave the details as an exercise.
Recall that a Hilbert column space is 
{\em injective} (being a `corner' of $B(H)$), so that
it is automatically completely contractively complemented.

(iv) This follows from (ii) and (iii).

(v) Follows from \cite[Lemma 6.5.6]{KR_II}
and Theorem \ref{IV.J.3} (i), as in the first paragraph
of the proof of (ii).           

(i) By (ii) we have $X \cong C^w_n(X_0)$. By the proof of (ii) and by the last assertion of 
Theorem \ref{IV.J.3}, we see that $\Al(X_0) = P_i\Al(X)P_i$ is abelian (recall that the $P_i$ here 
are abelian projections).
\end{proof}

One may easily deduce from the last Theorem
the fact mentioned earlier that if $X$ does not contain 
a copy of $C_2$, then all left $M$-projections on $X$ commute.
On the other hand, the opposite extreme to having all left $M$-projections
commute would be the following:
 
\begin{corollary}  \label{gof} 
Suppose that $X$ is a dual operator space with the property that there are no nontrivial 
left $M$-projections on $X$ which commute with all other left $M$-projections. 
Then one and only one of the following cases occurs:
\begin{itemize} 
\item [(i)] $X \cong C_\infty^w(X)$ completely isometrically and weak-* homeomorphically.
\item [(ii)] $X$ is of finite type, and for all finite $n$ there exists
a dual operator space $X_0$ such that we may write
$X \cong C_{n}(X_0)$ completely isometrically and weak-* homeomorphically. 
\item [(iii)] $X \cong C_n(Z)$ completely isometrically and weak-* homeomorphically
for a finite integer $n$, and a dual operator space $Z$ which has no nontrivial left 
$M$-projections at all.
\end{itemize} 
The $X_0$ and $Z$ here may be taken to be right $M$-summands of $X$.
\end{corollary}

\begin{proof}  The condition asserts that $\Al(X)$ is
a factor. If $\Al(X)$ is infinite, then it is properly infinite and we may appeal to 
Theorem \ref{tdec} (iii). If $\Al(X)$ is type $\II_1$ (and thus
finite), then we appeal to
Theorem \ref{tdec} (ii). If $X$ is finite type $\I$, then $\Al(X) \cong M_n$
for a finite $n$, and $X \cong C_n(X_0)$ as in 
Theorem \ref{IV.J.3}. By the 
last assertion
of Theorem \ref{IV.J.3}  
we have $M_n(\Al(X_0)) \cong M_n$. Thus $\Al(X_0)$ is one-dimensional.

Finally, since $\Al(C_n(X_0)) \cong M_n(\Al(X_0))$, 
we see that cases (ii) and (iii) are mutually 
exclusive.
   So is case (i), that 
being the properly infinite case.
\end{proof}

{\bf Remark:}  We have not explored the connection between our 
type decomposition and the weak-* TRO type decomposition
that can be found in the recent paper \cite{Rtyp}  
for example.  However we imagine that this connection is largely a 
formal one; we certainly do not have access in our general 
framework to the kinds of properties that weak-* TRO's possess.

\begin{example} \label{toexa} 
We give an (dramatic) illustration of the obstacle that arises in using such
type decompositions to study dual operator spaces. Namely, one may have a dual operator space 
$X$ with type decomposition $X = X_1 \oplus_{\rM} X_2$ for example,
where $X_1$ is a type $\I_1$ summand and $X_2$ is a type $\I_m$ summand ($m \neq 1$),
but as an operator space in its own right, $X_1$ is of type $\II_1$ say. This shows that a type $\I$ space can have 
a right $M$-summand which relative to itself may be of quite different type.
Of course this bad behavior cannot occur if the summands are
hereditary in the sense of Definition \ref{IV.B.4c}
(see also Proposition \ref{herind}).

To construct such an example, let ${\mathcal R}$ be the hyperfinite type $\II_1$ factor, 
and let $\tau$ be its trace.  We fix an integer $m > 1$ and identify $M_m \cong B(\mathbb{C}^m)$.   
Let $X$ be the subspace of $B(H \oplus \mathbb{C}^m \oplus \mathbb{C}^m)$
consisting of matrices of the form 
$$ \left[ \begin{array}{cccl} a & 0 & 0 \\ 0 & \tau(a) I_m  & 0\\ 0 & B & 0 \end{array} \right],$$
for $a \in {\mathcal R}, B \in M_m$. We will use basic 
facts about triples and the
triple envelope from \cite{Hamana}. The triple 
subsystem $Z$ generated by $X$ in 
$B(H \oplus \mathbb{C}^m \oplus  \mathbb{C}^m)$ is easily seen to be the subspace consisting of 
matrices of the form
$$ \left[ \begin{array}{cccl} a & 0 & 0 \\ 0 & C  & 0\\ 0 & B & 0 \end{array} \right],$$
for $a \in {\mathcal R}$ and $B, C \in M_m$. It is easy to show that $Z$ has only two 
nontrivial `triple ideals' (using the well known facts that ${\mathcal R}$ and $M_m$ are simple, 
and that triple ideals in $Z$ are in a 1-1 correspondence with two-sided ideals in the $C^*$-algebra $Z^* Z$),
namely the ones isomorphic to ${\mathcal R} \oplus 0$ and $0 \oplus M_{2m,m}$. From this it follows  
that $Z$ is the triple envelope of $X$. Thus the projections in $\Al(X)$ are easily characterized:
they are exactly all the projections of the form 
$\lambda 1_{\mathcal R} \oplus \lambda I_m \oplus P$, 
for $\lambda = 0$ or $1$ in $\mathbb{C}$, and projections $P \in M_m$.
Since $\Al(X)$ is the span of these projections we see that $\Al(X) \cong \mathbb{C} \oplus M_m$.
Thus in the type decomposition of $X$ we only have a type $\I_1$ and a type $\I_m$ summand. 
However if $P$ is the projection corresponding to the  type $\I_1$  summand, then we see that 
$P(X)$ is type $\II_1$ relative to itself.
\end{example} 

Finally we remark on what it means to say that an operator space $X$ is (left) type $\II_1$, for example.
At this point the main significance we see is as follows. Hitherto in operator space theory
the only `column sums' considered are of the form $C_I(X)$ or $C^w_I(X)$.  
We believe that a `left type $\II_1$ operator space', for example, is a generalized kind of
`column sum'. Dual operator spaces of types $\I$--$\III$ simultaneously generalize von 
Neumann algebras of these types, and also `operator space column sums'.  Such generalized 
column and row sums should play a further role in operator space theory.

\subsection{Morita Equivalence} \label{Meq}

In this section we highlight a connection between our work
and the very well understood
Morita equivalence of von Neumann algebras \cite{Rief}.
One of our key tools will be Theorem \ref{IV.I.0}.

\begin{definition} \label{lme}  Two dual operator spaces
$X$ and $Y$ are said to be {\em column Morita equivalent}
if there is a cardinal $I$ such that
$C^w_I(X) \cong C^w_I(Y)$ via a
weak-* continuous complete isometry.  There is a similar
definition for {\em row Morita equivalence}.  We say that
$X$ and $Y$ are  {\em Morita equivalent}    if there 
are cardinals $I, J$ such that $M^w_{I,J}(X) \cong M^w_{I,J}(Y)$ via a
weak-* continuous complete isometry.   
\end{definition}   

Clearly column or row Morita equivalence implies Morita equivalence.

\begin{proposition} \label{rme}  If  dual operator spaces
$X$ and $Y$ are Morita equivalent, then  the von
Neumann algebras $\Al(X)$ and $\Al(Y)$ are 
Morita equivalent in the von Neumann algebra sense 
of \cite{Rief}.  Thus $M^w_I(\Al(X)) \cong M^w_I(\Al(Y))$ as von
Neumann algebras, for some index set $I$.
\end{proposition}

\begin{proof}    It is well known that 
von Neumann algebras $M$ and $N$ are Morita equivalent 
in the sense of \cite{Rief} 
if and only if $M \bar{\otimes}
B(H) \cong N \bar{\otimes}
B(H)$ for some Hilbert space $H$. 
Thus the result follows 
from Theorem \ref{IV.I.0}.  \end{proof}

{\bf Remarks.} 1) 
It follows that von Neumann algebras are Morita equivalent in the
sense of Definition \ref{rme}
if and only if they are Morita equivalent in Rieffel's sense.  
 
2) By taking second duals it follows
that if operator spaces
$X$ and $Y$ are {\em column stably isomorphic}
(that is, if $C_I(X) \cong C_I(Y)$ completely isometrically),
then $X^{**}$ and $Y^{**}$ are column Morita equivalent.  One 
may deduce from this that in this case the right $M$-ideals of 
$X$ corresponding to projections in the center of 
$\Al(X^{**})$
are in a bijective correspondence with the corresponding set of
right $M$-ideals of $Y$.  We omit the details.  \\

Since we are not aware of the following result in the literature,
we include a proof of it supplied by David Sherman.

\begin{lemma} \label{smepr}
Let $M, N$  be von Neumann algebras with orthogonal
projections $p \in M, q \in N$ such that
$M \cong q N q$ and $N \cong p M p$ as von Neumann algebras. 
 Then $M$ and $N$ are Morita equivalent.  That is,
$M \bar{\otimes} B(H) \cong N \bar{\otimes} B(H)$ for 
some Hilbert space $H$.
\end{lemma} 
 
\begin{proof}  
Let $e, f$ be the central supports of $p, q$.  It is well
known that $p M p$ is Morita equivalent to 
$e M$.  Thus $N \bar{\otimes} B(H) \cong (e \otimes 1) (M
\bar{\otimes} B(H))$, for some Hilbert space $H$.    A similar relation holds 
for $M \bar{\otimes} B(K)$. We may assume without loss
of generality, by replacing $H$ and $K$ by Hilbert spaces
of larger dimension,  that $H = K$.  Replacing 
$M$ and $N$ by their tensor products with $B(H)$, 
we see that it suffices to prove our Lemma in the case that
$p, q$ are central projections.   So henceforth
suppose that $p, q$ are central, and we will in fact deduce that
in this case $M \cong N$ as von Neumann algebras.   

Let $f : M \to q N, g : N \to pM$ be the $*$-isomorphisms.
We claim that $f$ maps the center of 
$M$ into the center of $N$.  For 
if $x$ is in the center of
$M$ and $y \in N$ then $$f(x) y = f(x) q y = f(x) f(z) = 
f(x z) = f(z x) = y f(x),$$ for some $z \in M$.   
Similarly $g$ maps between centers.
With some work, one can check that the family $\{ (fg)^k(1_N - q) ,
(fg)^k f(1_M - p) : k =  0, 1, 2, \cdots \}$ consists of mutually 
orthogonal central projections.   Let $d = \sum_k (fg)^k(1_N - q)$ and
$e = \sum_k (fg)^k f(1_M - p)$, then $d e = 0$.
Also, $\{ (fg)^k(1_N) : k =  0, 1, 2, \cdots \}$
is a decreasing sequence of projections; let $r$ be its
infimum
or weak limit.
Since $\sum_{k=0}^n (fg)^k(1_N - q) + \sum_{k=0}^n (fg)^k f(1_M - p) = 1_N - (fg)^{n+1}(1_N)$,
it follows that $d + e + r = 1_N$. Thus we may centrally
 decompose $N$ as a direct sum of three
pieces $Nd, Ne, Nr$.  Similarly $M$ is
a direct sum of three matching pieces $Ma, Mb, Mc$.
It is easy to see that $g$ restricts to an isomorphism
of $Nd$ with $Mb$, and that $f^{-1}$ restricts to an isomorphism 
of $Ne$ with $Ma$. To complete the argument 
we will show that $g$ restricts to an isomorphism of $Nr$ 
with $M c$.    

To that end, note that $x \in Nr$ if and only if 
$x = (fg)^k(1_N) x$ for all $k =  0, 1, 2, \cdots$.
A similar statement holds for $M c$.
For such $x$ we have $g(x) = (g f)^k(g(1_N)) g(x)$.
Hence $(g f)^k(1_M) g(x) = g(x)$ for all $k$, so that
$g(x) \in M c$.    To see that $M c \subset g(Nr)$,
take $y \in M c$.   Thus $(g f)^k(1_M) y = y $ for all $k$.
Let  $x = g^{-1}(y)$, then 
$$g(x) = y = (g f)^k(1_M) y = g((fg)^{k-1}f(1_M) x).$$
Thus $x = (fg)^{k-1}f(1_M) x$ for all $k \in \mathbb{N}$,
so that $x \in Nr$. 
  \end{proof}  

The following was found during a conversation with
David Sherman: 

\begin{proposition} \label{meq}   Suppose that
$X$ is a dual operator space and that 
$Y$ is a hereditary (see
Definition \ref{IV.B.4c}) right $M$-summand of $X$
(or more generally of  $C^w_I(X)$ for a cardinal $I$). 
Then $\Al(X)$ is  Morita equivalent to 
$\Al(Y)$ if either of the following two conditions
hold:
\begin{itemize} 
\item [(i)]  The central support
of the left $M$-projection onto $Y$ is $\Id_X$ (or more generally $\Id_{C_I^w(X)}$).    
\item [(ii)]  There is a hereditary right $M$-summand
of $C^w_J(Y)$ which is completely isometric and 
weak-* homeomorphic to $X$.
 \end{itemize} 
\end{proposition}

\begin{proof} The hypothesis implies that
$\Al(Y) \cong P M P$ as 
von Neumann algebras, where $P$ is the 
left $M$-projection onto $Y$, and $M = \Al(C^w_I(X))$.   
If, in addition, condition (i) holds,
then $P$ is a so-called {\em full projection},
and it is well known that in this case
$P M P$ is Morita equivalent to $M$.  On the 
other hand $\Al(C^w_I(X)) \cong 
M^w_I(\Al(X))$, which is Morita equivalent to 
$\Al(X)$.

Under assumption
(ii) we also have $\Al(X)  \cong Q \Al(C^w_J(Y)) Q$,
and again we recall that $\Al(C^w_J(Y)) \cong
M^w_J(\Al(Y)) \cong \Al(Y) \bar{\otimes} B(H)$,
for some Hilbert space $H$.  
As in the first lines of the proof of Lemma \ref{smepr},
this implies that $\Al(X) \bar{\otimes} B(\tilde{H})$ is 
$*$-isomorphic to  
$e (\Al(Y) \bar{\otimes} B(\tilde{H}))$ for a central 
projection $e \in  \Al(Y) \bar{\otimes} B(\tilde{H})$.
A similar relation holds for 
$\Al(Y) \bar{\otimes} B(K)$ for some Hilbert space $K$.
We may assume that $\tilde{H} = K$, as in the last proof,
and then the result follows from Lemma \ref{smepr}.
\end{proof}

{\bf Remark.} 
Thus if $X$ and $Y$ are two operator spaces with $Y$
completely isometric to a hereditary right $M$-summand of $X$,
and $X$ completely isometric to a hereditary right $M$-summand of $Y$,
then $\Al(X)$ is  Morita equivalent to $\Al(Y)$.   \\
  
We conclude this section by mentioning a few 
properties of dual operator spaces that are invariant under 
Morita equivalence.
For example, it follows from Proposition \ref{rme}
and facts mentioned in the last section 
of \cite{Rief} that if $X$ and $Y$ are 
Morita equivalent, and if $X$  is type $\I$ (resp. type $\II$,
type $\III$), then $Y$ is type $\I$ (resp. type $\II$,
type $\III$).  If $X$ and $Y$ are properly infinite,
then  by Theorem \ref{tdec} (iii) it follows that 
$X \cong Y$ completely isometrically and 
weak-* homeomorphically if and only 
if $X$ and $Y$ are column Morita equivalent.
If $X$ and $Y$ are Morita equivalent, and
if all left $M$-projections on $X$ commute,
then $Y$ is type $\I$ (see \cite[Theorem 8.10]{Rief}. 
  Finally, it follows that
if $X$ and $Y$ are Morita equivalent, then 
the `central right $M$-summands' of $X$ (i.e.\ those
corresponding to projections in the center of $\Al(X)$)
are in bijective correspondence with the 
`central right $M$-summands' of $Y$.  This is because
Morita equivalent rings have isomorphic centers.  
If, further,  $X$ is 
an operator space of the type discussed in Corollary 
\ref{gof}, then so is $Y$.

\section{Central $M$-structure for Operator Spaces}

\subsection{The Operator Space Centralizer Algebra} \label{S71}  

We define the {\em operator space centralizer algebra} $Z(X)$ to
be the subset $\Al(X) \cap \Ar(X)$ inside
$CB(X)$.  
We remark that this space was used in the final
part of \cite{Dual} in the study of certain
operator modules, however we shall not need
anything from there.
We will see  that $Z(X)$ is a slight modification
of the classical
Banach space centralizer algebra
of $X$.
  We recall that 
the latter is a
commutative $C^*$-algebra inside $B(X)$
whose projections are exactly the $M$-projections of $X$
(see \cite{HWW}).  
Note that $\Vert T \Vert = \Vert T \Vert_{cb}$ for $T \in Z(X)$,
since this is true already for $\Al(X)$. 
We will write the Banach space centralizer algebra
as $\Cent(X)$.
Although $Z(X)$ may be developed entirely analogously to
the classical centralizer theory found in
\cite{AlfsenEffros,Behrends,HWW}, we will instead
emphasize the connections to
the one-sided theory, and use these to give
a swift development. The ensuing theory 
is in many ways less interesting than the
one-sided theory presented in previous
sections, precisely because it is so close to the
classical, commutative theory surveyed
in \cite{HWW}.   
We are not saying that $Z(X)$ is unimportant: in fact its
projections are precisely the complete $M$-projections of
Effros and Ruan \cite{ERcmp}.  The corresponding
complete $M$-ideals (or equivalently the 
$M$-ideals whose corresponding
projections in the second dual are in $Z(X^{**})$)
have powerful applications (see e.g\ \cite{ERcmp,MAII,PoonRuan}).
  
From the descriptions of $\Al(X)$ and $\Ar(X)$  from
\cite{Shilov} (or from \cite{BP01}), it  follows easily
that $S T = T S$ if $S \in \Al(X),
T \in \Ar(X)$.  Also, if $T \in Z(X)$ then
the involution $T^\star$ of $T$ in
$\Al(X)$ equals its involution  in $\Ar(X)$.
To see this, write $T^\dagger$ for the latter involution.
Then $T + T^\star$ and $T + T^\dagger$
are Hermitian in the Banach algebra $B(X)$.   Hence their
difference $T^\star - T^\dagger$ is Hermitian.
Similarly, by looking at $i(T - T^\star)$ and
$i(T - T^\dagger)$, we have that $i(T^\star - T^\dagger)$ is Hermitian.
From the proof of Lemma \ref{App.A.8} we deduce that
$T^\star = T^\dagger$.
 
 Thus the subalgebra
$Z(X) = \Al(X) \cap \Ar(X)$ of $\Al(X)$ is a $C^*$-subalgebra
of $\Al(X)$ which is also commutative.  Since the
intersection of weak-* closed subspaces of a dual space
is weak-* closed, it follows from
Theorem \ref{wst} and the abstract characterization
of $W^*$-algebras \cite{Sakai} that
if $X$ is a dual operator space then,            
$Z(X)$ is a commutative $W^*$-algebra.
 
\begin{lemma} \label{chpn}  If $X$ is an operator space
and if $P : X \to X$ is a linear idempotent map, then the following
are equivalent:
\begin{itemize} \item [(i)]  $P \in Z(X)$.
 \item [(ii)] $P$ is a complete $M$-projection
in the sense of Effros and Ruan \cite{ERcmp}.
 \item [(iii)]  $\Vert
[P(x_{ij}) + (I-P)(y_{ij}) ] \Vert \leq \max \{ \Vert
[x_{ij}] \Vert ,
\Vert [y_{ij}] \Vert \}$
for all $n \in
 \mathbb{N}$ and $[x_{ij}], [y_{ij}]$ in $M_n(X)$.
 \item [(iv)]  $ \max \{ \Vert
[P(x_{ij})] \Vert, \Vert[ (I-P)(x_{ij}) ] \Vert \} =  \Vert
[x_{ij}] \Vert$, for all $n \in
 \mathbb{N}$ and $[x_{ij}]$ in $M_n(X)$.
\end{itemize}
\end{lemma}
 
\begin{proof}
(i) $\Leftrightarrow$ (ii)
Clearly the projections in $Z(X)$ are
exactly the maps on $X$ which are
both a left and a right $M$-projection.
By \cite[6.1]{BEZ}
these are exactly the complete $M$-projections.
 
(ii) $\Rightarrow$ (iv)  This is the definition of a
complete $M$-projection.
 
(iii) $\Leftrightarrow$ (iv)  This is an easy
exercise.
 
(iv) $\Rightarrow$ (i)   Given (iv) we will use
Theorem \ref{left_M} (iv) to show that
$P \in \Al(X)$.  Similarly $P \in \Ar(X)$,
and then (i) follows.   Note that  
$$\Vert \nu^c_P(x) \Vert = \max \{  \Vert  P_{2,1}(\nu^c_P(x))  \Vert, \Vert (I-P)_{2,1}(\nu^c_P(x))  \Vert \}
= \max \{ \Vert  Px \Vert , \Vert (I-P)(x) \Vert
= \Vert x \Vert$$  for $x \in X$,
 and similarly for matrices.  \end{proof}
 
{\bf Remark.} Note that condition (i) above implies that
$P$ is a selfadjoint projection in the center of $\Al(X)$.
However the converse is false, as may be seen by looking at
\cite[Example 6.8]{Shilov}. There $X = D_n P$, where
$D_n$ is the diagonal $n \times n$ matrices, and $P$
is a suitable invertible positive matrix. Then
$\Al(X) = D_n$, which has plenty of nonzero central projections;
whereas $Z(X) = \mathbb{C}$, which has only one.

\begin{corollary} \label{Zinc} If $X$ is an operator space,
then $Z(X)$ is a $C^*$-subalgebra of
the Banach space centralizer algebra $\Cent(X)$.
Also, if $X$ is a dual operator space, then 
$\Cent(M^w_\infty(X)) = Z(M^w_\infty(X))$,
and this algebra is isomorphic to
$Z(X)$ as commutative $W^*$-algebras.
\end{corollary}
 
\begin{proof} If $X$ is a dual operator space,
then both $Z(X)$ and $\Cent(X)$ are densely spanned
by their projections.  Since every
complete $M$-projection is an $M$-projection, it
follows from the last Lemma that $Z(X)\subset \Cent(X)$
as unital subalgebras of $B(X)$. In fact, $Z(X)$
must be a $*$-subalgebra of $\Cent(X)$ (as may be seen
by considering
Hermitians for example), and indeed it is a
weak-* closed subalgebra.  To see this recall
that the weak-* topology on $\Cent(Z)$ is the relative
weak-* topology inherited from $B(X)$, and use
Theorem \ref{wst}.
If $R \in \Cent(M^w_\infty(X))$, then we claim that
$R$ is simply the `countably
infinite amplification' $T \otimes I_\infty$
of a map $T \in \Cent(X)$.
This is because of \cite[I.3.15]{HWW}, which says that
any operator in the centralizer commutes with every Hermitian   
operator. One may consider Hermitian
operators on $M^w_\infty(X)$ obtained by left or right multiplying
by a selfadjoint scalar matrix, and in particular Hermitian
permutation and diagonal matrices (such operators are clearly 
Hermitian in $B(M^w_\infty(B(H)))$, if $X \subset B(H)$,
and so their restrictions are Hermitian too). From 
such considerations we see that any
$R \in \Cent(M^w_\infty(X))$ is of the desired form.
Clearly the map $\theta : \Cent(M^w_\infty(X)) \to \Cent(X): T 
\otimes I_\infty \mapsto T$ is a unital 1-1 contractive  
homomorphism, which forces it to be a
$*$-isomorphism onto its range. It is also easy to
see that $\theta$ is weak-* continuous.
Hence by the Krein-Smulian theorem, $\ran(\theta)$ is a $W^*$-subalgebra of $\Cent(X)$.
On the other hand, Lemma \ref{chpn} implies that
$\theta$ takes projections in $\Cent(M^w_\infty(X))$ to
projections in $Z(X)$. Since the span of
the projections is dense, $\theta$ must map into
$Z(X)$.  Conversely, if $T \in Z(X)$, then
by a slight modification of Proposition \ref{IV.I.1}, we see that the
amplification $T \otimes I_\infty$
is both left and right adjointable on
$M^w_\infty(X)$, and hence lies in
$Z(M^w_\infty(X))$. Thus
 $Z(M^w_\infty(X)) = \Cent(M^w_\infty(X))$,
and $\theta$ is the desired
$*$-isomorphism.
 
If $X$ is a general operator space, it suffices to show that any $T$ in the
positive part of $\ball(Z(X))$ is also in $\Cent(X)$.
Now $T^{**}$ is in the
positive part of $\ball(Z(X^{**}))$
by Proposition \ref{IV.C.1} and its matching
`right-handed' version.  Hence $T^{**}$ is in the
positive part of $\ball(\Cent(X^{**}))$.
The classical theory (e.g.\ \cite[Proposition 
I.3.9]{HWW})
implies $T \in \Cent(X)$.
\end{proof}
 
Note that in general $Z(X) \neq \Cent(X)$.
To see this one need
only note that not every $M$-projection on an operator space is
a complete $M$-projection (see \cite{ERcmp}).

On the other hand, the method in the last proof yields 
for a general operator space $X$, a 1-1 *-homomorphism
$Cent(M^w_\infty(X)) \to Cent(X)$.   We shall see later
that the range of this map is $Z(X)$.

\begin{lemma} \label{chpp}  If $X$ is an operator space
and if $T : X \to X$ is linear, then the following
are equivalent:
\begin{itemize} \item [(i)]  $T$ is in the
positive part of $\ball(Z(X))$,
\item [(ii)]  $T^{**}$ is in the
positive part of $\ball(Z(X^{**}))$,
\item [(iii)]  $\Vert
[T(x_{ij}) + (I-T)(y_{ij}) ] \Vert \leq \max \{ \Vert
[x_{ij}] \Vert ,
\Vert [y_{ij}] \Vert \}$
for all $n \in \mathbb{N}$ and $[x_{ij}], [y_{ij}]$ in $M_n(X)$,
\item [(iv)]  Same as (iii) but
for all matrices in $K_\infty(X)$,
\item [(v)]  Same as (iii) but
for all matrices in $M^w_\infty(X)$.
\end{itemize}
\end{lemma}
 
\begin{proof}  First we show that (i) implies
(iii)--(v). It follows by an argument
towards the end of the proof of Corollary \ref{Zinc},
that suitable amplifications of $T$ are in
$Z(M_n(X))$, $Z(K_\infty(X))$, and $Z(M^w_\infty(X))$,
respectively. If,
further, $0 \leq T \leq 1$, then these amplifications
of $T$ are easily seen to also lie between $0$ and $1$, and hence
from Corollary \ref{Zinc} and \cite[Proposition I.3.9]{HWW} we obtain (iii)--(v). 
                                    
(v) $\Rightarrow$ (iv) $\Rightarrow$ (iii)  This is clear.

(ii)  $\Rightarrow$ (i)  The assumption implies that
$T^{**} \in \Al(X^{**})_{sa}$.  Hence $T \in \Al(X)_{sa}$
by Proposition \ref{IV.C.1}.  Similarly
$T \in \Ar(X)$, and the rest is clear.
 
(iii) $\Rightarrow$ (ii)  By taking the
second dual of $T$, we may suppose that
the hypothesis in (iii) holds for $T^{**}$.
Hence it is easy to see that
(v) holds for  $T^{**}$.  Thus by
\cite[Proposition I.3.9]{HWW} we have
$T^{**} \otimes I_\infty$ is in the positive part of the
ball of $\Cent(M^w_\infty(X^{**}))$. By the
previous Lemma it follows that
$T^{**}$ is in the positive part of the
ball of $Z(X^{**})$.
\end{proof}
 
We now establish a host of equivalent characterizations
of $Z(X)$.   We briefly explain the notation
in (vi) below.  Recall from \cite{Shilov}
that  if $({\mathcal T}(X),J)$ is
the triple envelope of $X$, then by one of the definitions
of $\Al(X)$ in that paper there exists  an operator
$\tilde{T}$ in the multiplier algebra
of the `left $C^*$-algebra'
${\mathcal E}(X)$ of
${\mathcal T}(X)$ such that
$\tilde{T} J(x) = J(T(x))$ and
$\tilde{T}^* J(x) = J(T^\star(x))$ for
all $x \in X$.

\begin{theorem}   If $X$ is an operator space
and if $T : X \to X$ is linear, then the following
are equivalent:
\begin{itemize} 
\item [(i)]   $T \in  Z(X)$.
\item [(ii)]   $T^{**}  \in  Z(X^{**})$.
\item [(iii)]   $T$ is a linear combination of four operators
satisfying the condition in
Lemma \ref{chpp} (iii).                    
\item [(iv)]  $T \otimes I_\infty \in \Cent(M^w_\infty(X))$.
Or equivalently, $T$ is in the range of the canonical
$*$-homomorphism $\Cent(M^w_\infty(X)) \rightarrow \Cent(X)$.
\item [(v)]  Same as (iii), with $M_\infty(X)$ replaced
by $K_\infty(X)$.
\item [(vi)]   $T \in \Al(X)$ and the map
$\tilde{T}$ is in the center
of $M({\mathcal E}(X))$ (see the discussion above the
Theorem).
\end{itemize}
\end{theorem}

\begin{proof}   (i) $\Rightarrow$ (ii)  See the last part of
the proof of Corollary \ref{Zinc}.

(ii) $\Rightarrow$ (i)   
Write $T^{**}$ as a norm limit 
of linear combinations of $M$-projections.
Then $T^*$ is a limit of linear combinations of 
$L$-projections.  Using results from 
\cite[p. 26]{Behrends84} we have that 
$T$ is a linear combination of two Hermitian operators.
The result now follows from Proposition
\ref{IV.C.1}, and its matching `right version'.  
    
(i) $\Leftrightarrow$ (iii)  This is evident from
Lemma \ref{chpp}.
 
(iii) $\Rightarrow$ (iv)  Also evident from
Lemma \ref{chpp} and \cite[Proposition I.3.9]{HWW}.
 
(iv)  $\Rightarrow$ (v)  Follows easily for example from
\cite[Proposition I.3.9]{HWW}.
 
(v) $\Rightarrow$ (iv)  We leave this as an
exercise, using  for example \cite[Proposition I.3.9]{HWW}.
 
(iv) $\Rightarrow$ (i)  Assuming (iv), 
write $T \otimes I_\infty$ as a linear combination
of four positive elements $R_i$ in the ball of
$\Cent(M^w_\infty(X))$. By the proof of Corollary
\ref{Zinc}, each $R_i = S_i \otimes I_\infty$
for an $S_i \in \Cent(X)$.  By \cite[Proposition I.3.9]{HWW}
it is easy to see that $S_i$ satisfies the conditions
of  Lemma \ref{chpp}.  Thus $S_i \in Z(X)$, and hence
so is $T$.

(vi)  $\Rightarrow$ (i)
Suppose that $T \in \Al(X)$, and that 
$\tilde{T}$ is the associated
implementing operator in the center of $M({\mathcal E}(X))$. 
 Now ${\mathcal T}(X)$ is
a strong Morita equivalence
${\mathcal E}(X)$-${\mathcal F}(X)$-bimodule,
where ${\mathcal F}(X)$ is the `right $C^*$-algebra of
${\mathcal T}(X)$'.  By a basic fact
about Rieffel's strong Morita equivalence, there
is an operator $S$ in the center of 
$M({\mathcal F}(X))$ such that
$\tilde{T} J(x) = J(T(x)) =
J(x) S$ and $\tilde{T}^* J(x) = J(T^\star(x)) =   
J(x)  S^*$ for all $x \in X$.   From this
it is clear that
$T \in \Ar(X)$.  So $T \in Z(X)$.
 
(i) $\Rightarrow$ (vi)  Suppose that  $T \in Z(X)$.
Then $\tilde{T} J(x) = J(T(x))$ and
$\tilde{T}^* J(x) = J(T^\star(x))$ for
all $x \in X$.  Similarly there exists an operator
$S$ such that $J(x) S = J(T(x))$ and $J(x) S^* = J(T^\star(x))$.
It follows that
$\tilde{T} J(x) J(y)^* =  J(x) S J(y)^* = J(x) J(y)^* \tilde{T}$.  Thus it follows that
$\tilde{T}$ is in the center of $M({\mathcal E}(X))$.
\end{proof}
 
{\bf Remark.} The space $Z(X)$ has other reformulations in terms of
structurally continuous functions, or $M$-boundedness,
just as in the classical case \cite{AlfsenEffros,Behrends,HWW}.
This may be seen from (iv) or (v) 
of the last Theorem; we omit the details.

\begin{corollary}
\label{weht}  Let $X$ be an operator space.
Then $Z(X) \cong \Cent(M^w_\infty(X)) \cong
\Cent(K_\infty(X))$ as commutative $C^*$-algebras.
\end{corollary}
                                       
\begin{proposition}
\label{cont}  \begin{itemize}
\item [(i)] If $X$ is a Banach space, then $\Cent(X)
= Z(\MIN(X))$.
\item [(ii)] If ${\mathcal A}$ is a $C^*$-algebra then
$\Cent({\mathcal A}) = Z({\mathcal A})$, and
this also equals the center of the multiplier
$C^*$-algebra $M({\mathcal A})$.
\item [(iii)] If ${\mathcal A}$ is
an operator algebra with a contractive approximate identity,
then $\Cent({\mathcal A}) = Z({\mathcal A})$, and
this also equals
the diagonal of the center of $M({\mathcal A})$.
\end{itemize}
\end{proposition}
 
\begin{proof}  (i)  This follows from
\cite[Corollary 4.22]{Shilov}.

(ii)  This is well known, and in any case follows from
(iii).
 
(iii) First suppose that ${\mathcal A}$ is a
unital dual operator algebra (a unital weak-* closed
subalgebra of a $W^*$-algebra).
Then  $M({\mathcal A}) = {\mathcal A}$,
and the three algebras which we are trying to show are
equal, are at least commutative $W^*$-algebras,
and may all be viewed as subalgebras of $B({\mathcal A})$.
They are therefore spanned by their projections. If $P$ is an
$M$-projection on ${\mathcal A}$, then from
\cite{ERns} we have $Px = ex$ for a
central projection $e \in {\mathcal A}$, so that
$P$ is a complete $M$-projection. Thus the         
first two algebras are equal. View  $Z({\mathcal A})
\subset \Al({\mathcal A}) \subset CB({\mathcal A})$.
We employ the canonical $*$-homomorphism from
$\Al({\mathcal A})$ into ${\mathcal A}$ given by $T \mapsto T(1)$,
which maps onto the `diagonal' ${\mathcal A} \cap {\mathcal A}^*$.
The restriction $\theta$ of this $*$-homomorphism to
$Z({\mathcal A})$ takes
the $M$-projection $Px = ex$ above to the
central projection $e \in {\mathcal A}$. By linearity and
density considerations, $\theta$ maps into the diagonal of the center of 
${\mathcal A}$, and has dense range. Thus $\theta$ is a *-isomorphism of
$Z(\mathcal{A})$ onto the diagonal of the center of $\mathcal{A}$.
 
Next suppose that ${\mathcal A}$ is not a dual algebra.
If $T \in \Cent({\mathcal A})$, then
$T^{**} \in \Cent({\mathcal A}^{**})$, and hence by the
last paragraph there
exists a unique $\eta$ in the diagonal of the center of
${\mathcal A}^{**}$ with $T^{**}(\nu) = \eta \nu = \nu \eta$ for all
$\nu \in {\mathcal A}^{**}$.
Since $T^{**}(a) = T(a) \in {\mathcal A}$ for $a \in {\mathcal A}$,
it follows that $\eta$ is in $M({\mathcal A})$, and indeed is in the
center of that algebra. The mapping $\pi : \Cent(\mathcal{A}) \to Z(M(\mathcal{A})): T \mapsto \eta$
is a  1-1  isometric homomorphism.
Thus $\pi$ is a $*$-monomorphism and maps into
the diagonal of the center of $M({\mathcal A})$.  Conversely, if
$\eta$ is a selfadjoint (Hermitian) element in
the center of $M({\mathcal A})$, then
we may regard $\eta$ as a
selfadjoint (Hermitian)
in ${\mathcal A}^{**}$ with $\eta a = a \eta$ for
all $a \in
 {\mathcal A}$.  It follows that $\eta$ is in the diagonal
of the center
of ${\mathcal A}^{**}$, so that by the last paragraph there is a
selfadjoint $S \in Z({\mathcal A}^{**})$  with
$S(\nu) = \eta \nu$ for all $\nu \in {\mathcal A}^{**}$.  Since     
$S(a) = \eta a \in {\mathcal A}$  for
all $a \in {\mathcal A}$, it follows
that $S \in Z({\mathcal A})$.
Hence $\pi$ is surjective, and is a $*$-isomorphism, and the
three algebras coincide as desired.
\end{proof}

We will not attempt to systematically transfer 
results from earlier sections to the case of $Z(X)$
and complete $M$-ideals. This is indeed quite 
routine, and almost identical to the classical case.
We do observe that complete $M$-projections
do not exhibit the bad behavior that we saw in 
Example \ref{IV.B.4}.  That is, if $P$ is a projection in $Z(X)$, and
if $Y = P(X)$, then $P \Ml(X) P = \Ml(Y), P \Al(X) P = \Al(Y)$,
and $P Z(X) = Z(Y)$.  
This is because in this case $X \cong P(X) \oplus_\infty (I-P)(X)$,
and one may appeal for example to \cite[Theorem A.13]{Shilov}
and the remark after it.

\subsection{Complete $L$-Projections and the
Cunningham Algebra.}
 
It is clear that $P$ is a complete $L$-projection
if and only if $P^*$ is a complete $M$-projection.
One may define the operator space version
of the classical Cunningham algebra
to be the closed linear span in
$B(X)$ of the complete $L$-projections on $X$.
The theory of the  operator space Cunningham algebra is almost
identical to that of the classical Cunningham algebra
(see e.g.\ \cite{Behrends,HWW}).  Because of
this, and since we have no applications for it,
we leave the details to the interested reader.
Indeed the operator space Cunningham algebra of an operator
space $X$ is a $C^*$-subalgebra of the classical Cunningham algebra
of $X$.  Also, it is 
$*$-isomorphic 
to $Z(X^*)$, and it
may be also be identified with the classical Cunningham algebra
of $S^1_\infty(X)$, where the latter is defined to
be the operator space projective 
tensor product of the predual of $M_\infty$ with $X$.

\section{Future directions}

We begin 
 by mentioning two of the most important 
remaining unsolved problems in the basic theory of 
one-sided $M$-ideals.  
Most importantly, we do not have 
a characterizations of one-sided M-ideals and summands in terms of ball
intersection properties, in the spirit of \cite{AlfsenEffros}
and later researchers (see the works of R. Evans, A. Lima,
and others cited in \cite{HWW}).
Probably related to this is the 
question of proximinality of one-sided M-ideals.  We 
recall that $J \subset X$ is proximinal if 
for all $x \in X$, the distance $d(x,J)$ is achieved.
We suspect that there are examples of  one-sided M-ideals
 which are not proximinal, in which case one could ask for 
sufficient conditions for proximinality.  Such criteria may well
be very useful in operator space theory.
In fact there probably are nonproximinal
left ideals with contractive
right approximate identity in a unital operator algebra
(see Subsection \ref{III.E}).   In this regard we mention 
the following construction suggested by George Willis.  If 
${\mathcal A}$ is a unital Banach algebra, and if $a \in 
{\mathcal A}$, then the closure
of ${\mathcal A}(1-a)$ will be a left ideal with a right
contractive approximate identity if the norms of the elements
$e_n =  1 - \frac{1}{n} \sum_{k=1}^n a^k$ approach 1 (or have
a subnet whose norms approach 1).   By a simple
algebraic computation $e_n \in {\mathcal A}(1-a)$.
It seems feasible that one could get a counterexample
to proximinality of this type, by taking a suitable operator
$a \in B(H)$, and letting  ${\mathcal A}$ be a closed unital
subalgebra of $B(H)$ containing $a$.   
 
In this paper we have only attempted to transfer the `basic theory'
of $M$-ideals to the one-sided situation.  The 
more advanced theory treated in later chapters of \cite{HWW}
remains wide open---for example the theory of spaces that 
are one-sided $M$-ideals in their bidual.   We will however
wait until applications arise that call for the 
development of such further theory.
 
The first author has recently introduced
one-sided multipliers between two different spaces
\cite{BSur}.
   We would guess that
the proofs of several of our results for multipliers in this paper
(for example the polar decomposition theory)
should adapt to give analogous results and applications
for the new class of
multipliers.  
In a similar spirit, there is a recent notion
of `quasimultipliers' introduced by Kaneda
and Paulsen \cite{Kthes}.  Quasimultipliers are intimately
related to the one-sided multipliers,
and may be developed in a somewhat analogous fashion.
Kaneda has suggested a `quasi'- or `inner' version
of one-sided $M$-ideals which would be interesting
to pursue.   We would guess that
the proofs of several of our results for multipliers and ideals in this paper
should adapt to give analogous results
for quasimultipliers.   However since the quasimultiplier space has 
so little structure in general, others of our results may fail, 
or be quite difficult, to transfer.                             

\newpage

\appendix

\section{Some Results from Banach Space Theory} \label{App.A}

In this section we collect some facts that we shall need from Banach space theory. We provide references 
when possible, and proofs otherwise.

\subsection{Dual spaces} 

\begin{theorem} \label{KS} {\rm (Krein-Smulian)}  
\begin{itemize}
\item [(i)]  Let $X$ be a dual Banach space, and $Y$ be a linear subspace of $X$.  Then
$Y$ is weak-* closed in $X$ if and only if  $\ball(Y)$ is closed in the  weak-* topology 
on $X$.  In this case, $Y$ is also a dual Banach space, with predual $X_*/ {}^\perp Y$, and the
inclusion of $Y$ into $X$ is weak-* continuous.
\item [(ii)]  A bounded linear map $T$ between dual Banach spaces is
weak-* continuous if and only if whenever $x_i \rightarrow x$
is a bounded net converging
weak-* in the domain space, it follows that
$T(x_i) \rightarrow T(x)$ weak-*.
\item [(iii)]   Let $X$ and $Y$ be dual Banach spaces, and
$T : X \rightarrow Y$ be a weak-* continuous linear isometry.
Then $T$ has weak-* closed range $V$ say, and $T$ is a
weak-* homeomorphism of $X$ onto $V$.
\end{itemize}
\end{theorem}                           

Variations of the Krein-Smulian theorem can be found in any functional analysis book. 
The above formulation and its proof appears in \cite{ERbimod,Dual}.

\begin{lemma} \label{App.A.1}
Let $X$ be a Banach space and $J$ be a closed linear subspace of $X$. If $J^{\perp\perp}$ is contractively
complemented in $X^{**}$ and $J$ is a dual Banach space, then $J$ is contractively complemented in $X$.
\end{lemma}

\begin{proof}
Let $Y$ be a Banach space such that there exists an isometric isomorphism $\Phi:Y^* \to J$. Let 
$P:X^{**} \to J^{\perp\perp}$ be a contractive surjection,
 and let $\Psi:J^{\perp\perp} \to J^{**}$ be the
canonical isometric isomorphism. For any Banach space $Z$, let $\iota_Z:Z \to Z^{**}$ be the canonical
isometric inclusion. Then the composition
\[
\xymatrix{
	X \ar[r]^{\iota_X} & X^{**} \ar[r]^P & J^{\perp\perp} \ar[r]^\Psi & J^{**} \ar[r]^{(\Phi^{-1})^{**}}
		& Y^{***} \ar[r]^{\iota_Y^*} & Y^* \ar[r]^\Phi & J
}
\]
is a contractive surjection.
\end{proof}

\begin{lemma} \label{App.A.2}
Let $X$ be a Banach space and let 
$J_a$, $a \in A$, be closed linear subspaces of $X$. Then
\begin{enumerate}
\item[(i)] $\left(\bigcup_{a \in A} J_a\right)^\perp = \bigcap_{a \in A} J_a^\perp$.
\item[(ii)] $\left(\bigcap_{a \in A} J_a\right)^\perp = 
\overline{\spn}^{\wks}\left\{\bigcup_{a \in A} J_a^\perp\right\}$.
\end{enumerate}
\end{lemma}

\begin{proof}
(i) This is clear. In fact, the $J_a$'s can be arbitrary sets in this case. 

(ii) Likewise, the inclusion 
$\overline{\spn}^{\wks}\left\{\bigcup_{a \in A} J_a^\perp\right\} \subset 
\left(\bigcap_{a \in A} J_a\right)^\perp$ is straightforward, even if the $J_a$'s are arbitrary sets. Now 
suppose that
 $f \in \left(\bigcap_{a \in A} J_a\right)^\perp$ but 
$f \notin \overline{\spn}^{\wks}\left\{\bigcup_{a \in A} J_a^\perp\right\}$. Then there exists a weak-* 
continuous linear functional $\Phi$ on $X^*$ such that $\Phi(f) \neq 0$ but $\Phi(g) = 0$ for all 
$g \in \overline{\spn}^{\wks}\left\{\bigcup_{a \in A} J_a^\perp\right\}$. Since $(X^*,\wks)^* = X$, there 
exists an $x \in X$ such that $f(x) \neq 0$ but $g(x) = 0$ for all 
$g \in \overline{\spn}^{\wks}\left\{\bigcup_{a \in A} J_a^\perp\right\}$. It follows that 
$x \notin \bigcap_{a \in A} J_a$, so that there exists a $b \in A$ with $x \notin J_b$. Because 
${}^\perp(J_b^\perp) = \overline{\spn}\{J_b\} = J_b$, there exists a 
$g \in J_b^\perp$ such that $g(x) \neq 0$, a contradiction.
\end{proof}

\begin{lemma} \label{App.A.3}
Let $X$ be a Banach space and let 
$J_a$, $a \in A$, be weak-* closed linear subspaces of $X^*$. Then
\begin{enumerate}
\item[(i)] ${}^\perp\left(\bigcup_{a \in A} J_a\right) = \bigcap_{a \in A} {}^\perp J_a$.
\item[(ii)] ${}^\perp\left(\bigcap_{a \in A} J_a\right) = 
\overline{\spn}\left\{\bigcup_{a \in A} {}^\perp J_a\right\}$
\end{enumerate}
\end{lemma}

\begin{proof}
(i) This is clear for arbitrary sets $J_a$. 

(ii) The inclusion 
$\overline{\spn}\left\{\bigcup_{a \in A} {}^\perp J_a\right\} 
\subset {}^\perp\left(\bigcap_{a \in A} J_a\right)$ is easy, even for arbitrary sets $J_a$. Now suppose that
$x \in {}^\perp\left(\bigcap_{a \in A} J_a\right)$ but 
$x \notin \overline{\spn}\left\{\bigcup_{a \in A} {}^\perp J_a\right\}$. Then there exists an $f \in X^*$ 
such that $f(x) \neq 0$ but $f(y) = 0$ for all 
$y \in \overline{\spn}\left\{\bigcup_{a \in A} {}^\perp J_a\right\}$. It follows that 
$f \notin \bigcap_{a \in A} J_a$, so that $f \notin J_b$ for some $b \in A$. Since 
$({}^\perp J_b)^\perp = \overline{\spn}^{\wks}\{J_b\} = J_b$, there exists a 
$y \in {}^\perp J_b$ such that $f(y) \neq 0$, a contradiction.
\end{proof}

\begin{lemma}[\cite{HWW}, Lemma I.1.14] \label{App.A.4}
Let $X$ be a Banach space, and $J$ and $K$ be closed linear subspaces of $X$. Then the following are equivalent:
\begin{enumerate}
\item[(i)] $J + K$ is norm closed in $X$.
\item[(ii)] $J^\perp + K^\perp$ is norm closed in $X^*$.
\item[(iii)] $J^\perp + K^\perp$ is weak-* closed in $X^*$.
\end{enumerate}
\end{lemma}

\begin{lemma}[\cite{Fabian}, Exercise 5.13] \label{App.A.5}
Let $X$ be a Banach space, and $J$ and $K$ be closed linear subspaces of $X$. Then the following are equivalent:
\begin{enumerate}
\item[(i)] $X = J \oplus K$.
\item[(ii)] $X^* = J^\perp \oplus K^\perp$.
\item[(iii)] $X^{**} = J^{\perp\perp} \oplus K^{\perp\perp}$.
\end{enumerate}
\end{lemma}

\begin{lemma}[\cite{Fabian}, Exercise 5.18] \label{App.A.6}
Let $X$ be a Banach space and $P:X^* \to X^*$ be a projection. If the range and kernel of $P$ are weak-*
closed, then $P$ is weak-* continuous.
\end{lemma}

\subsection{Miscellaneous facts}

\begin{lemma}[\cite{Fabian}, Exercise 5.27] \label{App.A.7}
Let $X$ be a Banach space, and $Y$ and $Z$ be closed linear subspaces of $X$. If $Y$ is finite-dimensional, then
$Y + Z$ is closed.
\end{lemma}

\begin{lemma} \label{App.A.8}
Let $\mathcal{B}$ be a unital Banach algebra, and $a, b, c, d \in \Her(\mathcal{B})$. If
$a + ib = c + id$, then $a = c$ and $b = d$.
\end{lemma}

\begin{proof}
It suffices to show that if $a, b \in \Her(\mathcal{B})$ and $a = ib$, then $a = b = 0$. 
For such $a, b$, and for any state $f$ on $\mathcal{B}$, we have
$f(a) = f(b) = 0$.  Thus the {\em numerical radii} are zero,
and it is well known that this implies that $a = b = 0$
(see e.g.\ \cite{BonsallDuncan}).
\end{proof}

\begin{lemma} \label{App.A.9}
Let $X$ be a Banach space and $P, Q:X \to X$ be bounded projections. If $\|PQ\| < 1$, then
$\ran(P) + \ran(Q)$ is closed.
\end{lemma}

\begin{proof}
We give Akemann's argument.
Assume that $\ran(P) + \ran(Q)$ is not closed. Consider the bounded linear map
\[
	T:\ran(P) \oplus_\infty \ran(Q) \to X:(x,y) \mapsto x + y.
\]
Since $\ran(T) = \ran(P) + \ran(Q)$, $T$ cannot be bounded below. Therefore, there exist sequences
$(x_n)$ in $\ran(P)$ and $(y_n)$ in $\ran(Q)$ such that $\|x_n\| = \|y_n\| = 1$ and
$\|x_n - y_n\| \leq 2^{-n}$ for all $n \in \mathbb{N}$. But then
\[
	1 = \|x_n\| = \|Px_n\| \leq \|P(x_n - y_n)\| + \|Py_n\| \leq \|P\|\|x_n - y_n\| + \|PQy_n\|
		\leq 2^{-n}\|P\| + \|PQ\|
\]
for all $n \in \mathbb{N}$. Choosing $n$ sufficiently large results in a contradiction. 
\end{proof}

\section{Infinite Matrices over an Operator Space} \label{App.B}

In this section we briefly discuss the theory of infinite matrices over an operator space. For a more
detailed discussion, see \cite{ERbook}, \S 10.1, and the references contained therein.\\

Let $X$ be an operator space and $I$ be a (typically infinite) index set. By $M_I^w(X)$ we mean the 
linear space of all $I \times I$ matrices
\[
	x = \begin{bmatrix} & \vdots & \\ ... & x_{ij} & ... \\ & \vdots & \end{bmatrix}_{i, j \in I}
\]
over $X$ such that
\[
	\|x\| \equiv \sup\left\{\left\|\begin{bmatrix}  & \vdots & \\ ... & x_{ij} & ... \\ & \vdots & 
		\end{bmatrix}_{i, j \in F}\right\|: F \subset I \text{ is finite}\right\} < \infty.
\]
It is not hard to verify that $M_I^w(X)$ is a Banach space. In fact, it is an operator space with respect 
to the matrix norms determined by the linear isomorphisms
\[
	M_n(M_I^w(X)) = M_I^w(M_n(X)).
\]

Similar definitions and observations pertain to `rectangular'
matrices $M_{I,J}^w(X)$.  We write $C_I^w(X)$ for 
$M_{I,1}^w(X)$, and $R_I^w(X)$ for $M_{1,I}^w(X)$. 
In fact $M_{I,J}^w(X) \cong R_J^w(C_I^w(X)) \cong C_I^w(R_J^w(X))$.  If $I = \aleph_0$
we simply write $C_\infty^w(X)$ and $R_\infty^w(X)$ for $C_I^w(X)$ and $R_I^w(X)$.

Some simple examples and observations are in order:
\begin{itemize}
\item $M^w_I(\mathbb{C}) = B(\ell^2(I))$. This space is typically abbreviated $M_I$.
\item More generally, if $\mathcal{R}$ is a $W^*$-algebra, then 
$M_I^w(\mathcal{R}) = M_I \stimes \mathcal{R}$. In particular, $M_I^w(\mathcal{R})$ is again a 
$W^*$-algebra.
\item If $\mathcal{A}$ is $C^*$-algebra, 
then $M_I^w(\mathcal{A})$ need not be a $C^*$-algebra. In fact, 
it need not even be an algebra.
\end{itemize}

Contained inside $M_I^w(X)$ are a number of distinguished operator subspaces:
\begin{itemize}
\item $M_I(X) = M_I \itimes X$, the norm closed linear span of $\{A \otimes x: A \in M_I, x \in X\}$ (the
`elementary tensors').
\item $K_I(X)$, the norm closure of the `finite rank' elements of $M_I^w(X)$.
\item $C_I^w(X)$ and $R_I^w(X)$, which may be 
identified with a fixed column and row of $M_I^w(X)$, respectively.
\item $C_I(X)$ and $R_I(X)$, which may be
identified with a fixed 
 column and row of $M_I(X)$, respectively. Equivalently, 
with a fixed  
column and row of $K_I(X)$, respectively.
\item $R_I(C_I^w(X))$, $C_I^w(R_I(X))$, $C_I(R_I^w(X))$, and $R_I^w(C_I(X))$.
\end{itemize}
The reader will have no trouble verifying the following diagram of inclusions:
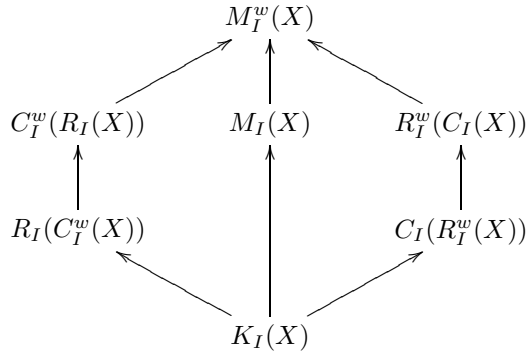
\[
\xymatrix{
	& M_I^w(X) &\\
	C_I^w(R_I(X))\ar[ur] & M_I(X)\ar[u] & R_I^w(C_I(X))\ar[ul]\\
	R_I(C_I^w(X))\ar[u] & & C_I(R_I^w(X))\ar[u]\\
	& K_I(X)\ar[ul]\ar[uu]\ar[ur] & \\
}
\]
And of course $C_I(X) \subset C_I^w(X)$ and $R_I(X)
\subset R_I^w(X)$. We write $\infty$ in place of the subscript $I$
if $I$ is countably infinite.

Some more examples and observations:
\begin{itemize}
\item  Taking $X = \mathbb{C}$ 
in the `column and row matrix' definitions
above, we have $C_I^w = C_I$ and $R_I^w = R_I$.  Indeed
if $H$ is a Hilbert space with orthonormal basis 
$\{e_i: i \in I\}$, then $H_c = C_I$ completely isometrically and
$H_r = R_I$ completely isometrically.
\item If $\mathcal{A}$ is a unital operator 
algebra, then $C_I(R_I^w(\mathcal{A}))$ and $R_I(C_I^w(\mathcal{A}))$ are operator algebras
with contractive approximate identities, and $R_I^w(C_I(\mathcal{A}))$ and $C_I^w(R_I(\mathcal{A}))$ are unital operator algebras.
\item  The second dual of $K_I(X)$ is $M^w_I(X^{**})$ 
(see \cite[Theorem 10.1.4]{ERbook}).  The 
second dual of $C_I(X)$ is $C^w_I(X^{**})$, and a similar assertion 
holds for rows.  \end{itemize}

As a final remark, we note that if $X$ is an operator space and $I$ is an index set, then $M_I^w(X^*)$ is a 
dual operator space. Indeed, $M_I^w(X^*) = CB(X,M_I)$ completely isometrically. A bounded net $(f_\alpha)$ in 
$CB(X,M_I)$ converges to $f \in CB(X,M_I)$ in the weak-* topology if and only if $f_\alpha(x) \to f(x)$ weak-* 
for all $x \in X$, which is the same as saying that $\tr(f_\alpha(x)A) \to \tr(f(x)A)$ for all $x \in X$ and 
all trace-class operators $A \in M_I$. This, in turn, is the same as saying that 
$f_\alpha(x)_{i,j} \to f(x)_{i,j}$ for all $x \in X$ and all $i, j \in I$. Therefore, weak-* convergence of 
bounded nets in $M_I^w(X^*)$ is the same as entry-wise weak-* convergence.


\newpage

\begin{thebibliography}{99}

\bibitem[Ake1]{Akemann69}
C. Akemann, \emph{The general Stone-Weierstrass problem}, J. Funct. Anal. \textbf{4} (1969),
277--294.

\bibitem[Ake2]{Akemann70}
\bysame, \emph{Left ideal structure of $C^*$-algebras}, J. Funct. Anal. \textbf{6}
(1970), 305--317.

\bibitem[AE]{AlfsenEffros}  
E. M. Alfsen and E. G. Effros, {\em Structure in real Banach spaces. I, II}, Ann. of Math. {\bf 96} (1972), 98--173.

\bibitem[AR]{MAII}  
A. Arias and H. P. Rosenthal, \emph{$M$-complete approximate identities in operator spaces},
Studia Math. {\bf 141} (2000), 143--200.

\bibitem[Arv1]{Arv1}  
W. B. Arveson, \emph{Subalgebras of $C^*$-algebras}, Acta Math. {\bf 123} (1969), 141--224.

\bibitem[Arv2]{Arv2}  
\bysame, \emph{Subalgebras of $C^*$-algebras. II}, Acta Math. \textbf{128} (1972), 271--308.

\bibitem[Beh1]{Behrends}
E. Behrends, \emph{$M$-structure and the Banach-Stone theorem}, Lecture Notes in Mathematics 736,
Springer-Verlag, Berlin-Heidelberg-New York, 1979.

\bibitem[Beh2]{Behrends84}
\bysame, \emph{Normal operators and multipliers on complex Banach spaces and a symmetry property
of $L^1$-predual spaces}, Israel J. Math. \textbf{47} (1984), 23--28.

\bibitem[Ble1]{Ble_cb}
D. Blecher, \emph{A completely bounded characterization of operator algebras}, Math. Ann. \textbf{303} (1995), 
227--239.

\bibitem[Ble2]{Bgeo}
\bysame, \emph{Geometry of the tensor product of $C^*$-algebras}, Math. Proc. Camb. Phil. Soc. 
\textbf{104} (1988), 119--127.

\bibitem[Ble3]{Dual}
\bysame, \emph{Multipliers and dual operator algebras}, J. Funct. Anal. \textbf{183} (2001), 
498--525.


\bibitem[Ble4]{BSur}  \bysame, \emph{Multipliers, $C^*$-modules,
and algebraic structure in spaces of Hilbert space operators},
 Preprint (2003). 

\bibitem[Ble5]{OpAlg} 
\bysame, \emph{One-sided 
ideals and approximate identities in operator algebras}, 
J. Austral. Math. Soc., accepted.

\bibitem[Ble6]{Shilov}
\bysame, \emph{The Shilov boundary of an operator space and the characterization theorems}, 
J. Funct. Anal. \textbf{182} (2001), 280--343.


\bibitem[BEZ]{BEZ}
D. Blecher, E. Effros, and V. Zarikian, \emph{One-sided $M$-ideals and multipliers in operator 
spaces. I}, Pacific J. Math. \textbf{206} (2002), 287--319.

\bibitem[BMP]{BMP}  D. Blecher, P. Muhly, and V. Paulsen, 
{\em  Categories of
operator modules - Morita
equivalence and projective modules}
Memoirs of the A.M.S., Vol. 143, number 681,
 January 2000.
 
\bibitem[BP]{BP01}
D. Blecher and V. Paulsen, \emph{Multipliers of operator spaces and the injective envelope}, 
Pacific J. Math. \textbf{200} (2001), 1--17.

\bibitem[BSZ]{BSZ}
D. Blecher, R. Smith, and V. Zarikian, \emph{One-sided projections on $C^*$-algebras}, 
J. Operator Theory, accepted.

\bibitem[Boh]{Bohnenblust}
F. Bohnenblust, \emph{A characterization of complex Hilbert spaces}, Portugaliae Math. \textbf{3} 
(1942), 103--109.

\bibitem[BD]{BonsallDuncan}  
F. F. Bonsall and J. Duncan, \emph{Complete normed algebras}, Springer-Verlag, New York-Heidelberg, 
1973.

\bibitem[Eff1]{Effros} 
E. G. Effros, \emph{Order ideals in a $C^*$-algebra and its dual}, Duke Math. J. \textbf{30} (1963), 391--412.



\bibitem[ER1]{ERcmp}   
E. G. Effros and Z-J. Ruan, {\em Mapping spaces and liftings for operator spaces},  Proc. London Math. Soc.
{\bf 69} (1994), 171--197.

\bibitem[ER2]{ERns}  
\bysame, {\em On non-self-adjoint operator algebras},  Proc. Amer. Math. Soc. 
{\bf 110} (1990), 915--922.

\bibitem[ER3]{ERoc}
\bysame, {\em Operator
convolution algebras: an approach to
quantum groups}, Unpublished (1991).

\bibitem[ER4]{ERbook}
\bysame, \emph{Operator spaces}, Oxford University Press, Oxford, 2000.

\bibitem[ER5]{ERbimod} \bysame, 
 {\em Representations of
operator bimodules and their applications, }J. Operator Theory {\bf 19 }%
(1988), 137-157.         

\bibitem[Fab]{Fabian}
Fabian et al., \emph{Functional analysis and infinite-dimensional geometry}, Springer-Verlag, New York,
2001.

\bibitem[GKS]{GKS}
G. Godefroy, N. Kalton, and P. Saphar, \emph{Unconditional ideals in Banach spaces}, Studia Math. \textbf{104} (1993), 
13--59.

\bibitem[Hal]{Halmos}
P. Halmos, \emph{Introduction to Hilbert space}, Chelsea Publishing Company, New York, 1957.
 
\bibitem[Ham]{Hamana}  
M. Hamana, {\em Triple  envelopes and Silov boundaries of operator spaces,} Math. J. Toyama University
{\bf 22} (1999), 77--93.
                  
\bibitem[HWW]{HWW}
P. Harmand, D. Werner, and W. Werner, \emph{$M$-ideals in Banach spaces and Banach algebras}, 
Lecture Notes in Mathematics 1547, Springer-Verlag, Berlin-Heidelberg-New York, 1993.

\bibitem[KR]{KR_II}
R. V. Kadison and J. R. Ringrose, \emph{Fundamentals of the theory of operator algebras. Vol. II: Advanced theory}, 
Amer. Math. Soc., Providence, 1997.

\bibitem[Kak]{Kakutani}
S. Kakutani, \emph{Some characterizations of Euclidean space}, Japan. J. Math. \textbf{16} (1939), 
93--97.

\bibitem[Kan]{Kthes}  
M. Kaneda, {\em Multipliers and algebrizations
of  operator spaces,}
Ph. D. Thesis University of Houston (2003).


\bibitem[Kir]{Kirchberg}  
E. Kirchberg,  {\em On restricted perturbations in inverse images and a description of normalizer algebras in 
$C^*$-algebras}, J. Funct. Anal. {\bf 129} (1995), 35--63.    
 
\bibitem[Lan]{La}     
E. C. Lance, {\em Hilbert C$^*$-modules -- A toolkit for operator algebraists,} London Math. Soc. Lecture Notes,
Cambridge University Press, Cambridge, 1995.

\bibitem[Pau]{Paulsen}  \emph{Completely bounded maps and operator algebras},
Cambridge University Press (2002).

\bibitem[Pis1]{Pisier}
G. Pisier, \emph{Introduction to operator space theory}, 
London Math. Soc. Lecture Note Series 294,
Cambridge University Press (2003).
                              

\bibitem[Pis2]{OH}  
\bysame, \emph{The operator Hilbert space $OH$, complex interpolation and tensor norms}, Mem. 
Amer. Math. Soc. \textbf{122} (1996), no. 585, 1--103.

\bibitem[PR]{PoonRuan}
Y-T. Poon and Z-J. Ruan, \emph{Operator algebras with contractive approximate identities}, 
Canadian J. Math. \textbf{46} (1994), 397--414.

\bibitem[Pow]{Powers}
R. T. Powers, \emph{Representations of uniformly hyperfinite algebras and their associated von 
Neumann rings}, Ann. of Math. \textbf{86} (1967), 138--171.

\bibitem[Pro]{Prosser}  
R. T. Prosser, {\em On the ideal structure of operator algebras},
Mem. Amer. Math. Soc. {\bf 45} (1963).    

\bibitem[Rie]{Rief}   
M. A. Rieffel, {\em Morita equivalence for }$C^{*}$-{\em algebras and} $W^{*}$-{\em algebras,} 
J. Pure Appl. Algebra {\bf 5 }(1974), 51--96.              

 
\bibitem[Rua]{Rtyp}  Z. J. Ruan, {\em   Type decomposition and
rectangular AFD property for $W^*$-TRO's,} To appear
Canadian J. Math..
   
\bibitem[RN]{RN}  B. Russo and M. Neal, {\em State spaces
of $JB^*$-triples,} Preprint (2003).  

\bibitem[Sak1]{Sakai}
S. Sakai, \emph{A characterization of $W^*$-algebras}, Pacific J. Math. \textbf{6} (1956), 763--773.

\bibitem[Sak2]{Sakaib}
\bysame, \emph{$C^*$-algebras and $W^*$-algebras}, Springer-Verlag, Berlin, 1971.  

\bibitem[SW]{SW}
R. Smith and J. Ward, \emph{$M$-ideal structure in Banach algebras}, J. Funct. Anal. \textbf{27} 
(1978), 337--349.

\bibitem[Weg]{Wegge}
N. E. Wegge-Olsen, \emph{$K$-theory and $C^*$-algebras}, Oxford University Press, Oxford, 1993.

\bibitem[KWer]{kw}  
K. H. Werner, {\em A characterization of $C^*$-algebras by nh-projections on matrix ordered spaces,} 
unpublished preprint, Saarbrucken.

\bibitem[WWer]{Wend}
W. Werner, \emph{Multipliers on matrix ordered operator spaces and some $K$-groups}, preprint (2002).

\bibitem[Zar1]{Zarikian}
V. Zarikian, \emph{Complete one-sided $M$-ideals in operator spaces}, Ph.D. thesis, UCLA, 2001.

\bibitem[Zar2]{ZarikianMpaper}
\bysame, \emph{Local characterizations of one-sided $M$-ideals (working title)}, in preparation.

\end{thebibliography}
\end{document}